\documentclass[twoside, 10pt]{article}

\usepackage[hmargin= 0.75in, height=9.4in]{geometry}
\usepackage{amssymb,amsthm, amsfonts,amsmath,mathrsfs, amsxtra, url, hyperref}

\usepackage{fancyhdr}
\fancypagestyle{plain}{
 \fancyhf{}

}

\renewcommand{\theequation}{\thesection.\arabic{equation}}
 \numberwithin{equation}{section}

\newtheorem {thm}{Theorem}[section]
\newtheorem {prop}{Proposition}[section]
\newtheorem {lemm}{Lemma}[section]
\newtheorem {deff}{Definition}[section]
\newtheorem {cor}{Corollary}[section]
\newtheorem {rem}{Remark}[section]

\makeatletter
\newenvironment{highlightequation}{%
  \def\tagform@##1{\maketag@@@{(\ignorespaces##1\unskip\@@italiccorr*)}}%
  \ignorespaces
}{%
  \def\tagform@##1{\maketag@@@{(\ignorespaces##1\unskip\@@italiccorr)}}%
  \ignorespacesafterend
}
\makeatother

\def\ba{\begin{array}}
\def\ea{\end{array}}
\def\bea{\begin{eqnarray}}
\def\eea{\end{eqnarray}}
\def\beas{\begin{eqnarray*}}
\def\eeas{\end{eqnarray*}}
\def\bi{\begin{itemize}}
\def\ei{\end{itemize}}
\def\bc{\begin{cases}}
\def\ec{\end{cases}}

\def\bhe{\begin{highlightequation}  }
\def\ehe{\end{highlightequation}  }

\def\a{\alpha}
\def\ga{\gamma}
\def\d{\delta}
\def\e{\varepsilon}
\def\z{\zeta}
\def\k{\kappa}
\def\l{\lambda}

\def\si{\sigma}

\def\o{\omega}

\def\D{\Delta}
\def\G{\Gamma}
\def\L{\Lambda}
\def\O{\Omega}

\def\Th{\Theta}
\def\U{\Upsilon}

\def\bF{{\bf F}}

\def\cD{{\cal D}}
\def\cE{{\cal E}}
\def\cF{{\cal F}}
\def\cG{{\cal G}}

\def\cI{{\cal I}}

\def\cK{{\cal K}}

\def\cN{{\cal N}}

\def\cR{{\cal R}}

\def\cT{{\cal T}}

\def\cX{{\cal X}}
\def\cY{{\cal Y}}
\def\cZ{{\cal Z}}

\def\hE{\mathbb{E}}

\def\hH{\mathbb{H}}

\def\hK{\mathbb{K}}
\def\hL{\mathbb{L}}

\def\hN{\mathbb{N}}

\def\hP{\mathbb{P}}
\def\hQ{\mathbb{Q}}
\def\hR{\mathbb{R}}
\def\hS{\mathbb{S}}

\def\hV{\mathbb{V}}

\def\sB{\mathscr{B}}

\def\sJ{\mathscr{J}}
\def\sK{\mathscr{K}}

\def\sN{\mathscr{N}}

\def\sP{\mathscr{P}}

\def\sY{\mathscr{Y}}
\def\sZ{\mathscr{Z}}

\def\fL{\mathfrak{L}}

\def\fh{\mathfrak{h}}
\def\fg{\mathfrak{g}}

\def\fc{\mathfrak{c}}

\def\({\left(}
\def\){\right)}
\def\[{\left[}
\def\]{\right]}
\def\lan{\langle}
\def\ran{\rangle}

\def\no{\noindent}

\def\ss{\smallskip}
\def\ms{\medskip}

\def\q{\quad}
\def\qq{\qquad}
\def\n{\negthinspace}
\def\dn{\n  \n }
\def\tn{\n  \n  \n }

\def\ol{\overline}
\def\ul{\underline}
\def\ua{\mathop{\uparrow}}
\def\da{\mathop{\downarrow}}

\def\wt{\widetilde}
\def\wh{\widehat}
\def\fra{\mathfrak{a}}

\def\dtp{{\hbox{$dt \otimes d \hP-$a.s.}}}
\def\dsp{{\hbox{$ds \otimes d \hP-$a.s.}}}
\def\pas{{\hbox{$\hP-$a.s.}}}

\def\hb{\hbox}
\def\dis{\displaystyle}
\def\cd{\cdot}
\def\cds{\cdots}
\def\sgn{\hb{sgn}}

\def\fa{\,\forall \,}

\def\es{\emptyset}

\def\b1{{\bf 1}}
\def\qed{\hfill $\Box$ \medskip}

\def\essinf{\mathop{\rm essinf}}
\def\esssup{\mathop{\rm esssup}}
\def\liminf{\mathop{\ul{\rm lim}}}
\def\limsup{\mathop{\ol{\rm lim}}}
\newcommand{\esup}[1]{ \underset{#1}{\esssup}\,}
\newcommand{\einf}[1]{ \underset{#1}{\essinf}\,}

\newcommand{\lmt}[1]{ \underset{#1}{\lim}}
\newcommand{\lmtu}[1]{ \underset{#1}{\lim} \n  \ua \,}
\newcommand{\lmtd}[1]{ \underset{#1}{\lim} \n  \da \,}

\begin{document}

 \title{\bf  Doubly   Reflected   BSDEs  with Integrable \\ Parameters and Related Dynkin Games
 \thanks{We would like to thanks the referees for their careful reading and helpful comments which helped us improve our paper.}}

\author{
  Erhan Bayraktar\thanks{ \noindent Department of
  Mathematics, University of Michigan, Ann Arbor, MI 48109; email:
{\tt erhan@umich.edu}.}  \thanks{E. Bayraktar is supported in part by the National Science Foundation  a Career grant DMS-0955463 and an Applied Mathematics Research grant DMS-1118673, and in part by the Susan M. Smith Professorship. Any opinions, findings, and conclusions or recommendations expressed in this material are
those of the authors and do not necessarily reflect the views of the National Science Foundation.} $\,\,$,
$~~$Song Yao\thanks{
\noindent Department of
  Mathematics, University of Pittsburgh, Pittsburgh, PA 15260; email: {\tt songyao@pitt.edu}. } }
\date{}

\maketitle

 \begin{abstract}

   We study   a doubly reflected backward stochastic differential  equation (BSDE)
 with integrable parameters and the related  Dynkin game.
 When the lower obstacle $L$ and the upper obstacle $U$ of the equation  are completely separated,
 we construct a   unique solution of the doubly reflected BSDE by pasting local solutions,
 and   show that  the $Y-$component of the unique solution represents the value process of the corresponding
 Dynkin game under   $g-$evaluation,
   a nonlinear expectation induced by BSDEs  with the same generator $g$ as the doubly reflected BSDE concerned.
 In particular, the first time  $\tau^*$  when process $Y $  meets $L$
 and the first time  $\ga^*$  when process $Y $  meets $U$ form a saddle point of the Dynkin game.


\end{abstract}

 \smallskip


  \section{Introduction}

  \label{sec:introduction}

 In this paper, we study a  doubly reflected backward stochastic differential equation
 with generator $g$, integrable  terminal data $\xi$ and two
   integrable obstacles $L$,  $ U $
       \bea   \label{DRBSDE}
    \left\{\ba{l}
 \dis  Y_t =   \xi  + \int_t^T g \( s,Y_s, Z_s \)  ds    +    K_T \n -\n  K_t
   -     J_T   +   J_t
   - \int_t^T Z_s d B_s    , \q    t \in [0,T] ,   \vspace{1mm}      \\
   \dis L_t   \le Y_t  \le  U_t        , \q    t \in [0,T] , \vspace{1mm} \q \\
  \dis       \int_0^T \n    (  Y_t   -  L_t   )     d K_t
  = \int_0^T \n    (  U_t - Y_t   )     d J_t   = 0       \q (\hb{flat-off conditions}) .
     \ea \right.
    \eea
    A solution of such an equation consists of four adapted processes:  a  continuous process $Y$,
    a locally square-integrable process $Z$ and two continuous increasing processes $K$ and $J$.
     Klimsiak \cite{Klimsiak_2013} studied the same problem  but assumed
     an extended {\it Mokobodzki's condition}:
   there exists a semi-martingale between $L$ and $U$, which is   practically difficult to verify.
   Instead, we only require   the two obstacles $L$, $U$ to be completely separable, i.e. $L_t<U_t$, $\fa t \in [0,T]$.

 Backward stochastic differential equations
 (BSDEs) were introduced in linear case by Bismut \cite{Bismut-73} as the adjoint equations for the stochastic
Pontryagin maximum principle in control theory. Later, Pardoux and Peng \cite{PP-90} extended them to a fully
nonlinear version
 \bea  \label{BSDE}
       Y_t=   \xi + \int_t^T g(s,Y_s, Z_s   )  ds
   -\int_t^T Z_s d B_s    , \q    t \in [0,T]  ,
    \eea
and showed that  the BSDE admits a unique solution $(Y,Z)$ when
   generator $g$ is  Lipschitz continuous in $(y,z)$
 and   terminal datum $\xi$ is square-integrable.
   Since then, the theory of   BSDEs has rapidly grown    and   been applied in
many areas   such as mathematical finance, theoretical economics, stochastic   control, stochastic
differential games, partial differential equations
(see e.g. the references in \cite{EPQ-97} or in \cite{Karatzas_Soner_1998}).

As a variation of BSDEs, a  BSDE   with one reflecting obstacle (say lower obstacle $L$)
 \bea \label{RBSDEL}
   \begin{cases}
   \dis      L_t \le  Y_t=   \xi + \int_t^T g(s,Y_s, Z_s   )  ds + K_T - K_t
   -\int_t^T Z_s d B_s    , \q    t \in [0,T]  ,          \vspace{1mm}     \\
   \dis \int_0^T  (Y_t-L_t) dK_t = 0   \q (\hb{flat-off condition}) .
   \end{cases}
    \eea
 was first studied by El Karoui et al. \cite{EKPPQ-1997}. If   $g$ is  Lipschitz continuous in $(y,z)$
 and if both  terminal datum $\xi$ and  lower obstacle $L$ are square-integrable,
 these authors showed that the reflected BSDE
 has a unique solution $(Y,Z,K)$ and that  the $Y-$component of the unique solution is the Snell envelope
 of the reward process $L$ in the related optimal stopping problem under  $g-$evaluation
(for a more general statement, see e.g.  appendix A of \cite{CTZ_2013},  Section 7 of \cite{QRBSDE}).
 As a nonlinear expectation induced by BSDEs  with the same generator $g$ as the   reflected BSDE,
 the $g-$evaluation  possesses many (martingale) properties of the classic linear expectation
 and   thus become a very useful tool in nonlinear analysis. In particular, the $g-$evaluation
 is closely related to risk measures in mathematical finance.

 Based on \cite{EKPPQ-1997}, Cvitani\'c and Karatzas \cite{Cvitanic_Karatzas_1996}
 extended the   research  of BSDEs to those with two reflecting obstacles.
 They showed that   a doubly reflected BSDE with Lipschitz generator, square-integrable terminal datum
 and square-integrable obstacles
 admits a unique solution under   Mokobodzki's condition (there exists a quasimartingale between two obstacles)
  or   certain regularity
 condition on one of the obstacles \big(see assumption (H) of \cite{HLM_1997} for a simplified form\big).
 Cvitani\'c and Karatzas also found  that   the $Y-$component of the unique solution
 is  exactly  the value process of the related Dynkin game,
 a zero-sum stochastic differential game of optimal stopping, under $g-$evaluation
  (for a more general statement, see e.g. \cite{DQS_2013}).
 From a perspective of mathematical finance,  this discovery
 is   significant  for the evaluation of American game options or {\it Israeli} options, see e.g.
 Hamad\`ene \cite{Hamadene_2006}. Later,
   Hamad\`ene et al. \cite{Hamadene_Lepeltier_Wu_1999, Hamadene_Lepeltier_2000, Hamadene_2006}
  added controls into  a doubly reflected BSDE  and the drift coefficient  of the associated  state  process   to
  analyze a {\it mixed}   zero-sum controller and stopper  game
    as well as the corresponding {\it saddle} point problem.
  For the literature
  and the recent advances of Dynkin games, see e.g. \cite{Karatzas_Wang_2001,
  Touzi_Vieille_2002, Hamadene_Zhang_2009, MR3162260}.
  As to the history
  and latest development of   controller and stopper games, see e.g.
  \cite{MR1872738, MR2435857, OSNE1, OSNE2, OS_CRM, BH11,  ETZ_2012, NZ_2012, ROSVU}.

    Among other development in   doubly reflected BSDEs,
     Lepeltier and San Mart\'in \cite{LS_2004} obtained  the existence result when
   $g$ is only continuous and has linear growth in variables $(y,z)$;
      Xu  \cite{Xu_2007} got the wellposedness result when
  the Lipschitz continuity of    $g$
  in $y-$variable  is relaxed to a monotonicity condition; 
    and Bahlali et al. \cite{Bahlali_Hamadene_Mezerdi_2005},
    Essaky et al. \cite{EHO_2011, Essaky_Hassani_2013} analyzed the existence of a maximal solution   when
     $g$ has quadratic growth in $z-$variable.

  All the above articles on doubly reflected BSDEs, except \cite{Hamadene_2006},
  assumed either (extended) Mokobodzki's condition
   or the aforementioned   regularity  condition.
  According to \cite{Hamadene_2006}'s observation that the existence of   local solutions of
  a doubly reflected BSDE   relies on neither of these two conditions,
 Hamad\`ene and Hassani \cite{Hamadene_Hassani_2005} pasted local solutions to form a unique solution
 of a doubly reflected BSDE with two distinct obstacles. Since then, the complete separation of   obstacles
   has  been   postulated by most of the subsequent papers
 including \cite{Buckdahn_Li_3, Lp_DRBSDE, Hamadene_Hdhiri_2006, HRZ_2009}
 as well as the present one.

       During the evolution of the BSDE theory, some efforts were made  to weaken   the square integrability
  on  terminal data so as to match up with the fact that linear BSDEs are well-posed for integrable terminal data:
   El Karoui et al. \cite{EPQ-97} demonstrated that for any $p-$integrable terminal datum with $p \in (1, \infty)$,
   a BSDE with Lipschitz generator  admits a unique $p-$integrable solution. This wellposedness  result was later upgraded
  by Briand et al. \cite{Briand_Carmona_2000, BH_Lp_2003} who reduced the  Lipschitz condition
  of generator $g$ on  $y-$variable   to a   monotonicity condition on $y$.
  After Hamad\`ene and Popier \cite{Lp_RBSDE}   extended \cite{BH_Lp_2003}'s results
  for reflected BSDEs,  Hamad\`ene et al. \cite{Lp_DRBSDE} make a
  further generalization for doubly reflected BSDEs with two completely separate obstacles.

  We   dedicate this paper to the solvability of the doubly reflected BSDE \eqref{DRBSDE}
  with integrable parameters   and will discuss the related Dynkin game.
  Besides the   monotonicity condition on $y-$variable and the Lipschitz condition on $z-$variable,
  if the generator $g$  additionally has a growth condition on $z-$variable
  of order $\a \in (0,1)$ \big(see (H7) of \cite{BH_Lp_2003} or (H5) in the current paper\big),
  then the BSDE with integrable terminal datum admits a unique  solution
  $(Y,Z)$ such  that   both $Y$ and $Z$ are $p-$integrable processes  for any $p \in (0,1)$
  and  that $Y $  is of class (D).
   So the corresponding $g-$evaluation is well-defined
   for each integrable random variable.
  Under the same hypotheses on generator $g$ as Section 6 of \cite{BH_Lp_2003},
  we will demonstrate  a similar wellposedness result  for  doubly  reflected BSDEs
  with integrable terminal data and two distinct integrable obstacles. Though we   follow
  the approach  of \cite{Hamadene_Hassani_2005,   Lp_DRBSDE} on pasting local solutions, the estimations
   used for $L^p-$solutions, $p>1$ are no longer valid in the $p=1$ or class (D) case. We managed to derive
   some novel estimation and approximation scheme.

 To construct  a unique solution  of a reflected BSDE  with integrable terminal datum $\xi$ and
 integrable lower obstacle $L$,
we   use  the {\it penalization} method introduced in \cite{EKPPQ-1997}
     together with a {\it localization} technique. This is because   the approximating solutions are
    only $p-$integrable ($\fa p  \n \in \n  (0,1)$): Given $n \in \hN$,
     we compensate  the generator $g$ by $n$ times the distance that $y-$variable is below   $L_t$,
     i.e.  $g_n (t,y,z) : = g (t,y,z) + n (y-L_t)^-$.  The BSDE with generator $g_n$ and
     terminal datum $\xi$
     \bea  \label{eq:b711}
     Y^n_t =   \xi  + \int_t^T g \( s,Y^n_s, Z^n_s \)  ds    +   n \int_t^T  \big(  Y^n_s -L_s  \big)^-  ds
   - \int_t^T Z^n_s d B_s    , \q    t \in [0,T] ,
     \eea
      has a unique $p-$integrable ($\fa p  \n \in \n  (0,1)$) solution $(Y^n,Z^n)$  such that   $Y^n$ is of class (D).
       The monotonicity of $\{g_n\}_{n \in \hN }$   implies that of $\{Y^n\}_{n \in \hN}$, thanks to
       a  general comparison result (Proposition \ref{prop_BSDE_comp_basic}). 
       Then we can find a stopping time $\tau_\ell$ such that $|Y^n|$ is uniformly bounded by $\ell$ over the
       stochastic interval $[ \n [ 0, \tau_\ell ] \n ]$.
       By a local estimation (Lemma \ref{lem_RBSDE_estimate}),
       the  local $\hL^2-$norms of $Z^n$'s are uniformly bounded by a multiple of $  \ell^2$. So up to a subsequence,
       $ Z^n $ weakly converges to some $\cZ^\ell$. Consequently, we can deduce that
       $K^n_t  \n := \n   n \int_0^t  \n   \big(  Y^n_s \n - \n L_s  \big)^-  ds $ converges
       to $\cK^\ell_t \n := \n  Y_0  \n - \n  Y_t  \n - \n  \int_0^t  \n  g(s, Y_s, \cZ^\ell_s) ds
        \n + \n  \int_0^t  \n  \cZ^\ell_s d B_s $ uniformly
        over $[ \n [ 0, \tau_\ell ] \n ]$. Letting $n \to \infty$ in \eqref{eq:b711} shows that
        $(Y, \cZ^\ell, \cK^\ell)$ is a local  solution of \eqref{RBSDEL} over $[ \n [ 0, \tau_\ell ] \n ]$.
  Pasting up   $ (Y, \cZ^\ell, \cK^\ell)  $'s over stochastic intervals $ ]\n] \tau_{\ell-1} ,  \tau_\ell  ]\n] $'s
  we obtain  a global $p-$integrable ($\fa p \n \in \n  (0,1)$) solution $(Y,Z,K)$ of    \eqref{RBSDEL}.
  The uniqueness of such a solution follows from a comparison result (Proposition \ref{prop_RBSDE_comp})
  of reflected BSDEs,   which is a corollary of Proposition \ref{prop_BSDE_comp_basic}.

 Applying Proposition \ref{prop_BSDE_comp_basic} again shows that with respect to the corresponding $g-$evaluation,
  the $Y-$component of the unique solution of
 \eqref{RBSDEL} is a supermartingale  and even a martingale up to the first time when process $Y$ meets
 the lower obstacle  $L$. Consequently,   $Y$ is the Snell envelope
 of the reward process $L$ in the related
 optimal stopping problem in which the player   is  trying to select a
 best exit time from the game so as to maximize her expected reward under $g-$expectation.

 Based on the wellposedness result   for reflected BSDEs with integrable parameters,
 we next   take   \cite{Hamadene_Hassani_2005}'s approach of pasting local solutions
 to construct a global solution of \eqref{DRBSDE}: Let $(Y^n,Z^n,K^n)$ be the
 unique $p-$integrable ($\fa p  \n \in \n  (0,1)$) solution of a reflected BSDE
 with the penalized generator $g_n$  and the upper obstacle $U$.
   We first show that the increasing limit $Y$ of $Y^n$'s, together with some processes $(Z^\ell,K^\ell)$,
     solves \eqref{RBSDEL} over some stochastic intervals $[\n[ \nu_\ell, \nu'_\ell ]\n]$ for any $\ell \in \hN$.
 A reverse conclusion can be obtained    for the limit $\wt{Y}$ of a decreasing scheme that involves
 reflected BSDEs with generator $\wt{g}_n (t,y,z) \n := \n g(t,y,z) \n - \n n (y-U_t)^+ $ and the lower obstacle $L$:
 For some processes $(\wt{Z}^\ell,\wt{J}^\ell)$,
 $(\wt{Y},\wt{Z}^\ell,\wt{J}^\ell)$ solves a reflected BSDE with upper obstacle $U$ over some stochastic interval
 $[\n[ \nu'_\ell, \nu_{\ell+1} ]\n]$ for any $\ell \in \hN$.  Then pasting
   $(Y,Z^\ell,K^\ell,0) $ and $(\wt{Y},\wt{Z}^\ell,0, \wt{J}^\ell) $
   alternatively over $[\n[ \nu_\ell, \nu'_\ell ]\n]$ and $[\n[ \nu'_\ell, \nu_{\ell+1} ]\n]$
  yields a global $p-$integrable ($\fa p  \n \in \n  (0,1)$) solution of the doubly reflected BSDE \eqref{DRBSDE}.

 Leveraging Proposition \ref{prop_BSDE_comp_basic} once again shows that with respect to the corresponding $g-$evaluation,
  the $Y-$component of the   solution of \eqref{DRBSDE} just constructed  is
   a  submartingale up to the first time $ \tau^*   $
   when $Y$ meets the lower obstacle $L$,  and   is a  supermartingale up to time $ \ga^*   $
   when $Y$ meets the upper obstacle $U$.
  Consequently,   $Y$ is the value process of the related Dynkin game under $g-$evaluation
  in which $L$ \big(resp. $U$\big) is the amount process a player will receive from  her opponent
   when she stops the game earlier \big(resp.   not earlier\big) than  her opponent.
   The uniqueness result of \eqref{DRBSDE} then easily follows.
  Moreover,  the pair $(\tau^*, \ga^*)$ forms a saddle point of such a Dynkin game.

 Since  dealing mostly   with $p-$integrable ($\fa p \n \in \n  (0,1)$) solutions,
 we can not apply   Doob's martingale inequality and many well-known estimates in BSDE theory
 without using localization first, which  increases the technical difficulty.
 Also, to overcome   technical subtleties we encounter when proving the
 $p-$integrability  ($\fa p \n \in \n  (0,1)$) of the limit $Y$ in the penalization scheme,
 we   appropriately exploit   Tanaka-I\^to's formula, Hypothesis  (H5) and other tricks,
 see in particular the proof of \eqref{eq:a073}.


    The rest of the paper is organized as follows:
     After listing necessary notations, we give the definition of doubly reflected BSDEs and make some
     assumptions on their generators $g$ in Section \ref{sec:introduction}. We   first present in Section
    \ref{sec:DRBSDE} the main result of our paper,
    a wellposedness result of doubly reflected BSDEs with integrable parameters  as well as
    the $g-$martingale characterization of the $Y -$component of the unique solution, the latter of
    which implies that $Y$ is a value process of the related Dynkin games under $g-$evaluation.
    Section \ref{sec:BSDE} recalls   a wellposedness result   of   BSDEs with integrable terminal data
    and gives a general comparison result     for BSDEs over stochastic intervals,
    which plays an important role in our analysis. The unique solutions of BSDEs with generator $g$
    and integrable terminal data induce a widely-defined nonlinear expectation, called ``$g-$evaluation/expectation",
    whose properties will be discussed in Section \ref{sec:g_evaluation}.
    In Section \ref{sec:RBSDE}, to construct a unique solution for a   reflected  BSDE  with integrable parameters
    as a preparation for our main result,   we   use the {\it penalization} method which involves
    two auxiliary monotonicity results. And we show that the $Y-$component of the unique solution
    of the reflected BSDE is exactly the Snell envelope in the related optimal stopping problem under  $g-$evaluation.
   Section \ref{sec:proofs} contains   proofs of our results while the demonstration of
   some technical claims   are deferred to   the  Appendix.

\subsection{Notation and Definitions} \label{subsec:notation}


    \ms  Throughout this paper,  we fix a   time horizon $T \in (0, \infty)$,
    and let $B$ be a $d-$dimensional standard Brownian Motion defined
    on a complete probability space $(\O,\cF, \hP)$.
    The augmented filtration generated by $B$
 \beas
 \bF= \left\{\cF_t := \si \( \si\(B_s; s\in [0,t]\) \cup \sN \) \right\}_{t  \in  [0, T] }
 \eeas
  satisfies the {\it usual hypothesis},  where $\sN$   collects all $\hP-$null sets in $\cF$.

  Let $\cT$ be the set of all $\bF-$stopping times $\tau$ taking values in  $[0 , T]$.
   For any $\nu, \tau \n \in \n  \cT$ with $\nu  \n \le \n  \tau$,
   we  set $\cT_{\nu,\tau}  \n := \n  \{ \ga  \n \in \n  \cT \n :
      \nu   \n  \le \n  \ga    \n  \le \n  \tau    \}$.
    An increasing sequence $\{\tau_n\}_{n \in \hN}$ in $\cT$ is called ``stationary" if
  for \pas ~ $\o \in \O$, $T \n = \n \tau_n (\o)$ for some  $ n  \n = \n  n(\o)  \n \in \n  \hN $.
  As usual,   we say that a  $\sB \big([0,T]\big) \otimes \cF-$measurable
  process $X$ is of {\it class $($D}), with respect to $(\cT,\hP)$,  if
  $  \{ X_\tau  \}_{\tau \in \cT}$ is $\hP-$uniformly integrable.
   Moreover, we let $\sP$ denote the $\bF-$progressively measurable $\si-$field on $ [0,T] \times \O$ and
   will  use the convention  $ \inf \es := \infty$.


  \ss  Let $p \in (0,\infty)$. It holds  for any finite subset  $\{a_1, \cds, a_n\} $ of $(0, \infty)$ that
   \bea   \label{eqn-d011}
  \( 1 \land n^{p-1} \) \sum_{i=1}^n a_i^p \le
\(\sum_{i=1}^n a_i \)^p \le \( 1 \vee n^{p-1} \) \sum_{i=1}^n a_i^p .
 \eea
 And for any $p' \in (p,\infty)$, one has
 \bea  \label{eqn-d011b}
 x^p \le 1+ x^{p'} , \q \fa x \in (0,\infty) .
 \eea

 The following  spaces  will be frequently used in the sequel.

 \ss \no 1)\,  For any sub$-\si-$field  $\cG$  of $\cF$, let  $L^0 (\cG )$  be  the space of all real-valued,
$\,\cG-$measurable random variables $\xi$ and set $L^p (\cG )  \n : = \n   \Big\{\xi  \n \in \n  L^0 (\cG )  \n : \,
 \|\xi\|_{ L^p(\cG )}  \n := \n
 \big\{ \hE  [  |\xi |^p  ] \big\}^{1 \land \frac{1}{p}}<\infty \Big\}$.

  \ms \no 2)\,     We need the following subspaces of $\hS^0$, which denotes all real-valued, $\bF-$adapted continuous processes:

 \no  $\bullet$ $\hS^p   \n :=  \n    \Big\{  X  \n \in \n  \hS^0  \n :\,   \|X\|_{\hS^p}
   \n := \n  \big\{ \hE \[ (X_*)^p  \]  \big\}^{1 \land \frac{1}{p}} \n < \n \infty \Big\}$,
    where $X_*  \n := \n  \underset{t \in [0,T]}{\sup}|X_t| $;

   \no    \ss $\bullet$ $\hS^p_+   \n := \n     \Big\{  X  \n \in \n  \hS^0  \n :\, X^+
   \n = \n X \vee 0  \n \in \n  \hS^p  \Big\}$ and
 $\hS^p_-   \n := \n     \Big\{  X  \n \in \n  \hS^0  \n :\, X^-
 \n = \n (-X) \vee 0  \n \in \n  \hS^p  \Big\}$;


  \no    \ss  $\bullet$ $\hV^0 := \big\{  X \in \hS^0 :\,  X$ is of  finite variation\big\};

  \no    \ss   $\bullet$  $  \hK^0  :=  \big\{  X \in \hS^0 :\,  X$
   is an increasing process with $X_0=0 \big\}  \subset \hV^0$;

  \no    \ss   $\bullet$  $  \hK^p := \big\{ X \in \hK^0 :\,   X_T \in \hL^p (\cF_T)     $\big\}.

\ss \no 3)\,   Let     $ \wt{\hH}^{2,0} $ (resp. $\hH^{2,0} $) denote the space of all $\hR^d-$valued,
 $\bF-$progressively measurable (resp. $\bF-$predictable) processes $X$ with  $\int_0^T |X_t|^2 dt < \infty$, \pas ~
 and set   $   \hH^{2,p}  : =  \left\{ X \in   \hH^{2,0}  :
   \|X\|_{\hH^{2,p}} := \left\{ \hE \[ \big(  \int_0^T \n  |X_t|^2     dt
   \big)^{p/2} \] \right\}^{1 \land (1/p)}<\infty \right\} $.

 \ss \no   In the above notations, if $p \n  \ge \n   1$, $\|\cd \|_{\Xi^p}$ is a norm on
 $\Xi^p  \n  = \n   L^p(\cG), \hS^p, \hH^{2,p}  $.
 And if $p  \n  \in \n   (0,1)$, $(X,X') \to \|X  \n  - \n   X'\|_{\Xi^p}$ defines a distance on $\Xi^p$,
 under which  $\Xi^p$ is a complete metric space.

 \ss  Let us recall the notions of backward stochastic differential equations (BSDEs), reflected BSDEs
 and doubly reflected BSDEs: A (basic) parameter pair $(\xi,g)$ consists of a real-valued, $\cF_T-$measurable
 random variable $ \xi  $ and a function $g  \n   :  [0,T] \n  \times \n   \O  \n  \times \n   \hR  \n  \times \n   \hR^d
 \to  \hR $ that is  $\sP  \n  \otimes \n   \sB(\hR)  \n  \otimes \n     \sB(\hR^d)/\sB(\hR) -$measurable.

    \begin{deff}

    Given a parameter pair $(\xi,g)$, let  $ L, U \in \hS^0 $ such that $\hP\{L_t \le U_t, \fa t \in [0,T] \} = 1$
    and $L_T \le \xi \le U_T$, \pas ~   We say that
  1\big) $(Y,Z) \in \hS^0 \times \wt{\hH}^{2,0} $ is
    a solution of  a  BSDE with   terminal data $\xi$ and generator $g$ \big(BSDE\,$(\xi,g)$ for short\big)
    if \eqref{BSDE} holds \pas ~
  2\big) A triplet $(Y,Z,K) \in \hS^0 \times \wt{\hH}^{2,0} \times \hK^0 $ is a solution of  a reflected BSDE
    with   terminal data $\xi$,  generator $g$ and \big(lower\big) obstacle $L$ \big(RBSDE\,$(\xi,g,L)$ for short\big)
    if \eqref{RBSDEL} holds \pas ~
  3\big)  A quadruplet     $(Y,Z,K,J)
\in \hS^0  \times \wt{\hH}^{2,0}  \times \hK^0  \times \hK^0  $ is a solution of a  doubly reflected  BSDE with
terminal data $\xi$, generator $g$,  lower obstacle $L$ and upper obstacle $U$
\big(DRBSDE\,$(\xi,g,L,U)$ for short\big) if \eqref{DRBSDE} holds \pas

\end{deff}

\begin{rem}  \label{rem_RBSDEU}
  Given a  parameter pair $(\xi,g)$,
  \bea  \label{eq:g_neg}
  g_-(t,\o,y,z) := - g(t,\o,-y,-z) , \q \fa (t,\o,y,z) \in [0,T] \times \O \times \hR \times \hR^d
  \eea
   clearly defines a  $\sP     \otimes     \sB(\hR)     \otimes     \sB(\hR^d)/\sB(\hR) -$measurable function.
   For any  $ L \in \hS^0 $ with  $L_T \le \xi  $, \pas,
  $(Y,Z,K) \in \hS^0 \times \wt{\hH}^{2,0} \times \hK^0 $ solves RBSDE$(\xi,g,L)$ if and only if
  $ (\wt{Y}, \wt{Z}, \wt{J} ) = (-Y,-Z,K) \in \hS^0 \times \wt{\hH}^{2,0} \times \hK^0  $
  is a solution of  the following reflected BSDE
  with   terminal data $\wt{\xi}=-\xi$,  generator $g_-$ and  upper obstacle $U=-L$:
 \bea \label{RBSDEU}
   \begin{cases}
   \dis      U_t \ge  \wt{Y}_t=   \wt{\xi} + \int_t^T g_-( s,\wt{Y}_s, \wt{Z}_s   )  ds - \wt{J}_T + \wt{J}_t
   -\int_t^T \wt{Z}_s d B_s    , \q    t \in [0,T]  ,         \vspace{1mm}     \\
   \dis \int_0^T  (U_t - \wt{Y}_t   ) d \wt{J}_t = 0 . \q (\hb{flat-off condition})
   \end{cases}
    \eea

\end{rem}

 Let  $g  : [0,T] \times \O \times \hR \times \hR^d \to \hR$ be a
   $\sP     \otimes     \sB(\hR)     \otimes     \sB(\hR^d)/\sB(\hR) -$measurable function.
 To study  doubly  reflected BSDEs with generator $g$ and integrable parameters  $(\xi,L,U)$,
 we will make the following assumptions on function   $g$:

\ss \no {\bf Standing assumptions on $g$.}

 \no  Let $ \k >0 $, $\l \in \hR$, $\a \in (0,1)$ and let $\{ h_t \}_{t \in [0,T]} $ be a non-negative
 integrable process \big(i.e. $h \in \hL^1 ([0,T] \times \O, \sB([0,T]) \otimes \cF, dt \otimes \hP ) $\big).
 It holds  \dtp~ that

 \ss \no {\bf (H1)}   $
   |g(t,\o,y,z )-g(t,\o,y,z')| \le \k |z -z'| $, $ \fa y \in \hR $, $ \fa z ,z' \in \hR^d $;

    \ss \no{\bf (H2)} $
  \sgn( y -y' ) \cd \( g(t,\o,y ,z)- g(t,\o,y',z) \)  \le \l |y -y'|  $, $ \fa y , y' \in \hR$,
  $ \fa z \in \hR^d $;

  \ss \no{\bf (H3)}   $y \to g(t,\o,y,z) $ is continuous,  $ \fa  z   \in   \hR^d $;

   \ss \no{\bf (H4)} $|g(t,\o,y,0)| \le h_t (\o) + \k |y|$,  $ \fa  y \in \hR$;

    \ss \no{\bf (H5)}
  $  |g(t,\o,y,z) - g(t,\o,y,0)| \le \k (h_t(\o)+|y|+|z|)^\a, \q \fa (y,z) \in \hR \times \hR^d $.

  \ss From now on,  for any  $p \n \in \n  [ 0, \infty) $ we let $C_p$ be a generic constant
  depending on $p, \k, \l^+,   T      $ and $ \hE \n \int_0^T \n  h_t dt$ \big(in particular,
  $C_0$ will denote a generic constant depending   on $ \k, \l^+,   T   $ and $ \hE \n \int_0^T \n  h_t dt$\big),
  whose form may vary from line to line. For convenience, we will call
  a function $g : [0,T] \times \O \times \hR \times \hR^d \to \hR$ a ``generator"
  if it is $\sP        \otimes        \sB(\hR)        \otimes        \sB(\hR^d)/\sB(\hR) -$ measurable
  and satisfies      $($H1$)$$-$$($H5$)$.

  \begin{rem} \label{rem_assum0}
   If a function $g : [0,T] \times \O \times \hR \times \hR^d \to \hR$ is   Lipschitz continuous in $y$
\big(i.e.  for some $\wt{\k} > 0$, it holds \dtp ~ that $|g(t,\o,y,z )-g(t,\o,y',z)|
\le \wt{\k} | y - y'| $, $ \fa y, y' \in \hR $, $ \fa z   \in \hR^d $\big),
 then  \big(H2\big) automatically holds and \big(H4\big) will be replaced by
 $  |g(t,\o,0,0)| \le h_t(\o)  $, \dtp
 \end{rem}

   \begin{rem} \label{rem_assum}
   Let $g$ be a generator.

\ss \no 1$)$ The function $g_-$ defined in \eqref{eq:g_neg} is also a  generator.

\ss \no 2$)$  Given $\tau \n \in \n  \cT$, since $\big\{\b1_{\{t \le \tau\}} \big\}_{t \in [0,T]}$
 is an $\bF-$adapted c\`agl\`ad process
 \big(and thus $\bF-$predictable\big),  the measurability of $g$ implies that
  \bea  \label{def_g_tau}
  g_\tau (t,\o,y,z) \n := \n  \b1_{\{t \le \tau (\o)\}} g (t,\o,y,z)  , \q
   \fa (t,\o,y,z)  \n \in \n  [0,T]  \n \times \n  \O  \n \times \n  \hR  \n \times \n  \hR^d
   \eea
  defines a $\sP    \n   \otimes    \n   \sB(\hR)    \n   \otimes   \n    \sB(\hR^d) / \sB(\hR) -$measurable function.
   And one can deduce that $ g_\tau $ also satisfies $($H1$)$$-$$($H5$)$
   \big(actually, it satisfies (H2) with $\wt{\l} \n = \n \l  \n \vee \n  0$\big).

 \ss \no 3$)$  If $  g' $ is another generator, so  is  $ ag + bg'$ for any $a,b > 0$.

 \ss \no 4$)$  Given  $L \in \hS^1_+$,
 $ g_L (t,\o,y) := \(y-L_t(\o)\)^- $, $ ( t,\o, y ) \in [0,T] \times \O \times \hR $
 is clearly a  $\sP  \n  \otimes \n   \sB(\hR) / \sB(\hR) -$measurable function that is
  Lipschitz continuous in $y$ 
 and satisfies  $ \hE \int_0^T  g_L (t,0) dt = \hE \int_0^T L^+_t dt
 \le T \|L^+\|_{\hS^1}  < \infty $. By Remark \ref{rem_assum0},   $g_L$ satisfies $($H2$)$$-$$($H4$)$.
  Then  part 3$)$ shows that for any $n \in \hN$
  \bea  \label{def_fn}
     g_n (t,\o, y,z) :=  g(t,\o, y,z) + n \(y-L_t(\o)\)^- ,
     \q \fa (t,\o, y,z) \in [0,T] \times \O \times \hR \times \hR^d
    \eea
  defines a  generator.

  \end{rem}

    \section{Main Result: Doubly Reflected BSDEs with Integrable Parameters and Related Dynkin Games}

    \label{sec:DRBSDE}

 The  contribution of  this paper is the following wellposedness result of a  doubly reflected
 BSDE  with integrable parameters  in which
 the $Y-$component of the unique solution represents the value of the related Dynkin game under
 a so-called ``$g-$evaluation" (see Section \ref{sec:g_evaluation}),
 a nonlinear expectation induced by   BSDEs  with the same generator $g$ as the doubly reflected BSDE.
  Like \cite{Hamadene_Hassani_2005}, we assume   the complete separation of  the lower and upper
 obstacles in the doubly reflected BSDE instead of
   the traditional Mokobodski condition which is quite difficult to check in  practice.

    \begin{thm} \label{thm_DRBSDE_exist}
   Let $g$ be a generator. For any $\xi \n  \in \n   L^1(\cF_T)$,
   $L  \n  \in \n   \hS^1_+$ and  $U  \n  \in \n   \hS^1_-$ such that
   $ \hP\{L_T  \n \le \n  \xi  \n \le \n  U_T\} \n = \n \hP\{L_t  \n < \n  U_t, \,  \fa t  \n \in \n  [0,T] \}  \n = \n  1$,
   DRBSDE\,$(\xi,g,L,U)$  admits a unique solution $(Y,Z, K, J)
     \n \in \n   \( \underset{p \in (0,1)}{\cap}  \hS^p \n \) \n
      \n \times \n   \hH^{2,0}  \n \times \n  \hK^0  \n \times \n  \hK^0   $
    such that $Y$  is of class $($D$)$.

 \ss   Define $R (\tau,\ga) \n := \n   \b1_{\{\tau < \ga\} }  L_\tau
    \n + \n \b1_{\{ \ga \le \tau \}\cap \{\ga < T\}}   U_\ga
    \n + \n   \b1_{\{\tau = \ga =  T \}} \xi $, $ \fa \tau,\ga  \n \in \n  \cT $.
     Let $\nu \in \cT $,
     $\tau^*_\nu  \n := \n      \inf \n \big\{t  \n \in \n  [\nu, T] \n
  :   Y_t \n = \n \b1_{\{ t < T\}} L_t  \n + \n  \b1_{\{ t = T\}} \xi \big\}   \n \in \n  \cT_{\nu,T}   $ and
  $\ga^*_\nu  \n := \n      \inf \n \big\{t  \n \in \n  [\nu, T]
  :   Y_t \n = \n \b1_{\{ t < T\}} U_t  \n + \n  \b1_{\{ t = T\}} \xi \big\}  \n \in \n  \cT_{\nu,T}   $.
   It holds for any $\tau, \ga  \in \cT_{\nu,T} $ that
   \bea    \label{eq:b671}
      \cE^g_{\nu, \tau \land  \ga^*_\nu }\big[Y_{\tau \land \ga^*_\nu  } \big]   \le  Y_\nu
      \le  \cE^g_{\nu, \tau^*_\nu \land \ga}\big[Y_{  \tau^*_\nu \land \ga  } \big]  , \q \pas
   \eea
  Consequently, it holds \pas  ~  that
   \bea    \label{eq:b811}
   \esup{\tau \in \cT_{\nu,T}} \cE^g_{\nu,  \tau  \land \ga^*_\nu} \big[R( \tau , \ga^*_\nu)\big]  \n =  \n
    Y_\nu \n =  \n  \cE^g_{\nu, \tau^*_\nu \land \ga^*_\nu}\[ R \big(\tau^*_\nu, \ga^*_\nu \big) \]
     \n = \n  \einf{\ga \in \cT_{\nu,T}} \cE^g_{\nu, \tau^*_\nu \land \ga} \big[R(\tau^*_\nu, \ga)\big] .
   \eea
   In particular, we have
\bea  \label{eq:b537}
      Y_\nu    \n =  \n   \esup{  \tau \in  \cT_{\nu,T}}  \einf{  \ga \in  \cT_{\nu,T}}
        \cE^g_{\nu, \tau \land \ga}\big[R(\tau,\ga) \big]
         \n =  \n   \einf{  \ga \in  \cT_{\nu,T}}  \esup{  \tau \in  \cT_{\nu,T}}
        \cE^g_{\nu, \tau \land \ga}\big[R(\tau,\ga) \big]
   ,  ~ \; \pas
     \eea

   \end{thm}

   \begin{rem}
    $($1$)$ For any $\nu, \z \n \in \n  \cT$ with $0  \n \le \n  \nu  \n \le \n  \z  \n \le \n  \tau^*_0$,
     it is clear that $\tau^*_\nu  \n = \n  \tau^*_0$, \pas ~
      Then \eqref{eq:b671} shows that
      $ Y_\nu  \n \le \n   \cE^g_{\nu, \tau^*_\nu \land \z}\big[Y_{\tau^*_\nu \land \z} \big]
  \n = \n  \cE^g_{\nu,   \ga}\big[Y_{  \ga} \big]  $, \pas, which shows that
   the $Y-$component of   the unique solution of DRBSDE$(\xi,g,L,U)$  is  a $g-$submartingale up to time $ \tau^*_0 $
   \big(see \eqref{def_g_martingale} for  definition of $g-$martingales\big).
  Similarly, $Y$ is a $g-$supermartingale up to time $ \ga^*_0 $.
  Consequently, $Y$ is a $g-$martingale up to time $ \tau^*_0 \n \land \n \ga^*_0 $.

 \ss \no  $($2$)$
   In \eqref{eq:b537}, if we regard $L$ \big(resp. $U$\big) as the amount process a player will receive from,
   or pay to if the   amount is negative,  her opponent
   when the time $\tau$ she chooses to stop the game is earlier \big(resp.  not earlier\big) than
   the stopping time $\ga$ selected by her opponent,
   then the $Y-$component of  the unique solution of DRBSDE$(\xi,g,L,U)$ is exactly
   the player's value   of the Dynkin game  under the   $g-$evaluation.
   If the game starts at   $\nu \in \cT$, \eqref{eq:b811} shows that
       the first time $\tau^*_\nu$ when the value process   $ Y$ meets     $L$  after $\nu$
   and the first time $\ga^*_\nu$ when    $ Y$ meets     $U$  after $\nu$
   form a saddle point of the game.

 \end{rem}

    \section{BSDEs with Integrable Parameters}

    \label{sec:BSDE}

    The derivation of Theorem \ref{thm_DRBSDE_exist}  is based on
    the wellposedness result of   BSDEs with integrable terminal data,
    i.e.  Theorem 6.2 and   6.3 of \cite{BH_Lp_2003}   cited below as Proposition \ref{prop_BSDE_exist}.
    Then in  Section \ref{sec:RBSDE}, we will exploit  the {\it penalization} method to construct  a unique solution
    of the corresponding reflected BSDEs with integrable parameters,
    with which we can adopt \cite{Hamadene_Hassani_2005}'s approach  of  pasting local solutions
    to obtain Theorem \ref{thm_DRBSDE_exist}.

\begin{prop} \label{prop_BSDE_exist}
 Let $g$ be a generator. For any $\xi \in L^1(\cF_T)$,   BSDE$(\xi,g)$   admits a unique solution
 $(Y,Z) \in   \underset{p \in (0,1)}{\cap} ( \hS^p \times  \hH^{2,p} ) $ such that $Y$ is of class $($D$)$.

\end{prop}

 This wellposedness result leads to  a general martingale representation theorem:

 \begin{cor}   \label{cor_martingale}
 For any $\xi \in L^1(\cF_T)$, there exists a unique $Z \in \underset{p \in (0,1)}{\cap}   \hH^{2,p} $
 such that \pas
 \bea   \label{eq:a235}
 \hE[\xi|\cF_t] = \hE[\xi] + \int_0^t Z_s dB_s, \q t \in [0,T] .
 \eea
 \end{cor}

 Proposition \ref{prop_BSDE_exist} also gives rise to  ``$g-$evaluation/expectation" (see  next section),
  a nonlinear expectation under which the value of  optimal stopping problem (resp. Dynkin game) solves
  the corresponding reflected  BSDE (resp. double reflected BSDE) with generator $g$,
  see \eqref{eq:b437} \big(resp. \eqref{eq:b537}\big).

 To derive a corresponding comparison result of Proposition \ref{prop_BSDE_exist} (which
 is crucial for  the penalty method in solving reflected BSDEs with integrable parameters),
 we need the following  mere  generalization of Lemma 2.2 of \cite{BH_Lp_2003}
 (cf. Corollary 1 of \cite{Lp_RBSDE}):

 \begin{lemm} \label{lem_p_power}
 Given $ V \in \hV^0 $,  if $(Y,Z) \in \hS^0 \times \wt{\hH}^{2,0}$
 satisfies that \pas
 \beas
 Y_t = Y_0 + V_t  - V_0 + \int_0^t Z_s d B_s, \q t \in [0,T] ,
 \eeas
 then it holds for any $p \in (1,\infty)$ that \pas
 \beas
 ~ |Y_t|^p  \n    = \n   |Y_0|^p  \n  + \n   p  \n   \int_0^t \sgn(Y_s) |Y_s|^{p-1}  d V_s
  \n  + \n   p  \n   \int_0^t \sgn(Y_s) |Y_s|^{p-1}   Z_s dB_s
  \n  + \n   \frac{p(p \n  - \n  1)}{2}  \n    \int_0^t \b1_{\{Y_s \ne 0\}} |Y_s|^{p-2}  |Z_s|^2 ds ,  \q t \in [0,T] .
 \eeas

 \end{lemm}

 With help of Lemma \ref{lem_p_power}, we can deduce a general comparison result for  BSDEs over stochastic intervals,
 which is critical   in proving Theorem \ref{thm_RBSDE_exist} and
 our main result, Theorem \ref{thm_DRBSDE_exist}:

        \begin{prop} \label{prop_BSDE_comp_basic}
    Given $\nu, \tau \n \in \n  \cT$ with $\nu  \n \le \n  \tau$,
    for $i \n = \n 1,2$    let  $g^i  \n   : \n [0,T] \n  \times \n   \O  \n  \times \n   \hR  \n  \times \n   \hR^d
 \n \to \n  \hR $ be an  $\sP  \n  \otimes \n   \sB(\hR)  \n  \otimes \n     \sB(\hR^d)  / \sB(\hR) -$   measurable function
 and    let  $ \(Y^i ,Z^i , V^i   \)  \n \in \n  \hS^0  \n \times \n  \hH^{2,0}  \n \times \n  \hV^0 $
    such that $\{Y^i_\ga\}_{\ga \in \cT_{\nu,\tau}}$ is uniformly integrable,     that
     $\hE   \Big[ \big( \int_\nu^\tau \n |Z^i_t|^2 \\ dt   \big)^{p/2} \Big] \n < \n \infty  $
     for some $p  \n \in \n  (\a,1)$, and that \pas
  \bea  \label{eq:b211}
  Y^i_t = Y^i_\tau + \int_t^\tau   g^i (s, Y^i_s, Z^i_s) ds + V^i_\tau - V^i_t
  - \int_t^\tau Z^i_s dB_s, \q \fa  t \in [\nu,\tau] .
  \eea
  Assume that  $      Y^1_\tau \le Y^2_\tau   $, \pas ~ and that \pas
  \bea \label{eq:b215}
  \int_t^s \b1_{\{Y^1_r > Y^2_r\}}  ( d V^1_r - d V^2_r ) \le 0  , \q \fa t, s \in [ \nu , \tau ] \hb{ with }  t < s   .
  \eea
  For either $i \n  = \n  1$ or $i \n  = \n  2$,      if \, $g^i$ satisfies \big(H1\big), \big(H2\big), \big(H5\big)
  and if $  g^1 (t,  Y^{3-i}_t ,Z^{3-i}_t ) \n  \le \n    g^2 (t, Y^{3-i}_t ,Z^{3-i}_t )  $,  \dtp ~
  on the stochastic interval $ [ \n   [ \nu,\tau ] \n ]
  \n := \n \big\{(t,\o)  \n \in \n  [0,T] \n \times \n \O  \n : \nu (\o)  \n \le \n  t \n \le  \n \tau (\o) \big\}$,
     then  it holds $\hP -$a.s. that   $ Y^1_t  \n  \le \n   Y^2_t$  for any $t \in [\nu,\tau]$.

    \end{prop}

 Applying Proposition \ref{prop_BSDE_comp_basic} over period  $[0,T]$ with $V^1 \n = \n V^2  \n \equiv \n  0$,
 we obtain  the following comparison result for BSDEs
 whose $Y-$solutions are of  class $($D$)$
 and whose $Z-$solutions are of $\hH^{2,p}$ for some $p \in (\a,1)$.

  \begin{prop}    \label{prop_BSDE_comp}
  For $i \n  = \n  1,2$,
    given  parameter pair   $ \(\xi_i,g^i  \) $ with $\xi_1 \n  \le \n   \xi_2$, \pas,
  let  $ \(Y^i ,Z^i    \)  $   be a solution
   of BSDE$ \(  \xi_i,g^i  \)$ such that $Y^i$ is of class $($D$)$
   and $Z^i \in \underset{p \in (\a,1)}{\cup}    \hH^{2,p}  $.
     For either $i \n  = \n  1$ or $i \n  = \n  2$,      if \, $g^i$ satisfies \big(H1\big), \big(H2\big), \big(H5\big)
  and if $  g^1 (t,  Y^{3-i}_t ,Z^{3-i}_t ) \n  \le \n    g^2 (t, Y^{3-i}_t ,Z^{3-i}_t )  $,  \dtp,
     then  it holds $\hP -$a.s. that   $ Y^1_t  \n  \le \n   Y^2_t$  for any $t \in [0,T]$.

    \end{prop}

\section{  $g-$Evaluations and $g-$Expectations }
\label{sec:g_evaluation}

 Let $g$ be a generator. For any $\tau \n \in \n  \cT$, since the function $g_\tau$ defined in \eqref{def_g_tau}
 is a generator, Proposition \ref{prop_BSDE_exist} shows that for any $\xi \n \in \n L^1(\cF_T)$,
 the BSDE$(\xi,g_\tau)$ admits a unique solution
 \bea \label{eq:xax011}
 \big(Y^{\tau, \xi}, Z^{\tau, \xi}\big)  \n \in  \n   \underset{p \in (0,1)}{\cap} ( \hS^p  \n \times \n  \hH^{2,p} )
 \eea
 such that $Y^{\tau, \xi}$ is of class (D).
 Then we can introduce the   notion of ``$g-$evaluation/expectation",
 which slightly generalizes the one initiated  in  \cite{Peng-97} and \cite{Peng_2004}:

 \begin{deff}

  A family of operators $  \cE^g_{\nu,\tau} \n :
  L^0(\cF_\tau)  \n \to \n  L^0(\cF_\nu) $, $\nu  \n \in \n  \cT,    \tau  \n \in \n  \cT_{\nu,T} $
  is called a ``$g-$evaluation" 
  if for any $\nu,\tau  \n \in \n  \cT$
  with $\nu  \n \le \n  \tau$ and any $\xi  \n \in \n  L^0(\cF_\tau)$,
   \beas
    \cE^g_{\nu,\tau}[\xi]   \n := \n
    \left\{
    \ba{ll}
     Y^{\tau, \xi}_\nu     \n \in \n  L^1(\cF_\nu)   & \hb{if } \xi \in  L^1(\cF_\tau) ; \\
  - \infty, & \hb{if } \hE[\xi^-] = \infty; \\
  \infty,   & \hb{if } \hE[\xi^-] < \infty \hb{ and } \hE[\xi^+] = \infty.
     \ea
     \right.
    \eeas
   In particular, for any $\nu \n \in \n  \cT$ and $\xi  \n \in \n  L^0(\cF_T)   $
    we  refer to  $\cE^g[\xi|\cF_\nu]  \n := \n  \cE^g_{\nu,T}[\xi]$
   as    ``$g-$expectation" of $\xi$ conditional on  the $\si-$ field   $   \cF_\nu $.

 \end{deff}

  \begin{rem} \label{rem_g_evalu}
  If $g$ is independent of $(y,z)$, i.e.,
  if $\{g_t\}_{t \in [0,T]}$ is an $\bF-$progressively measurable process
  with $\hE \n \int_0^T \n  |g_t| dt  \n < \n  \infty$,
  then   for any $\nu \n \in \n  \cT$, $\tau  \n \in \n  \cT_{\nu,T}$
    \bea \label{eq:b641}
    \cE^g_{\nu,\tau}[\xi ]  \n  = \n  \hE \n \[ \xi \n + \n  \int_\nu^\tau \n  g_t dt \, \bigg|\cF_\nu \] , ~ \pas ,
    \q \fa  \xi  \n \in \n  L^0(\cF_\tau)   .
    \eea
 When $g \equiv 0$,   the $g-$expectation degenerates into the classic linear expectation, i.e.
 for any $\nu \n \in \n  \cT$ and $\xi  \n \in \n  L^0(\cF_T)   $,
    $\cE^g[\xi|\cF_\nu]  \n  = \n  \hE [\xi|\cF_\nu]$, \pas

 \end{rem}

 In light of Proposition \ref{prop_BSDE_comp} and the uniqueness
 result in Proposition  \ref{prop_BSDE_exist}, one can  deduce that
 $g-$evaluation with domain $L^1(\cF_T)$  inherits the following
 basic properties from  the classic linear expectation:
 Let $\nu$,  $\tau \n \in \n  \cT$ with $\nu  \n \le \n    \tau$

 \ss \no (1)  ``Monotonicity": For any $\xi,\eta \in L^0(\cF_\tau)$ with $\xi \le  \eta$, \pas ~
  we have $\cE^g_{\nu,\tau}[\xi] \le \cE^g_{\nu,\tau}[\eta]$, \pas;

  \ss \no (2)  ``Time-consistency":  For any   $\ga  \in \cT_{\nu,\tau}$  and $\xi \in L^1(\cF_\tau)$,
$\cE^g_{\nu,\ga}\big[\cE^g_{\ga,\tau}[\xi]\big]=\cE^g_{\nu,\tau}[\xi]$, \pas;

  \ss \no  (3)  ``Constant-Preserving": If it holds \dtp ~ that $g(t,y,0) \n = \n  0$, $\fa y  \n \in \n  \hR$,
  then  $\cE^g_{\nu,\tau}[\xi] \n = \n \xi$, $\pas$ for any $  \xi  \n \in \n  L^1(\cF_\nu)$;

  \ss \no  (4)  ``Zero-one Law": For any $\xi  \n \in \n  L^1(\cF_\tau)$
and $A  \n \in \n  \cF_\nu$, we have $\b1_A\cE^g_{\nu,\tau}[\b1_A \xi] \n = \n \b1_A
\cE^g_{\nu,\tau}[\xi]$, \pas; In addition, if $g(t,0,0) \n = \n 0$, \dtp,
then $\cE^g_{\nu,\tau}[\b1_A \xi] \n = \n \b1_A \cE^g_{\nu,\tau}[\xi]$, \pas;

  \ss \no  (5)  ``Translation Invariant": If $g$ is independent of $y$, then
  $ \cE^g_{\nu,\tau}[\xi \n + \n \eta] \n = \n \cE^g_{\nu,\tau}[\xi] \n + \n \eta$, \pas ~
  for any $\xi  \n \in \n  L^0(\cF_\tau) $ and $ \eta  \n \in \n  L^1(\cF_\nu)$.

 We can define the corresponding $g-$martingales as usual:
 A $\sB \big([0,T]\big) \otimes \cF-$measurable process $X$ of class (D) is called
 an $g-$submartingale \big(resp. $g-$supermartingale or $g-$martingale\big) if for any $0 \le t \le s \le 0  $
  \bea   \label{def_g_martingale}
  \cE^g_{t,s} [X_s] \ge (\hb{resp.}\; \le \, \hb{or} \, =) \, X_t , \q \pas
  \eea
 The   $g-$martingales  possess many classic martingale  properties such as {\it Upcrossing
inequality}, {\it Optional sampling theorem},  {\it Doob-Meyer decomposition} and etc,
which relate the $g-$evaluation closely to risk measures in mathematical finance
\big(see \cite{pln}, \cite{Eman} for the case of Lipschitz $g-$evaluation with domain $L^2(\cF_T)$ and
see \cite{MaYao_2010}, \cite{HMPY-07} for the case of quadratic $g-$evaluation with domain $L^\infty(\cF_T)$\big).
 Due to the page limitation, we will   elaborate neither on  the martingale properties of our $g-$evaluation
 with domain $L^1(\cF_T)$ nor on the connection of this $g-$evaluation to risk measures  in the present paper.

   \section{Reflected BSDEs with Integrable Parameters and Related Optimal Stopping Problems}
   \label{sec:RBSDE}

 With Proposition \ref{prop_BSDE_exist} and Proposition \ref{prop_BSDE_comp},
 we can  employ the {\it penalization} method to obtain,
 as an intermediate step towards  our goal (Theorem \ref{thm_DRBSDE_exist}),
  the following wellposedness result of a   reflected  BSDE  with integrable parameters,
   in which  the $Y-$component of the unique solution stands for the value of the related
 optimal stopping problem under   $g-$evaluation.

\begin{thm} \label{thm_RBSDE_exist}
 Let $g$ be a generator. For any $\xi \n  \in \n   L^1(\cF_T)$
 and $L  \n  \in \n   \hS^1_+$ with $ L_T  \n  \le \n   \xi $, \pas,
   RBSDE$(\xi,g,L)$  admits a unique solution $(Y,Z,K)
    \in   \underset{p \in (0,1)}{\cap} ( \hS^p \times  \hH^{2,p} \times \hK^p ) $
    such that $Y$  is of class $($D$)$.

 Define $\cR_t \n := \n  \b1_{\{ t < T\}} L_t  \n + \n  \b1_{\{ t = T\}} \xi $, $  t  \n \in \n  [0,T]$.
 Let $\nu \in \cT $ and   $\tau_\sharp (\nu)  \n := \n      \inf\big\{t  \n \in \n  [\nu, T]
  :   Y_t \n = \n \cR_t\big\}  \n \in \n  \cT_{\nu,T}   $.      It holds for any $\ga \in \cT_{\nu,T} $ that
   \bea \label{eq:b661}
   \cE^g_{\nu, \ga}\big[Y_\ga \big]
   \le  Y_\nu =  \cE^g_{\nu, \tau_\sharp(\nu) \land \ga}\big[Y_{\tau_\sharp(\nu) \land \ga} \big] ,
   \q \pas
   \eea
   In particular, we have
\bea  \label{eq:b437}
      Y_\nu =   \esup{  \ga \in  \cT_{\nu,T}}    \cE^g_{\nu, \ga}\big[\cR_\ga \big]
  = \cE^g_{\nu, \tau_\sharp(\nu)}\big[ \cR_{\tau_\sharp(\nu) } \big]  ,  \q \pas
     \eea

\end{thm}

  \begin{rem}
 $($1$)$ In view of \eqref{eq:b661},
 the $Y-$component of   the unique solution of RBSDE$(\xi,g,L)$  is a $g-$supermartingale.
 For any $\nu, \tau \n \in \n  \cT$ with $0  \n \le \n  \nu  \n \le \n  \tau  \n \le \n  \tau_\sharp (0)$,
 it is clear that $\tau_\sharp (\nu)  \n = \n  \tau_\sharp (0)$, \pas ~
 Then we have $ Y_\nu  \n = \n   \cE^g_{\nu, \tau_\sharp(\nu) \land \ga}\big[Y_{\tau_\sharp(\nu) \land \ga} \big]
  \n = \n  \cE^g_{\nu,   \ga}\big[Y_{  \ga} \big]  $, \pas, which shows that
  $Y$ is  a $g-$martingale up to time $ \tau_\sharp (0) $.

 \ss \no  $($2$)$    In \eqref{eq:b437}, if we regard $\cR$ as a reward process that include
   a running reward $L$ and a terminal reward $\xi$, then
    the $Y-$component of   the unique solution of RBSDE$(\xi,g,L)$
    is exactly   the Snell envelope of  $\cR$ under the   $g-$evaluation.
   Given a start time   $\nu  \n \in \n  \cT$,  the first time $\tau_\sharp (\nu)$  when  $ Y$ meets   $\cR$
   after $\nu$   is an   optimal   stopping time for a player to choose
   if she is aimed to maximize her expected reward   under $g-$expectation.

 \end{rem}

    To derive the existence result in Theorem \ref{thm_RBSDE_exist},
    we  will   use  penalization  method which can be summarized in the following two monotonicity results:

         \begin{prop}   \label{prop_monotone_result_0}

   Let $L  \n  \in \n   \hS^1_+$ and let
    $g  \n  : [0,T]  \n \times \n  \O  \n \times \n  \hR  \n \times \n  \hR^d \to \hR$ be a
   $\sP    \n   \otimes   \n    \sB(\hR)    \n   \otimes   \n    \sB(\hR^d)/\sB(\hR) -$measurable function
   satisfying   $($H1$)$, $($H4$)$ and $($H5$)$.
  For any $n \n \in \n  \hN$, consider the function $g_n$   defined in \eqref{def_fn} and
   let   $(Y^n \n ,Z^n \n , J^n) \n \in  \n \(\underset{p \in (0,1)}{\cap} \hS^p \)
    \n \times \n  \hH^{2,0} \n \times \n \hK^0 $
  such that $Y^n$ is of class $($D$)$ and   that \pas ~
  \beas
  Y^n_t = Y^n_T + \int_t^T g_n (s, Y^n_s, Z^n_s) ds - J^n_T + J^n_t
  - \int_t^T Z^n_s dB_s, \q    t \in [0, T] .
  \eeas
     If $ \{ Y^n \}_{ n   \in   \hN}$ is an increasing sequence of processes,
   then its limit    $ Y_t  \n  := \n  \lmtu{n \to \infty}     Y^n_t   $,
    $ t  \n  \in \n  [0,T]   $ is an $\bF-$predictable process of class $($D$)$ that  satisfies
    $ \hE \n \[ \underset{t \in [0,T]}{\sup} |Y_t|^p \]  \n < \n \infty$, $\fa p  \n \in \n  (0,1)$.

   \end{prop}

     \begin{prop}   \label{prop_monotone_result}

   Let  $L  \n  \in \n   \hS^1_+$,    let
    $g  : [0,T] \times \O \times \hR \times \hR^d \to \hR$ be a
   $\sP     \otimes     \sB(\hR)     \otimes     \sB(\hR^d)/\sB(\hR) -$measurable function satisfying
   $($H1$)$$-$$($H4$)$, and let $\nu, \tau \in \cT$ with $\nu \le \tau$.
  For any $n \n \in \n  \hN$,  consider the function  $g_n$ defined in \eqref{def_fn} and
  let   $(Y^n,Z^n) \n \in  \n   \hS^0  \n \times \n  \hH^{2,0}  $ satisfies that \pas ~
  \bea  \label{eq:a243}
  Y^n_t = Y^n_\tau + \int_t^\tau g_n (s, Y^n_s, Z^n_s) ds
  - \int_t^\tau Z^n_s dB_s, \q \fa  t \in [\nu,\tau] .
  \eea
     If $\big\{\b1_{\{t \ge \nu \}} Y^n_{ \tau \land t }\big\}_{t \in [0,T]}$, $n \n \in \n \hN$
     is an increasing sequence of processes
   whose limit    $ Y_t  \n  := \n  \lmtu{n \to \infty}  \b1_{\{t \ge \nu \}}  Y^n_{ \tau \land t }   $,
    $ t  \n  \in \n  [0,T]   $
  satisfies $ \hP \{ Y_\tau \n \ge \n L_\tau \} \n = \n \hP \Big\{
   \underset{t \in [\nu,\tau]}{\sup} \((Y^1_t)^-  \n + \n  Y^+_t\)  \n < \n  \infty \Big\} \n = \n 1 $,
  then process  $\{Y_{\nu \vee t} \}_{t \in [0,T]}   $ has \pas ~ continuous paths
  and there exist $(Z,K)  \n \in \n  \hH^{2,0}  \n \times \n  \hK^0 $
  such that \pas
  \bea  \label{RBSDEL_tau}
  \begin{cases}
  \dis  L_t \le  Y_t
  = Y_\tau + \int_t^\tau g(s, Y_s, Z_s) ds
  + K_\tau - K_t
  - \int_t^\tau Z_s dB_s, \q \fa t \in [\nu,\tau] ,   \vspace{1mm}  \\
  \dis \int_\nu^\tau (Y_t-L_t) d K_t = 0 .
  \end{cases}
  \eea

   \end{prop}

 On the other hand,  the uniqueness result in Theorem \ref{thm_RBSDE_exist}
 follows from the following   comparison result  for reflected BSDEs
 whose $Y-$solutions are of  class $($D$)$
 and whose $Z-$solutions are of $\hH^{2,p}$ for some $p \in (\a,1)$.

      \begin{prop} \label{prop_RBSDE_comp}
    For $i=1,2$, given   parameter pair $\(\xi_i,g^i  \)$ and $L^i \n  \in \n   \hS^1_+ $ such that
    $\hP\{  L^i_T  \n  \le \n   \xi_i  \} = \hP\{   \xi_1 \le \xi_2   \}
    = \hP\{   L^1_t \le  L^2_t , \; \fa t \in [0,T]  \} = 1$,
    let  $ \(Y^i ,Z^i , K^i   \)  $   be a solution
   of RBSDE$ \(  \xi_i,g^i,L^i  \)$ such that $Y^i$ is of class $($D$)$
   and $Z^i \in \underset{p \in (\a,1)}{\cup}    \hH^{2,p}  $.
     For either $i \n  = \n  1$ or $i \n  = \n  2$,      if \, $g^i$ satisfies \big(H1\big), \big(H2\big), \big(H5\big)
  and if $  g^1 (t,  Y^{3-i}_t ,Z^{3-i}_t ) \n  \le \n    g^2 (t, Y^{3-i}_t ,Z^{3-i}_t )  $,  \dtp,
     then  it holds $\hP -$a.s. that   $ Y^1_t  \n  \le \n   Y^2_t$  for any $t \in [0,T]$.

    \end{prop}

 \begin{rem} \label{rem_RBSDEU2}
 By   Remark  \ref{rem_assum} $($1$)$,
 one can apply Theorem \ref{thm_RBSDE_exist}, Proposition \ref{prop_monotone_result} and
 Proposition \ref{prop_RBSDE_comp} to $g_-$ \big(defined in \eqref{eq:g_neg}\big) to obtain a version of them
    for the reflected BSDE with upper obstacle like \eqref{RBSDEU}.

 \end{rem}

\section{Proofs}

\label{sec:proofs}

 \subsection{Proofs of the results in Section \ref{sec:BSDE} and \ref{sec:g_evaluation} }

    \no {\bf Proof of Proposition \ref{prop_BSDE_exist}:}
    As condition (H7) of \cite{BH_Lp_2003} is automatically satisfied,
 it suffices to verify  condition (H5) therein, i.e.,
  Given $r \n  \ge \n   0$,
   \beas
  \hb{the process \; $  \psi^r_t(\o)  \n  := \n   \underset{|y| \le r}{\sup}
  |g(t,\o,y,0)  \n  - \n   g(t,\o,0,0) | $,   $ (t,\o)  \n  \in \n   [0,T]  \n  \times \n   \O$ \; is integrable. }
  \eeas
 By (H3), it holds $dt    \otimes    d\hP-$a.s. that
 $\psi^r_t(\o)  \n  = \n   \underset{ y  \in [- r,r] \cap \hQ}{\sup}
 |g(t,\o,y,0)  \n  - \n   g(t,\o,0,0) | $,
 which implies that  $\psi^r$ is $\bF-$progressively measurable. Also, (H4) shows that $dt \n  \otimes \n  d\hP-$a.s.,
  $ \psi^r_t(\o)  \n  \le \n    |  g(t,\o,0,0) |  \n  + \n   \underset{|y| \le r}{\sup} |g(t,\o,y,0)|
   \n  \le \n   2 h_t(\o)  \n  + \n   \k r  $.
 It follows that $\psi^r_t$   belongs to $\hL^1 ([0,T]  \n  \times \n   \O, \sB([0,T])  \n  \otimes \n   \cF,
  dt  \n  \otimes \n   \hP )$.   \qed

    \ss \no {\bf Proof of Corollary \ref{cor_martingale}:}
 Clearly, $ \fg (t,\o,y,z) \n : =  \n  0 $, $\fa (t,\o,y,z)  \n \in  \n
  [0,T]  \n \times \n  \O  \n \times \n  \hR  \n \times \n  \hR^d $
 is a generator. In light of Proposition \ref{prop_BSDE_exist},
    BSDE$(\xi,0)$ admits a unique solution
 $(Y,Z) \in   \underset{p \in (0,1)}{\cap} ( \hS^p \times  \hH^{2,p} ) $ such that $Y$ is of class $($D$)$.
 For any $n  \n  \in \n   \hN$, we  define   stopping time
 $  \tau_n  \n  := \n \inf\big\{ t  \n \in \n  [0,T] \n : \int_0^t \n  |Z_s|^2 ds  \n > \n  n \big\}
  \n \land \n  T  \n \in \n  \cT $, and see from  $Z \in \underset{p \in (0,1)}{\cap}   \hH^{2,p} \subset \hH^{2,0} $
  that   $\{\tau_n\}_{n \in \hN}$ is stationary.

  Let $t  \n \in \n  [0,T]$ and $n  \n  \in \n   \hN$.   Since
  $  Y_{\tau_n \land t}  \n = \n  Y_{\tau_n}  \n - \n  \int_{\tau_n \land t}^{\tau_n} Z_s dB_s$, \pas,
  taking conditional  expectation $\hE[ \cd |\cF_t]$ yields that
  \bea \label{eq:a233}
    Y_{\tau_n \land t}  \n = \n  \hE[Y_{\tau_n} |\cF_t ]  , \q  \pas ~
 \eea
  As    $\{\tau_n\}_{n \in \hN}$ is stationary,
 letting $n \to \infty $ in \eqref{eq:a233}, we can deduce from the continuity of $Y$
 and the uniform  integrability of $\{Y_\tau\}_{\tau \in \cT}$   that
  $   Y_t = \hE[ Y_T |\cF_t ] = \hE[ \xi |\cF_t ] $,   \pas ~
  In particular, $Y_0 = \hE[\xi]$. Then
  \beas
  \hE[ \xi |\cF_t ] = Y_t = Y_0 + \int_0^t Z_s dB_s = \hE[\xi] + \int_0^t Z_s dB_s , \q \pas
  \eeas
  This together with the continuity of processes $\big\{\hE[ \xi |\cF_\cd ]\big\}_{t \in [0,T]}$
   and $\big\{\int_0^t Z_s dB_s \big\}_{t \in [0,T]}$ leads  to  \eqref{eq:a235} while the
   uniqueness of process $Z$ is clear.    \qed

       \ss \no {\bf Proof of Proposition \ref{prop_BSDE_comp_basic}:}
  Without loss of generality,     suppose that $g^1$ satisfies (H1), (H2), (H5)  and  that
  \bea   \label{eq:a021}
         g^1 (t, Y^2_t,Z^2_t)  \n  \le \n   g^2 (t, Y^2_t,Z^2_t)    , \q   \dtp \hb{ on } [ \n   [ \nu,\tau ] \n ] .
  \eea
  Set     $(    \cY ,  \cZ  )  \n  : =  \n  (  Y^1 \n  - \n  Y^2, Z^1 \n  - \n  Z^2   )$
  and $q:=p/\a \in (1,1/\a)$.

 \ss \no {\bf (1)} {\it We first show that
 $ \hE \[ \underset{t \in [\nu, \tau]}{\sup}  \( \cY^+_t \)^q  \] < \infty $}.

 \ss  Since $\hE[| \cY_\nu|]  \n \le \n  \hE \big[|  Y^1_\nu|\big]  \n + \n  \hE\big[| Y^2_\nu|\big] \n < \n \infty$
 by the uniform integrability of $\{Y^i_\ga\}_{\ga \in \cT_{\nu,\tau}}$, $i \n = \n 1, 2$,
 Corollary \ref{cor_martingale} implies that
 there exists a unique $\wt{\cZ} \in \underset{p' \in (0,1)}{\cap}   \hH^{2,p'} $
 such that $ \hP \big\{ \hE[\cY_\nu|\cF_t] \n = \n  \hE[\cY_\nu]
  \n + \n  \int_0^t  \n \wt{\cZ}_s dB_s, ~ \fa t  \n \in \n  [0,T] \big\}  \n = \n  1 $.
 This together with     \eqref{eq:b211} shows  that \pas
 \bea
 \hspace{-3mm}  \wt{\cY}_t & \tn \dn: =& \tn \dn \hE[\cY_\nu|\cF_{\nu \land t}]
  \n  + \n  \cY_{\nu \vee (\tau \land t)}   \n - \dn  \cY_\nu
   \n =  \n  \hE[\cY_\nu]  \n + \dn  \int_0^{\nu \land t} \n  \wt{\cZ}_s dB_s
  \n - \dn  \int_\nu^{\nu \vee (\tau \land t)}  \n  \D g_s  ds
  \n - \n  V^1_{\nu \vee (\tau \land t)} \n + \n  V^1_\nu
  \n  +  \n    V^2_{\nu \vee (\tau \land t)}  \n - \n  V^2_\nu
  \n + \n  \int_\nu^{\nu \vee (\tau \land t)}  \n   \cZ_s  d B_s \nonumber \\
  & \tn \dn  = & \tn \dn  \hE[\cY_\nu]
 \n - \n  \int_0^t \n  \b1_{\{\nu < s \le \tau\}}   \D g_s    ds
  \n - \n  \int_0^t \n  \b1_{\{\nu < s \le \tau\}}   ( d V^1_s - d V^2_s )
      \n + \n  \int_0^t  \n  \( \b1_{\{s \le \nu \}} \wt{\cZ}_s
      \n + \n    \b1_{\{\nu < s \le \tau\}}  \cZ_s \) d B_s   , \q t \in [0,T] ,     \label{eq:b241}
 \eea
 where $  \D g_s  \n  : = \n  g^1 (s, Y^1_s,Z^1_s)  \n  - \n   g^2 (s, Y^2_s,Z^2_s) $.
 So $\wt{\cY}$ is an $\bF-$adapted continuous process, i.e. $\wt{\cY}  \n \in \n  \hS^0$.
 Applying It\^o-Tanaka's formula to process $\wt{\cY}^+$ yields that \pas
 \bea
   \wt{\cY}_t^+   & \tn   =    & \tn     \big(\hE[\cY_\nu]\big)^+  \n  -  \n    \int_0^t \n \b1_{\{ \wt{\cY}_s > 0 \}}
 \b1_{\{\nu < s \le \tau\}}   \D g_s ds   \n  +  \frac12  \fL_t
  \n  - \n       \int_0^t \n \b1_{\{ \wt{\cY}_s > 0 \}} \b1_{\{\nu < s \le \tau\}}   ( d V^1_s - d V^2_s ) \nonumber \\
   & \tn  & \tn  +   \int_0^t  \n \b1_{\{ \wt{\cY}_s > 0 \}} \( \b1_{\{s \le \nu \}} \wt{\cZ}_s
           \n   +    \n    \b1_{\{\nu < s \le \tau\}}  \cZ_s \) d B_s   ,   \q     t  \n  \in \n   [0,T ] , \label{eq:a031}
 \eea
 where $\fL$ is an $\bF-$adapted,  continuous increasing process known as the ``local time" of $\wt{\cY}$ at $0$.

 \ss  Let  $n \n \in \n  \hN$. We   define a stopping time
 $  \tau_n \n := \n  \inf\big\{ t  \n \in \n  [\nu,\tau] \n :   \int_\nu^t  \n |\cZ_s|^2 ds  \n > \n  n  \big\}
   \n \land \n  \tau  \n \in \n  \cT_{\nu,\tau} $,  and  integrate by parts the process
 $ \big\{\dis e^{\l^+  (\tau_n \land t)} \wt{\cY}^+_{\tau_n \land t} \big\}_{t \in [0,T]} $ to obtain that \pas
 \bea
 && \hspace{-1.2cm}  e^{\l^+  (\tau_n \land t)} \wt{\cY}^+_{\tau_n \land t}
     \n  =        \n     e^{\l^+  \tau_n  } \cY^+_{ \tau_n }
   \n  +  \n   \int_{\tau_n \land t}^{\tau_n} \b1_{\{ \cY_s > 0 \}}  \b1_{\{s > \nu\}}
   e^{\l^+ s}   \D g_s  ds
  \n  + \n   \int_{\tau_n \land t}^{\tau_n} \n \b1_{\{ \cY_s > 0 \}} \b1_{\{s > \nu\}}
  e^{\l^+  s  } ( d V^1_s - d V^2_s ) \nonumber \\
 &&  -   \frac12 \int_{\tau_n \land t}^{\tau_n} e^{\l^+ s} d  \fL_s
  \n  - \n  \l^+ \int_{\tau_n \land t}^{\tau_n} e^{\l^+  s  } \wt{\cY}^+_s ds
  \n  - \n   \int_{\tau_n \land t}^{\tau_n}  \n   \b1_{\{ \wt{\cY}_s > 0 \}} e^{\l^+ s} \n
  \( \b1_{\{s \le \nu \}} \wt{\cZ}_s
           \n   +    \n    \b1_{\{s > \nu\}}  \cZ_s \)  d B_s    ,
    \q      t  \n  \in \n   [0,T ] .    \q   \label{eq:a023}
 \eea
 Here we used the fact that  $\wt{\cY}_{\nu \vee (\tau \land t)}
   = \hE[\cY_\nu|\cF_\nu  ]  \n  + \n  \cY_{\nu \vee (\tau \land t)}  \n - \n  \cY_\nu
   \n = \n \cY_{\nu \vee (\tau \land t)} $, $\fa t \in [0,T]$, i.e.
   \bea \label{eq:b221}
    \wt{\cY}_t = \cY_t , \q \fa t \in [\nu, \tau] .
   \eea
 Since $g^1$ satisfies (H2) and (H5), it holds   \dsp ~ on $ [ \n   [ \nu,\tau ] \n ] $ that
 \bea
 \hspace{-5mm} \b1_{\{ \cY_s (\o) > 0 \}}   \Big( g^1  \n  \(s, \o, Y^1_s(\o),Z^1_s(\o)\)
  \n   - \n   g^1  \n  \( s, \o, Y^2_s(\o),Z^1_s(\o) \) \n \Big)
  \n   \le    \n  
    \b1_{\{ \cY_s (\o) > 0 \}} \l   \cY^+_s (\o)
 \n  \le   \n     \l^+    \cY^+_s (\o),
   \label{eq:a035}
\eea
 and that
 \beas
 \q \big| g^1 \n \(s, \o, Y^2_s(\o),Z^1_s(\o)\) \n   - \n   g^1 \n  \( s, \o, Y^2_s(\o),Z^2_s(\o) \) \big|
   \n  \le \n   \k  \n \(h_s(\o) \n  + \n  |Y^2_s(\o)| \n  + \n  |Z^1_s(\o)|\)^\a  \n  + \n
  \k  \n \(h_s(\o) \n  + \n  |Y^2_s(\o)| \n  + \n  |Z^2_s(\o)|\)^\a .
 \eeas
   Plugging them back into \eqref{eq:a023} and taking   $t  \n = \n  \nu   \n  \vee  \n t $ there,
   we see from  \eqref{eq:a021} and \eqref{eq:b215} that \pas
   \bea  \label{eq:b217}
     e^{\l^+  (\nu \vee (\tau_n \land t))} \cY^+_{\nu \vee (\tau_n \land t)}
     \n  \le         \n     e^{\l^+  \tau_n  } \cY^+_{ \tau_n }
   \n  +  \n  2  \k e^{\l^+ T } \eta
     - \n   \int_{\nu \vee (\tau_n \land t)}^{\tau_n}  \n   \b1_{\{ \cY_s > 0 \}} e^{\l^+ s}  \cZ_s  d B_s    ,
    ~ \;       t  \n  \in \n   [0,T ] ,
   \eea
   where  $\eta := \int_\nu^\tau \n  \(h_t  \n  + \n  |Y^2_t | \n  + \n  |Z^1_t | \n  + \n  |Z^2_t |\)^\a dt $.

 \ss  Let $ t  \n  \in \n   [0,T ] $.      Taking conditional expectation
   $\hE \big[\cd|\cF_{\nu \vee (\tau_n \land t)} \big]$  in \eqref{eq:b217}     yields   that
 $    e^{\l^+  (\nu \vee (\tau_n \land t))} \cY^+_{\nu \vee (\tau_n \land t)}   \le
       \hE \[ e^{\l^+ \tau_n } \cY^+_{ \tau_n }
   \n  +  \n  2  \k e^{\l^+ T } \eta  \big| \cF_{\nu \vee (\tau_n \land t)} \] $,   \pas,
 and it follows that
 \bea  \label{eq:a025}
 \b1_{\{\nu \le t \le \tau_n\}}  e^{\l^+  t} \cY^+_t   \le  \b1_{\{\nu \le t \le \tau_n\}}
       \hE \[ e^{\l^+ \tau_n } \cY^+_{ \tau_n }
   \n  +  \n  2  \k e^{\l^+ T } \eta  \big| \cF_t \] , \q \pas
 \eea
 By \eqref{eqn-d011} and H\"older's inequality,
  \beas
  \q  \eta    \le     \int_\nu^\tau \n  \( h_t    \n  + \n       |Y^2_t | \)^\a dt
  \n  +    \sum^2_{i=1}  \int_\nu^\tau  \n   | Z^i_t |^\a dt
      \le      T^{1-\a} \(\int_\nu^\tau \n  \(  h_t \n  + \n  |Y^2_t | \) dt \)^\a
   \n  + \n  T^{1- \a/2 } \sum^2_{i=1}       \(  \int_\nu^\tau  \n   | Z^i_t |^2 dt \)^{ \a/2 } , \q \pas
  \eeas
    Fubini's Theorem and the uniform integrability of $ \big\{ Y^2_\ga \big\}_{\ga \in \cT_{\nu,\tau}}$ imply that
     \beas
     \hE \n  \int_\nu^\tau  \n   |Y^2_t | dt   \n   = \n
     \hE \n  \int_\nu^\tau  \n    |Y^2_{\nu \vee (\tau \land t)} | dt
     \le \hE \n  \int_0^T  \n   \big|Y^2_{\nu \vee (\tau \land t)} \big| dt
      \n   = \n   \int_0^T  \n  \hE \[ \big|Y^2_{\nu \vee (\tau \land t)} \big| \] dt
   \n  \le \n   T \underset{\ga \in \cT_{\nu,\tau}}{\sup} \hE \[ |Y^2_\ga | \]  \n  < \n   \infty  .
   \eeas
  As $ q =p/\a $,   applying \eqref{eqn-d011} and H\"older's inequality again yields that
  \bea  \label{eq:b219}
  \hE [\eta^q] \le 3^{q-1} T^{(1-\a)q} \bigg\{ \hE  \int_\nu^\tau \n
  \( h_t \n  + \n  |Y^2_t | \) dt   \bigg\}^{p}
  +  3^{q-1} T^{(1-\a/2)q}  \sum^2_{i=1}  \,  \hE\[  \(  \int_\nu^\tau  | Z^i_t |^2 dt \)^{p/2} \]  < \infty .
  \eea

  \ss  We see from $ \hE   \Big[ \big( \int_\nu^\tau \n |Z^i_t|^2   dt   \big)^{p/2} \Big] \n < \n \infty $,
  $i \n = \n 1,2$
  that for \pas ~ $\o  \n \in \n  \O$, $\tau(\o)  \n = \n  \tau_{N_\o} (\o)$ for some $N_\o  \n \in \n  \hN$.
  For any $t  \n \in \n  [0,T]$, since the uniform integrability of
  $ \big\{ Y^i_\ga \big\}_{\ga \in \cT_{\nu,\tau}} $, $i \n = \n 1,2$ implies
 that of  $\big\{ e^{\l^+ \ga} \, \cY^+_\ga \big\}_{\ga \in \cT_{\nu,\tau}} $,
   letting $n  \n \to \n  \infty$ in  \eqref{eq:a025} yields that \pas
  \beas
  \b1_{\{\nu \le t \le \tau \}}  \cY^+_t   \n  \le  \n   \b1_{\{\nu \le t \le \tau \}}  e^{\l^+  t} \cY^+_t
   \n  \le  \n   \b1_{\{\nu \le t \le \tau \}}  2  \k e^{\l^+ T }
   \hE  \n  \[  e^{\l^+ \tau  } (Y^1_\tau  \n - \n  Y^2_\tau)^+   \n + \n \eta   | \cF_t  \]
    \n = \n  \b1_{\{\nu \le t \le \tau \}}  2  \k e^{\l^+ T } \hE  [  \eta   | \cF_t  ]  .
    \eeas
  Using the continuity of   $\cY^+$   and
 that of process $ \big\{ \hE  [  \eta   | \cF_t  ] \big\}_{t \in [0,T]} $, one gets
 $   \hP \big\{ \cY^+_t   \n  \le  \n  2  \k e^{\l^+ T }  \hE  [  \eta   | \cF_t  ]  ,
  ~   \fa t  \n \in \n  [\nu, \tau] \big\}  \n = \n  1 $.
 Then Doob's martingale inequality and \eqref{eq:b219} lead to that
 \bea    \label{eq:a041}
 \hE \[ \underset{t \in [\nu, \tau]}{\sup}  \( \cY^+_t \)^q  \]
 \le ( 2  \k )^q e^{q\l^+ T }  \hE\[ \underset{t \in [0,T]}{\sup} \(\hE  [  \eta   | \cF_t  ]\)^q  \]
 \le \( \frac{q}{q-1} \)^q( 2  \k )^q e^{q\l^+ T } \hE [    \eta^q   ] < \infty .
 \eea

 \ss \no {\bf (2)} {\it Next, we show that $    \hE \[   \underset{t \in [\nu, \tau]}{\sup}  \( \cY^+_t \)^q  \]  =0  $
 indeed; then the conclusion easily follows}.

  \ss   According to \eqref{eq:a031}, applying Lemma \ref{lem_p_power} yields that \pas
     \beas
   \(\wt{\cY}^+_t\)^q  & \tn = & \dn \dn \(\big(\hE[\cY_\nu]\big)^+\)^q
  \n  - \n   q  \n   \int_0^t  \b1_{\{ \wt{\cY}_s > 0 \}}  \b1_{\{\nu < s \le \tau\}}
   \(\wt{\cY}^+_s\)^{q-1}   \D g_s ds
 \n +  \n  \frac{ q }{2}  \n   \int_0^t \(\wt{\cY}^+_s\)^{q-1}   d \fL_s  \\
  & \tn      & \tn
       -     q  \n   \int_0^t  \b1_{\{ \wt{\cY}_s > 0 \}}  \b1_{\{\nu < s \le \tau\}}  \(\wt{\cY}^+_s\)^{q-1}
   ( d V^1_s \n - \n  d V^2_s )
  \n  +     q  \n   \int_0^t  \b1_{\{ \wt{\cY}_s > 0 \}}  \(\wt{\cY}^+_s\)^{q-1} \n
  \( \b1_{\{s \le \nu \}} \wt{\cZ}_s   \n   +    \n    \b1_{\{\nu < s \le \tau\}}  \cZ_s \) dB_s  \\
  & \tn      & \tn      +     \frac{q(q \n  - \n  1)}{2}  \n
   \int_0^t \b1_{\{\wt{\cY}_s > 0\}} \(\wt{\cY}^+_s\)^{q-2} \n
   \( \b1_{\{s \le \nu \}} \big| \wt{\cZ}_s \big|^2 \n + \n \b1_{\{\nu < s \le \tau\}} |\cZ_s|^2 \) ds
  ,  \q t \in [0,T] .
  \eeas
 Set $ \dis a \n :=   \n  \l^+  \n +    \frac{  \k^2}{1  \n \land \n  (q \n - \n 1)} $
 and let  $n  \n \in \n  \hN$. We define a stopping time
  $ \dis  \ga_n  \n := \n  \inf\bigg\{ t  \n \in \n  [\nu,\tau] \n : \n \underset{s \in [\nu, t]}{\sup}  \cY^+_s
    + \n  \int_\nu^t |\cZ_s|^2 ds  \n > \n  n  \bigg\}     \land    \tau  \n \in \n  \cT_{\nu,\tau} $,
    and  integrate by parts the process
 $ \Big\{     e^{a q (\ga_n \land t)} \( \wt{\cY}^+_{\ga_n \land t} \)^q \Big\}_{t \in [0,T]} $
 to obtain that \pas
 \beas
     & \tn    \n  & \hspace{-1cm}   e^{a q (\ga_n \land t)} \( \wt{\cY}^+_{\ga_n \land t} \)^q  \n  + \n
  \frac{q(q \n  - \n  1)}{2}  \n    \int_{\ga_n \land t}^{\ga_n}  \b1_{\{\wt{\cY}_s > 0\}}
   e^{a q s}  \(\wt{\cY}^+_s\)^{q-2} \n
   \( \b1_{\{s \le \nu \}} \big| \wt{\cZ}_s \big|^2 \n + \n \b1_{\{s > \nu\}} |\cZ_s|^2 \) ds  \nonumber \\
   &\tn    \n  =    &\tn   \n    e^{a q  \ga_n  } \(  \cY^+_{\ga_n  } \)^q
   \n  +  \n  q \int_{\ga_n \land t}^{\ga_n} \b1_{\{ \cY_s > 0 \} \cap \{s > \nu\}}
   e^{a q s} \(\cY^+_s\)^{q-1}   \D g_s ds
  \n  - \n  a q \int_{\ga_n \land t}^{\ga_n}    e^{a q  s  } \( \wt{\cY}^+_s \)^q   ds
  \n  -    \frac{ q }{2}  \n   \int_{\ga_n \land t}^{\ga_n} e^{a q s} \(\wt{\cY}^+_s\)^{q-1}   d \fL_s  \nonumber \\
  & \tn      & \tn
   +     q  \n   \int_{\ga_n \land t}^{\ga_n}   \b1_{\{ \cY_s > 0 \} \cap \{s > \nu\}}
     e^{a q  s  }  \n \(\cY^+_s\)^{q-1} \n ( d V^1_s \n - \n  d V^2_s )
  \n  -  \n   q  \n   \int_{\ga_n \land t}^{\ga_n}
   \b1_{\{ \wt{\cY}_s > 0 \}} e^{a q s}  \n  \(\wt{\cY}^+_s\)^{q-1}  \n
   \( \b1_{\{s \le \nu \}} \wt{\cZ}_s   \n   +    \n    \b1_{\{s > \nu\}}  \cZ_s \) dB_s
  ,  ~  t \n \in \n  [0,T] .
 \eeas
Then \eqref{eq:b221}, \eqref{eq:a021},  \eqref{eq:a035} and \eqref{eq:b215} imply that \pas
   \beas
 \q    & \tn    \n  & \hspace{-1.2cm}   e^{a q t} \( \cY^+_t \)^q  \n  + \n
  \frac{q(q \n  - \n  1)}{2}  \n    \int_t^{\ga_n}  \b1_{\{\cY_s > 0\}}
   e^{a q s}  \(\cY^+_s\)^{q-2}   |\cZ_s|^2   ds  \nonumber \\
   &\tn    \n  =    &\tn   \n    e^{a q  \ga_n  } \( \cY^+_{\ga_n  } \)^q
   \n  +  \n  q \int_t^{\ga_n} \b1_{\{ \cY_s > 0 \}}
   e^{a q s} \(\cY^+_s\)^{q-1}      \D g_s ds
  \n  - \n  a q \int_t^{\ga_n}    e^{a q  s  } \( \cY^+_s \)^q   ds
  \n  - \n  \frac{ q }{2}  \n   \int_t^{\ga_n} e^{a q s} \(\cY^+_s\)^{q-1}   d \fL_s    \nonumber \\
  & \tn      & \tn
  \n  +     q  \n   \int_t^{\ga_n}  \b1_{\{ \cY_s > 0 \}}  e^{a q  s  } \(\cY^+_s\)^{q-1}
      ( d V^1_s \n - \n  d V^2_s )
  \n  - \n   q  \n   \int_t^{\ga_n}
   \b1_{\{ \cY_s > 0 \}} e^{a q s}  \(\cY^+_s\)^{q-1}     \cZ_s   dB_s  \nonumber \\
   &\tn    \n  \le    &\tn   \n      e^{a q  \ga_n  } \( \cY^+_{\ga_n  } \)^q
   \n  +  \n  q \int_t^{\ga_n} \b1_{\{ \cY_s > 0 \}}
   e^{a q s} \(\cY^+_s\)^{q-1}  \( g^1 (s, Y^2_s,Z^1_s)  \n  - \n   g^1 (s, Y^2_s,Z^2_s) \) ds
  \n  - \n  \frac{q \k^2}{1  \n \land \n  (q \n - \n 1)}  \int_t^{\ga_n}  e^{a q  s } \( \cY^+_s \)^q  ds  \nonumber  \\
  & \tn      & \tn
  \n  -    q  \n   \int_t^{\ga_n}  \b1_{\{ \cY_s > 0 \}} e^{a q s}  \(\cY^+_s\)^{q-1}  \cZ_s dB_s
  ,  \q \fa t \in [\nu,\ga_n] .
 \eeas
   Since
$ \b1_{\{ \cY_t > 0 \}}  \k  \(\cY^+_t \)^{q-1}  |  \cZ_t   | \n \le \n
   \frac{q-1}{4} \b1_{\{ \cY_t > 0 \}} \n  \(\cY^+_t \)^{q-2} | \cZ_t  |^2
    \n + \n  \frac{\k^2}{q-1}  \n \(\cY^+_t \)^q $, $ \fa  t  \n \in \n  [\nu,\tau] $,
 we can deduce from (H1)   that \pas
 \beas
     e^{a q t} \( \cY^+_t \)^q   \n  + \n   \frac{q(q \n  - \n  1)}{4}  \n
 \int_t^{\ga_n}  \n   \b1_{\{\cY_s > 0\}} e^{a q s} \(\cY^+_s\)^{q-2}  |\cZ_s|^2 ds
     \n  \le     \n     e^{a q  \ga_n  } \( \cY^+_{\ga_n  } \)^q
    \n  -  q  \n   \int_t^{\ga_n}  \n   \b1_{\{ \cY_s > 0 \}} e^{a q s} \(\cY^+_s\)^{q-1}  \cZ_s dB_s
    ,       ~      \fa t \in [\nu,\ga_n] .
 \eeas

  Taking expectation   for $t=\nu$ shows that
   \bea  \label{eq:a039}
 \frac{q(q \n  - \n  1)}{4} \hE
 \int_\nu^{\ga_n}  \n   \b1_{\{\cY_s > 0\}} e^{a q s} \(\cY^+_s\)^{q-2}  |\cZ_s|^2 ds
 \le \hE \Big[   e^{a q  \ga_n  } \( \cY^+_{\ga_n  } \)^q   \Big]   .
   \eea
   On the other hand, the Burkholder-Davis-Gundy inequality implies that
   \beas
   \hE \[ \underset{t \in [\nu, \ga_n]}{\sup} \big( e^{a t} \cY^+_t \big)^q      \]
  & \le & \hE \Big[    e^{a q  \ga_n  } \( \cY^+_{\ga_n  } \)^q    \Big]  + q
   \hE \[ \underset{t \in [0, T ]}{\sup}    \bigg| \int_t^T  \n   \b1_{\{\nu \le s \le \ga_n \}}
       \b1_{\{ \cY_s > 0 \}} e^{a q s} \(\cY^+_s\)^{q-1}  \cZ_s dB_s \bigg| \, \] \\
     &  \le & \hE \Big[   e^{a q  \ga_n  } \( \cY^+_{\ga_n  } \)^q   \Big]  + C_q
   \hE \[  \( \underset{t \in [\nu, \ga_n]}{\sup} \big( e^{a t}  \cY^+_t \big)^{q/2} \) \cd
     \( \int_\nu^{\ga_n}   \n
       \b1_{\{ \cY_s > 0 \}} e^{  a q s} \(\cY^+_s\)^{q-2} | \cZ_s |^2 d s \)^{1/2} \] \\
     &  \le & \hE \Big[   e^{a q  \ga_n  } \( \cY^+_{\ga_n  } \)^q   \Big]
       + \frac12 \hE  \[  \underset{t \in [\nu, \ga_n]}{\sup} \big( e^{a t}  \cY^+_t \big)^q    \] +
        C_q  \hE    \int_\nu^{\ga_n}  \n   \b1_{\{ \cY_s > 0 \}} e^{  a q s} \(\cY^+_s\)^{q-2} | \cZ_s |^2 d s     .
   \eeas
  As   $  \hE \[ \underset{t \in [\nu, \ga_n]}{\sup} \big( e^{a t} \cY^+_t \big)^q  \] \le
  e^{a q T} \hE \[ \underset{t \in [\nu, \tau]}{\sup} \big(   \cY^+_t \big)^q   \] < \infty $ by \eqref{eq:a041},
   it follows from \eqref{eq:a039} that
  \beas
  \hspace{-3mm}
  \hE  \n  \[ \underset{t \in [\nu, \ga_n]}{\sup} \big(   \cY^+_t \big)^q    \]
  \n   \le  \n   \hE  \n   \[ \underset{t \in [\nu, \ga_n]}{\sup} \n \big( e^{a t} \cY^+_t \big)^q    \]
   \n  \le \n   2 \hE     \Big[    e^{a q  \ga_n  } \( \cY^+_{\ga_n  } \)^q    \Big]  \n  + \n
  C_q  \hE    \int_\nu^{\ga_n}  \n   \b1_{\{ \cY_s > 0 \}} e^{  a q s} \(\cY^+_s\)^{q-2} | \cZ_s |^2 d s
   \n  \le \n   C_q \hE     \Big[   e^{a q  \ga_n  } \( \cY^+_{\ga_n  } \)^q   \Big].
  \eeas
  Because of \eqref{eq:a041} and $ \hE   \Big[ \big( \int_\nu^\tau \n |Z^i_t|^2   dt   \big)^{p/2} \Big]
  \n < \n \infty $,  $i \n = \n 1,2$,
  it holds for \pas ~ $\o  \n \in \n  \O$ that $\tau(\o)  \n = \n  \ga_{N'_\o} (\o)$ for some $N'_\o  \n \in \n  \hN$.
  Letting $n \to \infty$ in the above inequality,
  we can deduce from the monotone convergence theorem, \eqref{eq:a041} and
  dominated convergence theorem that
  \beas
  \hspace{1cm} \hE \n  \[ \underset{t \in [\nu, \tau]}{\sup} \big(   \cY^+_t \big)^q   \]
  \n   =  \n  \lmtu{n \to \infty}  \hE \n  \[ \underset{t \in [\nu, \ga_n]}{\sup} \big(   \cY^+_t \big)^q   \]
   \n  \le \n   C_q \, \lmt{n \to \infty} \hE \Big[   e^{  a q \ga_n  } \n  \( \cY^+_{ \ga_n } \)^q   \Big]
   = C_q   \hE \n \[ e^{  a q \tau  } \n  \( ( Y^1_\tau  \n - \n  Y^2_\tau)^+  \)^q   \] = 0 .  \hspace{1cm} \hb{\qed}
  \eeas

\ss \no {\bf Proof of Remark \ref{rem_g_evalu}:}
Let   $\nu \n \in \n \cT$, $\tau  \n \in \n  \cT_{\nu,T}$.
It suffices to show \eqref{eq:b641} for    $ \xi  \n \in \n  L^1(\cF_\tau) $.
Given $n \n \in \n \hN$,  we still define the stopping time $\ga_n$ as in \eqref{eq:b407}.
 As $ Y^{\tau,\xi}_{\nu \land \ga_n} \n = \n Y^{\tau,\xi}_{  \ga_n}
 \n + \n  \int_{\nu \land \ga_n}^{  \ga_n} \n \b1_{\{s \le \tau\}} g_s d  s
    \n - \n  \int_{\nu \land \ga_n}^{  \ga_n} \n  Z^{\tau,\xi}_s d B_s  $, \pas,
 similar to \eqref{eq:b417}, taking conditional expectation $\hE \big[\cd |\cF_{\nu \land \ga_n}\big]$ yields that
 $  Y^{\tau,\xi}_{\nu \land \ga_n} \n = \n
   \b1_{\{\nu \le  \ga_n\}}   \hE \big[  Y^{\tau,\xi}_{  \ga_n}
    \n + \dn  \int_{\nu \land \ga_n}^{ \tau \land \ga_n} \n      g_s d  s  \big|\cF_\nu \big]
    \n + \n  \b1_{\{\nu >  \ga_n\}}  \Big(   Y^{\tau,\xi}_{  \ga_n}
    \n + \dn  \int_{\nu \land \ga_n}^{ \tau \land  \ga_n} \n     g_s d  s  \Big)  $,  \pas ~
   Since  $\{\ga_n\}_{n \in \hN}$ is stationary,
 letting $n  \n \to \n  \infty$, we can deduce from the uniform integrability of
 $\big\{Y^{\tau,\xi}_\ga\big\}_{\ga \in \cT}$ that
 \beas
 \hspace{0.7cm}   \cE^g_{\nu,\tau}[\xi]   \n = \n
 Y^{\tau,\xi}_\nu  = \b1_{\{\nu \le  T\}} \hE \[  Y^{\tau,\xi}_T
 \n + \n  \int_\nu^\tau  \n      g_s d  s   \bigg|\cF_\nu \]
   + \b1_{\{\nu  >  T\}}   \(  Y^{\tau,\xi}_T \n + \n  \int_\nu^\tau  \n      g_s d  s \)
   = \hE \[  \xi \n + \n  \int_\nu^\tau  \n      g_s d  s   \bigg|\cF_\nu  \]   , \q \pas  \hspace{0.7cm} \hb{\qed}
 \eeas

 \subsection{Proofs of  the results in  Section \ref{sec:RBSDE} }

   \ss \no {\bf Proof of Proposition \ref{prop_monotone_result_0}:}
  {\bf (1)} {\it We first show that $\hE \n \[  (  Y_*  )^p  \] \n < \n  \infty$,    $\fa p  \n \in \n  (0,1)$.}

    As the limit of $\bF-$adapted continuous processes $Y^n$'s (thus $\bF-$predictable),
   $Y$ is also an $\bF-$predictable process.
   For any  $(t,\o) \n  \in \n  [0,T] \n \times \n \O$,
   $Y_t (\o)  \n = \n  \lmtu{n \to \infty} Y^n_t (\o) $ implies 
   $Y^+_t (\o)  \n = \n  \lmtu{n \to \infty} Y^{n,+}_t (\o) $.
   Then one can deduce that
    \bea    \label{eq:a045}
       Y^+_* (\o)  \n = \n   \underset{t \in [0,T]}{\sup} Y^+_t (\o)
        \n = \n  \underset{t \in [0,T]}{\sup} \, \underset{n \in \hN}{\sup} Y^{n,+}_t (\o)
        \n = \n   \underset{n \in \hN}{\sup} \, \underset{t \in [0,T]}{\sup}   Y^{n,+}_t (\o)
        \n = \n  \underset{n \in \hN}{\sup} \,    Y^{n,+}_* (\o)
        \n = \n  \lmtu{n \to \infty}  Y^{n,+}_* (\o)  , \q \fa \o  \n \in \n  \O .
    \eea
  For any $n  \n  \in \n   \hN$,  the continuity of process $Y^n$ shows that \pas,
  $  Y^{n,+}_*  \n   = \n    \underset{t \in [0,T]}{\sup} Y^{n,+}_t
   \n  = \n    \underset{t \in [0,T] \cap \hQ}{\sup} Y^{n,+}_t \n \in \n  \cF_T  $,
  which implies that  $Y^{n,+}_*$ is $\cF_T-$measurable.
  Then we see from \eqref{eq:a045} that    $Y^{+}_*$ is also $\cF_T-$measurable.


 \ss Let $p \n  \in \n  (\a,1)$ and
 set   $\eta \n  := \n   \xi^+  \n  + \n      \int_0^T  h_s  ds     \n  + \n    L^{+}_* \in L^1(\cF_T)$.
 Given $n \in \hN$,   we claim that   \pas
     \bea  \label{eq:a073}
   ( Y^n_t )^+    \le     C_\a \hE \[ 1 \n  + \n      \eta
  \n  + \n     \(   Y^{n,+}_*    \)^\a  \Big|\cF_t \]   , \q  t \in [0,T] ,
  \eea
  (which will be shown in the last part of this proof).
  Since $ M^n_t: = \hE \[1 \n  +\n  \eta \n  + \n   \(  Y^{n,+}_*  \)^\a  \Big| \cF_t  \] $, $t \in [0,T]$
  is a uniformly integrable martingale,
  applying Lemma 6.1 of \cite{BH_Lp_2003},  we can deduce from \eqref{eq:a073},  \eqref{eqn-d011},  \eqref{eqn-d011b},
   H\"older's inequality and   Young's inequality that
 \beas
  \hE \[    (   Y^{n,+}_*     )^p  \]
 & \tn   \le  & \tn     C^p_\a \, \hE \[ \underset{t \in [0,T]}{\sup} (M^n_t)^p \]
  \n  \le \n   \frac{C^p_\a}{1-p} \(\hE [M^n_T]\)^p
  \n  \le \n   \frac{C^p_\a}{1-p} \Big\{ 1 +
 \( \hE  [  \eta   ]  \)^p   +  \(  \hE \[  (  Y^{n,+}_*  )^\a \]  \)^p \Big\}   \\
 & \tn   \le & \tn         \frac{C^p_\a}{1-p} \Big\{ 2 +
   \hE  [  \eta   ]       +  \(  \hE \[   \(   Y^{n,+}_*    \)^p \]  \)^\a \Big\}
          \le     C_{\a,p} \big\{ 1 +  \hE  [  \eta   ]  \big\}
   + \frac12 \hE \[   \(   Y^{n,+}_*    \)^p \] .
     \eeas
  As $\hE \[    (   Y^{n,+}_*     )^p  \] \le \hE \[    (   Y^{n}_*     )^p  \] < \infty$,
  we see that   $ \hE \[    (   Y^{n,+}_*     )^p  \]
 \le   C_{\a,p}     \big\{  1 +   \hE  [  \eta    ]  \big\}  $.
 When $n \to \infty$,  \eqref{eq:a045} and the monotone convergence theorem  yield that
 $ \hE \[    (   Y^{+}_*     )^p  \]
 \le   C_{\a,p}     \big\{  1 +   \hE  [  \eta    ]  \big\} $.
   Since
 \bea   \label{eq:a077}
 |Y_t| =  Y^-_t + Y^+_t \le (Y^1_t)^-  + Y^+_t  \le |Y^1_t|  + Y^+_t , \q \fa t \in [0,T] ,
 \eea
    \eqref{eqn-d011} implies that
 $ \hE \[  (  Y_*  )^p  \] \le \hE \[    (   Y^{1}_*  )^p  \]
 + \hE \[    (   Y^{+}_*     )^p  \] < \infty $.

 Moreover, for any $\wt{p} \n \in \n  (0,\a] $, \eqref{eqn-d011b} shows that
 $ \hE \n \[  \(  Y_*  \)^{\wt{p}} \, \]
  \n \le \n  1  \n + \n   \hE \n \[  \(  Y_*  \)^{\frac{\a+1}{2}}  \]  \n < \n  \infty $.
 Hence, $\hE \n \[  (  Y_*  )^p  \] \n < \n  \infty$,    $\fa p  \n \in \n  (0,1)$.

  \ss \no {\bf (2)} {\it Next, let us show  that  $Y$ is of class  $($D$)$. }

   Since  $ \hE \[    (   Y^{+}_*     )^\a  \] \n   < \n   \infty$,
 letting $n \to \infty$ in \eqref{eq:a073},
 one can deduce   from   \eqref{eq:a045} and the monotone convergence theorem   that
 for any $t  \n  \in \n   [0,T]$,
 $   Y^+_t    \n    \le  \n      C_\a \hE \[ 1   \n   +  \n        \eta
     +   \n    (  Y^{+}_*  )^\a  \big|\cF_t \]   $, \pas ~
     Using  the continuity of   process $ Y^+ $  and process
  $\big\{\hE \[ 1 \n  + \n    \eta  \n  + \n    (  Y^{+}_*  )^\a  \big|\cF_t \]\big\}_{t \in [0,T]}$,
  we see from  \eqref{eq:a077}   that \pas
  \beas
      |Y_t| \le  |Y^1_t|  + Y^+_t    \n    \le   \n   |Y^1_t|  +  C_\a \hE \[ 1 \n  + \n      \eta
  \n  + \n   (  Y^{+}_*  )^\a  \big|\cF_t \] , \q    t  \n  \in \n   [0,T]   .
  \eeas
   This  implies that $ Y $ is of class (D) as $Y^1$ is of class (D).

  \ss \no {\bf (3)} {\it It remains to demonstrate claim \eqref{eq:a073}.}

  For any $t \n   \in  \n   [0,T]$, the continuity of process $L$ shows that \pas,
  $  \G_t  \n  := \n    \underset{s \in [0,t]}{\sup} L^+_s
   \n  = \n   \underset{s \in [0,t] \cap \hQ}{\sup} L^+_s \n \in \n  \cF_t $,
  which implies that  $\G$ is an $\bF-$adapted, continuous increasing process  with   $\hE [\G_T]  \n  < \n  \infty$.

 \ss  Let $n  \n \in \n  \hN$. Since $
     \int_0^T \n  \b1_{\{ Y^n_t > \G_t \}} \( Y^n_t  \n - \n  L_t \)^- d  t  \n = \n 0  $,
  applying It\^o$-$Tanaka's formula to process $( Y^n  \n - \n  \G )^+$ yields that
 \beas
   (Y^n_t  \n  -   \n   \G_t  )^+       \n  =  \n     (Y^n_T    \n  -  \n    \G_T )^+
     \n + \n     \int_t^T   \n    \b1_{\{Y^n_s > \G_s\}} \( g(s, Y^n_s, Z^n_s) ds    \n  -   \n     d J^n_s
    \n  -  \n      Z^n_s d B_s \)
     \n  +  \n    \int_t^T   \n    \b1_{\{Y^n_s > \G_s\}} d \, \G_s
      \n  -  \n  \frac12  (\fL^n_T   \n   -  \n   \fL^n_t)    ,  ~  t   \n \in \n    [0,T] ,
\eeas
 where $\fL^n$ is   the ``local time" of $ Y^n - \G $ at $0$.

 Set $a \n := \n  2 ( \k  \n + \n    \k^2 ) $. Given $ j \in \hN$, we define a stopping time
 $  \ga_j   \n = \n  \ga^n_j  \n := \n  \inf\{ t  \n \in \n  [0,T] \n : \int_0^t \n  |Z^n_s|^2 ds  \n > \n  j \}
 \land T  \n \in \n  \cT $, and integrate by parts the process
 $\big\{ e^{a (\ga_j \land t) }  (Y^n_{\ga_j \land t}- \G_{\ga_j \land t}  )^+ \big\}_{t \in [0,T]}$
 to obtain that \pas
\bea
  \hspace{-3mm}   e^{a (\ga_j \land t) } \n \(Y^n_{\ga_j \land t} \n   - \n   \G_{\ga_j \land t} \)^+
  \dn  + \n   a  \dn   \int_{\ga_j \land t}^{\ga_j}  \n   e^{a s} (Y^n_s  \dn  - \n   \G_s )^+ ds
  & \tn \dn  =  & \tn \dn    e^{a \ga_j   } (Y^n_{\ga_j}  \n  - \n   \G_{\ga_j} )^+
  \n  + \dn   \int_{\ga_j \land t}^{\ga_j}  \n  \b1_{\{Y^n_s > \G_s\}} e^{a s}
  \( g(s, Y^n_s, Z^n_s) ds    \n  -   \n     d J^n_s
    \n  -  \n      Z^n_s d B_s \)    \nonumber  \\
  & \tn \dn    & \tn \dn
     +     \int_{\ga_j \land t}^{\ga_j}  \n  \b1_{\{Y^n_s > \G_s\}} e^{a s} d \, \G_s
  \n - \n  \frac12  \n  \int_{\ga_j \land t}^{\ga_j}  e^{a s} d \fL^n_s , \q t \in [0,T] .  \label{eq:a051}
\eea
 Since (H4), (H5),  \eqref{eqn-d011} and  \eqref{eqn-d011b}  imply that
 \bea
 |g(t, Y^n_t, Z^n_t)| & \tn  \le  & \tn   |g(t, Y^n_t, 0)| + |g(t, Y^n_t, Z^n_t)-g(t, Y^n_t, 0)|
 \le h_t + \k |Y^n_t| + \k(h_t + |Y^n_t| + |Z^n_t|)^\a \nonumber \\
  & \tn  \le & \tn   h_t + \k |Y^n_t| + \k(h_t + |Y^n_t|)^\a + \k |Z^n_t|^\a
    \le    h_t + \k |Y^n_t|   + \k (1+ h_t + |Y^n_t|) + \k |Z^n_t|^\a  \qq  \qq  \label{eq:b123} \\
  & \tn  \le & \tn  \k +   (1+\k) h_t + 2 \k |Y^n_t-\G_t| + 2 \k  \G_t +   \k |Z^n_t|^\a  , \qq \dtp ,  \nonumber
 \eea
   taking conditional expectation $\hE[\cd|\cF_t]$ in \eqref{eq:a051},
 we can deduce  from H\"older's inequality   that for any $t \in [0,T]$
  \bea
 \hspace{-0.5cm} e^{a (\ga_j \land t) } \(Y^n_{\ga_j \land t} \n  - \n   \G_{\ga_j \land t} \)^+
  & \tn   \dn   \le   & \tn  \dn  \k T e^{a T}
  \n  +  \n   \hE \bigg[ e^{a \ga_j  } (Y^n_{\ga_j}  \n  - \n   \G_{\ga_j} )^+
   \n  + \n  (1 \n  + \n  \k) e^{a T}  \int_0^T h_s ds
   \n  + \n  ( 1 \n  + \n  2 \k T    ) e^{a T}  \G_T  \nonumber \\
 & \tn   \dn  &  \n  +    \k T^{1   -   \a/2} e^{(1-\a)aT} \( \int_{\ga_j \land t}^{\ga_j} \b1_{\{Y^n_s > \G_s\}}
 e^{2 a s} | Z^n_s |^2 ds \)^{\a/2}  \bigg| \cF_t \bigg]  , \q \pas,   \label{eq:a055}
  \eea
  where we used the fact that $ \b1_{\{Y^n_t > \G_t\}} |Y^n_t-\G_t| = (Y^n_t-\G_t)^+ $.

 Applying It\^o's formula to process  $\Big\{ e^{2 a (\ga_j \land t) }
  \(  (Y^n_{\ga_j \land t}- \G_{\ga_j \land t}  )^+ \)^2 \Big\}_{t \in [0,T]}$ in \eqref{eq:a051} yields that
 \bea
 &   & \hspace{-2cm}   e^{2a (\ga_j \land t) } \(  (Y^n_{\ga_j \land t}- \G_{\ga_j \land t}  )^+ \)^2
 + 2 a  \int_{\ga_j \land t}^{\ga_j} e^{2 a s} \( (Y^n_s - \G_s )^+ \)^2 ds
 +  \int_{\ga_j \land t}^{\ga_j} \b1_{\{Y^n_s > \G_s\}} e^{2 a s} |Z^n_s|^2 ds  \nonumber \\
 & = & \tn        e^{2 a \ga_j   } \( (Y^n_{\ga_j} \n  - \n   \G_{\ga_j} )^+ \)^2
  \n  + \n   2 \int_{\ga_j \land t}^{\ga_j} \b1_{\{Y^n_s > \G_s\}} e^{2 a s} (Y^n_s  \n  - \n   \G_s )^+
  \( g(s, Y^n_s, Z^n_s) ds    \n  -   \n     d J^n_s   \n  -  \n      Z^n_s d B_s \)
        \nonumber  \\
 & &   \n +     2 \int_{\ga_j \land t}^{\ga_j} \b1_{\{Y^n_s > \G_s\}} e^{2 a s} (Y^n_s  \n  - \n   \G_s )^+ d \G_s
 \n - \n \int_{\ga_j \land t}^{\ga_j}  e^{2 a s} (Y^n_s - \G_s )^+ d \fL^n_s
     , \q t \in [0,T] .     \label{eq:a053}
 \eea
  Since (H1) and (H4) imply that    \dtp ~
 \beas
 |g(t, Y^n_t, Z^n_t)|  \n  \le   \n     |g(t, Y^n_t, 0)|  \n  + \n   |g(t, Y^n_t, Z^n_t) \n  - \n  g(t, Y^n_t, 0)|
  \n  \le \n   h_t  \n  + \n   \k |Y^n_t|  \n  + \n   \k  |Z^n_t|
   \n  \le  \n      h_t  \n  + \n     \k |Y^n_t-\G_t|  \n  + \n     \k  \G_t  \n  +  \n    \k |Z^n_t|  , \q
 \eeas
 it holds \dtp ~ that
 \beas
  \q  \b1_{\{Y^n_t > \G_t\}} (Y^n_t \n   - \n    \G_t )^+  g(t, Y^n_t, Z^n_t)
  \n   \le \n    (Y^n_t \n   - \n    \G_t )^+ h_t
  \n   + \n    ( \k \n   + \n   2 \k^2)  \( (Y^n_t  \n   - \n    \G_t )^+ \)^2
  \n   + \n    \frac14 \b1_{\{Y^n_t > \G_t\}}   \G^2_T
  \n   + \n    \frac14 \b1_{\{Y^n_t > \G_t\}} |Z^n_t|^2  .
 \eeas
     Set   $\Psi^n_t := \underset{s \in [t,T]}{\sup} ( Y^n_s - \G_s )^+ $, $t \in [0,T]$.  It then
     follows from \eqref{eq:a053} that
 \beas
 \q  && \hspace{-1cm}  e^{2a (\ga_j \land t) }
 \( \(Y^n_{\ga_j \land t} \n  - \n   \G_{\ga_j \land t} \)^+ \)^2
  \n  + \n   \frac12  \int_{\ga_j \land t}^{\ga_j} \b1_{\{Y^n_s > \G_s\}} e^{2 a s} |Z^n_s|^2 ds   \\
  & &      \le \n    e^{2 a T  }  \( (Y^n_{\ga_j}    )^+ \)^2
  \n  + \n   2 e^{2 a T} \Psi^n_t \int_t^T      h_s    ds
       \n  + \n    \frac12 T e^{2 a T} \G^2_T
    + 2 e^{2 a T} \Psi^n_t  \G_T
   - 2 \int_{\ga_j \land t}^{\ga_j} \b1_{\{Y^n_s > \G_s\}}
  e^{2 a s}  (Y^n_s - \G_s )^+  Z^n_s d B_s \\
   & &      \le \n    e^{2 a T  } \( (Y^n_{\ga_j}    )^+ \)^2
  \n  + \n  (\Psi^n_t)^2 \n  + \n   C_0 \( \int_0^T      h_s    ds \)^2
       \n  + \n    C_0 \G^2_T
   +   \bigg| 2  \int_{\ga_j \land t}^{\ga_j} \b1_{\{Y^n_s > \G_s\}}
  e^{2 a s}  (Y^n_s - \G_s )^+  Z^n_s d B_s \bigg| , \q t \in [0,T] ,
\eeas
 Taking powers of order $\a/2$ on both sides, we see from  \eqref{eqn-d011} that
 \bea
 & \tn  & \tn  \hspace{-1.2cm} 2^{\a/2-1} e^{\a a (\ga_j \land t) } \n
 \( \(Y^n_{\ga_j \land t}- \G_{\ga_j \land t} \)^+ \)^\a
 + \frac12 \(  \int_{\ga_j \land t}^{\ga_j} \b1_{\{Y^n_s > \G_s\}} e^{2 a s} |Z^n_s|^2 ds \)^{\a/2}
 \le \n    e^{\a a T  } \( (Y^n_{\ga_j}    )^+ \)^\a
  \n  + \n  ( \Psi^n_t )^\a  \nonumber   \\
 & &      +     C_\a  \n   \( \int_0^T      h_s    ds \)^\a
 \n + \n    C_\a \G^\a_T  \n  + \n    \bigg| 2  \n  \int_t^T  \n  \b1_{\{s \le \ga_j \}}  \b1_{\{Y^n_s > \G_s\}}
  e^{2 a s}  (Y^n_s \n - \n  \G_s )^+  Z^n_s d B_s \bigg|^{\a/2} , ~ t  \n \in \n  [0,T] . \label{eq:a511}
 \eea

 Let $t \n \in \n  [0,T]$. For any $A  \n \in \n  \cF_t $, since
 \beas
 && \hspace{-1cm} \b1_A  \bigg|   \int_t^T  \n  \b1_{\{s \le \ga_j \}}  \b1_{\{Y^n_s > \G_s\}}
  e^{2 a s}  (Y^n_s \n - \n  \G_s )^+  Z^n_s d B_s  \bigg|^{\a/2}
  =   \bigg|   \int_t^T  \n \b1_A \b1_{\{s \le \ga_j \}}  \b1_{\{Y^n_s > \G_s\}}
  e^{2 a s}  (Y^n_s \n - \n  \G_s )^+  Z^n_s d B_s  \bigg|^{\a/2} \\
 && \hspace{1cm} =   \bigg|   \int_0^T  \n   \b1_A \b1_{\{t \le s \le \ga_j \}}  \b1_{\{Y^n_s > \G_s\}}
  e^{2 a s}  (Y^n_s \n - \n  \G_s )^+  Z^n_s d B_s  \bigg|^{\a/2} ,
 \eeas
  multiplying $\b1_A$ to \eqref{eq:a511} and   taking   expectation,
 we can deduce from   the Burkholder-Davis-Gundy inequality and  \eqref{eqn-d011b}
 \beas
 & & \hspace{-0.8cm} \frac12 \hE \n  \[ \b1_A
 \(  \int_{\ga_j \land t}^{\ga_j} \b1_{\{Y^n_s > \G_s\}} e^{2 a s} |Z^n_s|^2 ds \)^{\a/2}   \]
 \le \n   C_\a \hE   \Bigg[   \b1_A \n  \( (Y^n_{\ga_j}    )^+ \)^\a
  \n  + \n   \b1_A   ( \Psi^n_t )^\a \n  + \n   \b1_A  \n \( \int_0^T  h_s  ds \)^\a   \n  + \n   \b1_A   \G^\a_T \\
 &  &   \q   +      \bigg(   \int_0^T  \n   \b1_A \b1_{\{t \le s \le \ga_j \}} \b1_{\{Y^n_s > \G_s\}}
  e^{4 a s} \( (Y^n_s \n  - \n   \G_s )^+ \)^2 | Z^n_s |^2 d s \bigg)^{\a/4} \,  \Bigg] \\
  &  &  \le  \n    C_\a \hE \n \[ \b1_A  \n  + \n   \b1_A    ( Y^n_{\ga_j} )^+
  \n  + \n  \b1_A  ( \Psi^n_t )^\a \n  + \n   \b1_A  \n  \int_0^T  h_s  ds     \n  + \n   \b1_A   \G_T
    \n  + \n      \b1_A  ( \Psi^n_t )^{\a/2} \n \cd \n \(   \int_t^T  \n     \b1_{\{  s \le \ga_j \}}
     \b1_{\{Y^n_s > \G_s\}} e^{2 a s}  | Z^n_s |^2 d s \)^{\a/4}   \] \\
  &  &  \le  \n    C_\a \hE \n \[ \b1_A  \n  + \n   \b1_A    ( Y^n_{\ga_j} )^+
  \n  + \n  \b1_A  ( \Psi^n_t )^\a \n  + \n   \b1_A  \n  \int_0^T  h_s  ds     \n  + \n   \b1_A   \G_T   \]
    \n  + \n    \frac14  \hE \n \[  \b1_A  \n  \(   \int_{\ga_j \land t}^{\ga_j}   \n
     \b1_{\{Y^n_s > \G_s\}} e^{2 a s}  | Z^n_s |^2 d s \)^{\a/2}   \] .
 \eeas
 Since $ \hE \n  \[
 \(  \int_{\ga_j \land t}^{\ga_j} \b1_{\{Y^n_s > \G_s\}} e^{2 a s} |Z^n_s|^2 ds \)^{\a/2}   \]
 \le e^{\a a T} j^{\a/2} $ and since
 \beas
  \hE \[  ( \Psi^n_0 )^\a  \]
  = \hE \[  \underset{t \in [0,T]}{\sup}  \( ( Y^n_t \n  - \n  \G_t )^+ \)^\a  \]
  \le \hE \[  \underset{t \in [0,T]}{\sup} \( ( Y^n_t )^+ \)^\a  \] \le \| Y^n \|_{\hS^\a} < \infty ,
 \eeas
  letting $A $ vary over $\cF_t$ yields that
 \beas
   \hE \n  \[   \(  \int_{\ga_j \land t}^{\ga_j} \b1_{\{Y^n_s > \G_s\}} e^{2 a s} |Z^n_s|^2 ds \)^{\a/2} \Bigg| \cF_t  \]
 \le  C_\a \hE \n \[ 1  \n  + \n        ( Y^n_{\ga_j} )^+
  \n  + \n     ( \Psi^n_t )^\a \n  + \n        \int_0^T  h_s  ds     \n  + \n    \G_T  \bigg| \cF_t   \] .
 \eeas
 Then we see from   \eqref{eq:a055}   that
  \bea  \label{eq:a071}
  \(Y^n_{\ga_j \land t} \n  - \n   \G_{\ga_j \land t} \)^+ \n  \le  \n
  e^{a (\ga_j \land t) } \(Y^n_{\ga_j \land t} \n  - \n   \G_{\ga_j \land t} \)^+
  \n   \le     \n     C_\a \hE \[ 1 \n  + \n    (Y^n_{\ga_j}     )^+
  \n  + \n   ( \Psi^n_0 )^\a \n  + \n      \int_0^T  h_s  ds     \n  + \n    \G_T \bigg|\cF_t \]    , \q \pas
  \eea

  \ss The uniform integrability of $\{Y^n_\ga\}_{\ga \in \cT}$   implies that of
  $ \big\{ (Y^n_\ga     )^+ \big\}_{\ga \in \cT}$.
   As  $Z^n \in \underset{p \in (0,1)}{\cap} \hH^{2,p} \subset  \hH^{2,0}$, $\{\ga_j\}_{j \in \hN}$ is stationary.
   So letting $j \to \infty$ in \eqref{eq:a071},  one can deduce from  the continuity of process $Y^n$   that
   \beas
  (Y^n_t)^+ \n  \le \n   \G_t  \n  + \n   \(Y^n_t \n  - \n   \G_t \)^+
  \n   \le     \n   \G_t  \n  + \n    C_\a \hE \[ 1 \n  + \n      \xi^+
  \n  + \n   ( \Psi^n_0 )^\a \n  + \n      \int_0^T  h_s  ds     \n  + \n    \G_T \bigg|\cF_t \]
   \n  \le \n   C_\a \hE \[ 1 \n  + \n      \eta  \n  + \n   ( \Psi^n_0 )^\a   \big|\cF_t \]
       , \q \pas
  \eeas
  Then claim \eqref{eq:a073} follows from the continuity of process $ Y^{n, +} $  and of process
  $\big\{\hE \[ 1 \n  + \n    \eta  \n  + \n   ( \Psi^n_0 )^\a  \big|\cF_t \]\big\}_{t \in [0,T]}$.   \qed

     \ss \no {\bf Proof of Proposition \ref{prop_monotone_result}:} The proof is relatively lengthy,
     see our introduction for a sketch.
     We will defer the demonstration of some technicalities  (those equations with starred labels) to the appendix.

   \ss \no {\bf (1)} For any $n \in \hN$,
    $K^n_t :=  n  \int_0^t \b1_{\{\nu < s \le \tau\}} ( Y^n_s-L_s  )^- ds $,
    $t \in [0,T]$   is clearly  a process of $\hK^0$ satisfying
    \bea  \label{eq:a245}
    K^n_t \n = \n  0 , \q \fa t  \n \in \n  [0,\nu] .
    \eea
  As    $    K^n_{\tau }  \n - \n  K^n_t
     \n = \n  n \n \int_t^{\tau } \b1_{\{\nu < s \le \tau\}} ( Y^n_s \n - \n L_s  )^- ds
     \n = \n  n \int_t^{\tau }   ( Y^n_s \n - \n L_s  )^- ds $, $  \fa  t  \n \in \n  [\nu,\tau ]  $,
      \eqref{eq:a243} shows that \pas
     \bea  \label{eq:a333a}
    Y^n_t    \n  =  \n Y^n_{\tau} \n + \n  \int_t^{\tau}  g(s, Y^n_s, Z^n_s) ds
    + K^n_{\tau} - K^n_t
  - \int_t^{\tau} Z^n_s dB_s, \q  \fa  t \in [\nu,\tau] .
    \eea

       Since $\big\{\b1_{\{t \ge \nu \}} \big\}_{t \in [0,T]}$ is an $\bF-$adapted c\`adl\`ag process
   and since   each   $\big\{ Y^n_{ \tau \land t } \big\}_{t \in [0,T]}$ is an $\bF-$adapted continuous process,
   we see that $Y$ is an $\bF-$optional process.
   (It takes some  effort to show the continuity of $Y$ between $\nu$ and $T$,
   see \eqref{eq:a205} for an intermediate result.)
   By    the Debut theorem,
   \bea  \label{eq:a273}
     \tau_\ell   \n  :=   \n
      \inf \bigg\{ t \in [\nu,\tau]:     (Y^1_t)^-  \n + \n  Y^+_t  \n + \n  L^+_t
    \n + \n  \int_\nu^t     h_s ds  \n > \n  \ell  \bigg\}    \n \land \n   \tau, \q \ell \in \hN
    \eea
  are stopping times with $\nu \le \tau_\ell \le \tau$, i.e. $\tau_\ell \in \cT_{\nu,\tau}$.
  As   $ \hE \n \[     L^{+}_*  \dn + \dn \int_0^T \n h_t dt \] \dn < \n \infty    $  and $\hP \Big\{
   \underset{t \in [\nu,\tau]}{\sup} \((Y^1_t)^-  \n + \n  Y^+_t\)  \n < \n  \infty \Big\} \n = \n 1$,
  it holds for  any $\o  \n \in \n  \O$ except on a   $\hP-$null set $ \cN_1 $ that
  \beas    
      \tau(\o) = \tau_{N_\o} (\o) \hb{ for some }  N_\o  \n \in \n  \hN .
  \eeas

 Now,    let us fix  $\ell \n \in \n  \hN$ for this part as well as next two parts.
 Let $\cN_2  \n := \n  \underset{n \in \hN}{\cup}
 \{\o  \n \in \n  \O \n : \hb{ the path $ Y^n_\cd (\o) $ is not continuous} \}$
 (which is clearly a $\hP-$null set) and set  $A_\ell  \n := \n  \{\nu  \n < \n  \tau_\ell\} \cap \cN^c_2
    \n \in \n  \cF_{\nu \land \tau_\ell}   \n \subset \n  \cF_\nu $.
 Given $\o \in A_\ell$, for any $n \in \hN$
   we can deduce from \eqref{eq:a273} that
    $  |Y^n_t(\o)|  \n   \le \n  \ell $, $\fa t \in \big[\nu(\o), \tau_\ell(\o)\big) $,
    and the continuity of each $Y^n$ implies that $  |Y^n_t(\o)|  \n   \le \n  \ell $, $\fa t \in \big[\nu(\o), \tau_\ell(\o)\big] $. Then it follows from the monotonicity of
   $\{Y^n\}_{n \in \hN}$ that
  \bea  \label{eq:a095}
  \underset{n \in \hN}{\sup} \, \big|Y^n_t (\o)\big| \le \big|Y^1_t(\o)\big| \vee \big|Y_t(\o)\big| \le \ell ,
  \q \fa t \in \big[\nu(\o), \tau_\ell(\o)\big], \; \fa \o \in A_\ell .
  \eea

  Let $n \n \in \n  \hN$.  As $ \hE \big[ | \b1_{A_\ell} Y^n_\nu| \big]  \n \le \n  \ell  $, 
   Corollary \ref{cor_martingale} shows that
 there exists a unique $\wt{Z}^{\ell,n}  \n \in \n  \underset{p \in (0,1)}{\cap}   \hH^{2,p} $
 such that $\hP \big\{ \hE \big[\b1_{A_\ell} Y^n_\nu|\cF_t \big]
 \n = \n  \hE \big[\b1_{A_\ell} Y^n_\nu \big]  \n +  \n \int_0^t \n  \wt{Z}^{\ell,n}_s dB_s,
 ~ \fa t  \n \in \n  [0,T] \big\} =1 $.
 Similar to \eqref{eq:b241},  we can deduce from    \eqref{eq:a333a}   that \pas
 \bea
 Y^{\ell,n}_t & \dn \dn: =& \dn \dn \hE \big[ \b1_{A_\ell} Y^n_\nu |\cF_{\nu \land t}\big]
 \n  + \n  Y^n_{\nu \vee (\tau_\ell \land t)}  \n - \n  Y^n_\nu
 \n = \n \hE \big[\b1_{A_\ell} Y^n_\nu \big]
 \n - \n  \int_0^t \n  \b1_{\{\nu < s \le \tau_\ell\}}     g(s, Y^n_s, Z^n_s) ds
 \n - \n  K^n_{\nu \vee (\tau_\ell \land t)}  \n + \n  K^n_\nu \nonumber \\
  & \dn \dn    & \dn \dn
  \n + \n  \int_0^t  \n  \( \b1_{\{s \le \nu \}} \wt{Z}^{\ell,n}_s
  \n + \n \b1_{\{\nu < s \le \tau_\ell\}}  Z^n_s \) d B_s   , \q t \in [0,T] . \label{eq:a519}
 \eea
 Thus $Y^{\ell,n}$ is an $\bF-$adapted continuous process (i.e. $Y^{\ell,n} \in \hS^0$) that satisfies
   \bea  \label{eq:a281}
     Y^{\ell,n}_t = \hE \big[ \b1_{A_\ell} Y^n_\nu |\cF_\nu  \big]  \n  + \n  Y^n_t  \n - \n  Y^n_\nu
   = \b1_{A_\ell} Y^n_\nu \n  + \n \b1_{A_\ell} ( Y^n_t  \n - \n  Y^n_\nu )
   = \b1_{A_\ell} Y^n_t  , \q \fa t \in [\nu, \tau_\ell] ,
   \eea
 which together with   \eqref{eq:a519} shows that    \pas ~
\bea  \label{eq:a333}
      \hspace{-0.5cm}  Y^{\ell,n}_t \dn - \n  Y^{\ell,n}_{\tau_\ell}
  \n - \n  K^n_{\tau_\ell} \dn  + \n  K^n_t   \dn  + \dn  \int_t^{\tau_\ell}  \dn  Z^n_s dB_s
        \n =   \dn             \int_t^{\tau_\ell}  \dn  g(s, Y^n_s, Z^n_s) ds
       \n =  \n  \b1_{A_\ell}  \dn  \int_t^{\tau_\ell}  \dn  g(s, Y^{\ell,n}_s, Z^n_s) ds
     \n =  \dn    \int_t^{\tau_\ell}  \dn  g(s, Y^{\ell,n}_s, Z^n_s) ds ,
     ~  \fa  t  \n \in \n  [\nu,\tau_\ell] . ~ \;
    \eea

  Since     $\hE \big[ \big|  Y^{\ell,n}_\nu\big|\big] 
   \n \le \n  \ell$
  by \eqref{eq:a281}, \eqref{eq:a095} and since   $ K^n_\nu  \n = \n  0$ by \eqref{eq:a245},
  applying Lemma \ref{lem_RBSDE_estimate}  with $(Y,Z,K)  \n = \n  (Y^{\ell,n} , Z^n , K^n)$
  and $( \tau, p)  \n = \n  (  \tau_\ell, 2)$, we see from \eqref{eq:a281}, \eqref{eq:a095} and  \eqref{eq:a273}  that
      \bea \label{eq:a093}
  \hspace{-0.5cm}   \hE \n \int_\nu^{\tau_\ell  } \n |Z^n_t|^2 dt   \n + \n   \hE \n  \[  ( K^n_{\tau_\ell }   )^2  \]
     \n \le \n  C_0   \hE  \n  \[ \underset{t \in [\nu, \tau_\ell ]}{\sup} \big| Y^{\ell,n}_t \big|^2 \]
      \n  + \n  C_0   \hE  \n  \[  \(   \int_\nu^{\tau_\ell}  \n  h_t dt  \)^2 \, \]
    \n \le \n  C_0  \ell^2         .
    \eea
   It then follows from (H1) that
 $  \hE \n  \int_\nu^{\tau_\ell  } \n  |g(t,Y^{\ell,n}_t,Z^n_t)  \n - \n g(t,Y^{\ell,n}_t,0) |^2 dt
 \n \le \n  \k^2 \, \hE \n  \int_\nu^{\tau_\ell  } \n  |Z^n_t|^2 dt
   \n \le \n  C_0  \ell^2    $.
  In virtue of Theorem 5.2.1 of \cite{FA_Yosida},
  $ \big\{ \b1_{\{\nu < t \le \tau_\ell\}} Z^n_t \big\}_{t \in [0,T]}$, $n    \n \in \n     \hN$
  has a weakly convergent subsequence
  \big(we   still denote it by $\{\b1_{\{\nu < t \le \tau_\ell\}} Z^n_t\}_{t \in [0,T]}$, $n  \n  \in  \n  \hN$\big)
  with   limit $\cZ^\ell \n  \in \n  \hH^{2, 2} $;
  and  $\big\{\b1_{\{\nu < t \le \tau_\ell\}} (g(t,Y^{\ell,n}_t,Z^n_t)
    \n  -  \n   g(t,Y^{\ell,n}_t,0)) \big\}_{t \in [0,T]}$,
  $n \n  \in \n   \hN$
  has a weakly convergent subsequence   \big(we   still denote it by
   $ \big\{\b1_{\{\nu < t \le \tau_\ell\}} (g(t,Y^{\ell,n}_t,Z^n_t)
    \n  -  \n   g(t,Y^{\ell,n}_t,0)) \big\}_{t \in [0,T]}$,
  $n   \n   \in   \n   \hN$\big)   with   limit $\wt{h}^\ell  \n \in \n  \hH^{2, 2} $.
  It is easy to deduce    that
  \bea \label{eq:a353}
   \cZ^\ell_t \n = \n  \b1_{\{\nu < t \le \tau_\ell\}} \cZ^\ell_t  \q
  \hb{and} \q  \wt{h}^\ell_t  \n = \n  \b1_{\{\nu < t \le \tau_\ell\}}  \wt{h}^\ell_t  , \q  \dtp ~
  \eea

 \ss  The $\bF-$optional measurability   of $Y$ implies that of stopped processes
    $\big\{Y_{\nu \land t}\big\}_{t \in [0,T]} $
   and  $\big\{Y_{\tau_\ell \land t}\big\}_{t \in [0,T]} $
   (see e.g. Corollary 3.24 of \cite{HWY_1992}).
   As $A_\ell \cap \{t >  \nu\} 
    \in \cF_t$ for any $t \in [0,T]$,
     $\big\{ \b1_{A_\ell \cap \{t \ge  \nu\}} \big\}_{t \in [0,T]}$ is an $\bF-$adapted c\`adl\`ag process.   Then
   \bea  \label{eq:a527}
   Y_{\nu \vee (\tau_\ell \land t)} - Y_\nu = \b1_{A_\ell \cap \{t \ge  \nu\}} \( Y_{\tau_\ell \land t} - Y_\nu     \)
   = \b1_{A_\ell \cap \{t \ge  \nu\}} \( Y_{\tau_\ell \land t} - Y_{\nu \land t}    \) ,   \q t \in [0,T]
   \eea
   is an $\bF-$optional process  and it follows that
     \bea   \label{eq:a133}
     \wt{K}^\ell_t   \n  :=   \n    Y_\nu \n - \n  Y_{\nu \vee (\tau_\ell \land t)}
     \n   - \n     \int_0^t \n  \b1_{\{ \nu < s \le  \tau_\ell   \}} \(  g(s,  Y_s, 0) \n  + \n  \wt{h}^\ell_s \) ds
    \n  + \n   \int_0^t \n  \b1_{\{ \nu < s \le  \tau_\ell   \}}  \cZ^\ell_s d B_s
  , \q  t \in [0,T] \qq
   \eea
   also defines an $\bF-$optional process.
   Since \eqref{eq:a527},  \eqref{eq:a095},  (H4), \eqref{eq:a273}  and H\"older's inequality imply  that
 \beas
       \big| \wt{K}^\ell_t \big|  &  \tn   \le   &  \tn
    \b1_{A_\ell \cap \{t \ge  \nu\}} \( \big| Y_{\tau_\ell \land t} \big| +  | Y_\nu  |    \)
    \n  + \n  \b1_{A_\ell} \n  \int_\nu^{\tau_\ell  }
   \dn  \( h_t  \n  + \n   \k |  Y_t |  \n  + \n  | \wt{h}^\ell_t | \) dt
    \n  + \n   \bigg| \int_0^t \n \b1_{\{\nu < s \le  \tau_\ell   \}}  \n
     \cZ^\ell_s d B_s \bigg|  \\
       &  \tn   \le   &  \tn   3 \ell  \n  + \n   \k \ell T
      \n  + \n     \( T \int_\nu^{\tau_\ell   }  \n   | \wt{h}^\ell_t |^2 dt  \)^{\dn 1/2}
     \n  + \n  \underset{t \in [0,T]}{\sup} \bigg|  \int_0^t \n \b1_{\{\nu < s \le  \tau_\ell   \}}
           \cZ^\ell_s d B_s \bigg|  , \q \fa t \in [0,T] ,
 \eeas
  Doob's martingale inequality and \eqref{eqn-d011} show  that
  \bea   \label{eq:a135}
  \hspace{-5mm} \hE \n \[ \big( \wt{K}^{\ell}_* \big)^2 \]  \n \le \n  C_0 \ell^2
    \n + \n  3 T \hE  \dn \int_\nu^{\tau_\ell}   \n  | \wt{h}^\ell_t |^2   dt
    \n + \n  3 \hE  \n \[  \underset{t \in [0,T]}{\sup} \bigg|  \int_0^t \n \b1_{\{\nu < s \le  \tau_\ell   \}}
           \cZ^\ell_s d B_s \bigg|^2 \]
    \n \le \n  C_0 \ell^2    \n + \n  C_0 \hE \int_\nu^{\tau_\ell}  \n
    \Big( \big| \wt{h}^\ell_t \big|^2  \n + \n  \big| \cZ^\ell_t \big|^2 \Big) dt  \n < \n  \infty .
  \eea

 We next claim that
   \bhe \bea
     \wt{K}^\ell \hb{ satisfies the conditions of Lemma \ref{lem_increasing}
  and is thus an increasing process. } \label{eq:a187}
  \eea \ehe
  As   $ \hE \Big[ \big( \wt{K}^\ell_T \big)^2 \Big] \n < \n \infty $ by \eqref{eq:a135}, it holds \pas ~ that
 $ \wt{K}^\ell_T  \n < \n  \infty $.
 Then applying Lemma 2.2 of \cite{Peng_1999} to  \eqref{eq:a133} shows that
 both process $\wt{K}^\ell$ and  process  $\big\{Y_{\nu \vee (\tau_\ell \land t)}   \big\}_{t \in [0,T]}$
 have \pas ~ c\`adl\`ag  paths.

 \ss  \no  {\bf (2)} By  H\"older's inequality and \eqref{eq:a093},
 $    \hE  \n  \int_\nu^{\tau_\ell}       (Y^n_t \n - \n  L_t)^- dt
 \n = \n  \frac{1}{n} \hE \big[    K^n_{\tau_\ell }       \big]
  \n \le \n  \frac{1}{n} \big\{  \hE \big[  ( K^n_{\tau_\ell }    )^2  \big] \big\}^{1/2}
    \n \le \n  \frac{1}{n}   C_0 \ell     $, $ \fa   n  \n \in \n  \hN $.
     Letting $n  \n \to \n  \infty$, we know from   the monotone convergence theorem 
        that
    \beas
      \hE  \n  \int_\nu^{\tau_\ell}  \n  (Y_t  \n  - \n  L_t)^- dt
    \n = \n  \lmtd{n \to \infty} \hE  \n  \int_\nu^{\tau_\ell}  \n  (Y^n_t  \n - \n  L_t)^- dt  \n = \n  0 \, ,
    \eeas
  so it holds \dtp ~ that
  $ \b1_{\{\nu < t < \tau_\ell\}} \( Y_t  \n - \n  L_t \)^-  \n = \n 0 $.
   Since   $\big\{Y_{\nu \vee (\tau_\ell \land t)}\big\}_{t \in [0,T]}$ has \pas ~ c\`adl\`ag  paths by part (1)
  and $L$ has \pas ~continuous paths, one can deduce that for
  any    $\o  \n \in \n  A_\ell $ except a $\hP-$null set $\wt{\cN}_\ell$,
  $   Y_t (\o) \n  \ge \n  L_t (\o)  $  for any   $  t  \n \in \n   \big[\nu(\o),\tau_\ell(\o)  \big)    $.
  Given $\o \in  \{\nu < \tau \} \cap \cN^c_1 \cap \cN^c_2 \cap \( \underset{\ell \in \hN}{\cup} \wt{\cN}_\ell \)^c $,
  there exists an $n_\o \in \hN$ such that $   \tau_{n_\o} (\o) = \tau (\o) > \nu (\o) $. So
  $ \o \in A_{n_\o} \cap \wt{\cN}^c_{n_\o}
  = \big\{\o' \in \O: \nu (\o') <   \tau_{n_\o} (\o')\big\} \cap \cN^c_2 \cap \wt{\cN}^c_{n_\o} $
  and   $   Y_t (\o) \n  \ge \n  L_t (\o)  $ holds
  for any   $  t  \n \in \n   \big[\nu(\o),\tau_{n_\o} (\o)  \big)  = \big[\nu(\o),\tau  (\o)  \big)  $.
  In summary, it holds  for \pas ~ $\o \in \{\nu < \tau \}$ that $   Y_t (\o) \n  \ge \n  L_t (\o)  $
  for any   $  t  \n \in \n \big[\nu(\o),\tau  (\o)  \big) $, which together with
  $\hP \{ Y_\tau \n \ge \n L_\tau \} \n = \n 1$ shows that for any $\o \in \{\nu < \tau \} $
  except on a $\hP-$null set $\wh{\cN}$
  \bea  \label{eq:a181}
  Y_t (\o)  \ge  L_t(\o) , \q \fa t \in \big[\nu(\o),\tau(\o)\big]  .
  \eea

  Now we freeze the parameter $\ell$ again and let $ \o \in A_\ell \cap \wh{\cN}^c   $.
  As $A_\ell \subset \{ \nu < \tau \} \cap \cN^c_2$,
  we see from \eqref{eq:a181} that
  $Y_t (\o)  \ge  L_t(\o) $ for any $ t \in \big[\nu(\o),\tau_\ell (\o)\big]  $.
  Since    continuous function
   $    \( Y^n_t  \n - \n  L_t \)^- (\o) $,
   $ t  \n \in \n  \big[\nu(\o),\tau_\ell (\o)\big] $ is decreasing  to
   $ \( Y_t  \n - \n  L_t \)^-  (\o)  \n  = \n  0 $,
   $t  \n \in \n  \big[\nu(\o),\tau_\ell (\o)\big]$ when $ n  \n \to \n  \infty $,
     Dini's theorem shows that
   \beas
       \lmtd{n \to \infty} \underset{t \in [\nu(\o),\tau_\ell (\o)]}{\sup}
   \( Y^n_t  \n - \n  L_t \)^- (\o)    =    0  .
   \eeas
    As $ \b1_{A_\ell} \underset{t \in [\nu,\tau_\ell]}{\sup}
   \( Y^n_t  \n - \n  L_t \)^-  \n \le \n
    \b1_{A_\ell} \, \underset{t \in [\nu,\tau_\ell]}{\sup} \( L^+_t
    \n + \n  |Y^n_t|\)  \n \le \n  2 \ell $,
   $\fa n  \n \in \n  \hN$ by \eqref{eq:a273}, \eqref{eq:a281} and  \eqref{eq:a095},
      an application of  the  bounded convergence theorem yields  that
   \bea  \label{eq:a151}
   \lmtd{n \to \infty} \hE \[ \b1_{A_\ell} \, \underset{t \in [\nu,\tau_\ell]}{\sup}
   \(  ( Y^n_t - L_t  )^- \)^2 \] = 0 .
   \eea

  Similar to the arguments used in \cite{EKPPQ-1997} (see pages 21-22 therein),    we can
  deduce from \eqref{eq:a151}  that
  \bhe  \bea \label{eq:a189}
  \hb{ $ \Big\{ Y^{\ell,n}  \Big\}_{n  \in    \hN} $    is a  Cauchy sequence in
  $\hS^2$ and $ \big\{ \b1_{\{\nu < t \le \tau_\ell  \}} Z^n_t \big\}_{t \in [0,T]} $,
  $n \n  \in \n  \hN$ is a  Cauchy sequence in   $ \hH^{2,2} $. } \q
  \eea  \ehe
   Let $\cY^\ell   \n \in \n \hS^2$ and $\wt{\cZ}^\ell \n \in \n \hH^{2,2}$ be their limits respectively, i.e.
  \bea  \label{eq:a161}
  \lmtd{n \to \infty}  \hE \[ \underset{t \in [0,T]}{\sup} \big| Y^{\ell,n}_t
  \n  - \n   \cY^\ell_t  \big|^2 \] \n  + \n
  \lmt{n \to \infty} \hE  \int_0^T  \n    \big| \b1_{\{\nu < t \le \tau_\ell \}} Z^n_t
  \n  - \n  \wt{\cZ}^\ell_t  \big|^2 dt
   \n   = \n   0 .
  \eea
  Up to a subsequence of $\Big\{Y^{\ell,n}\Big\}_{n \in \hN}$,   
  one has $\lmtd{n \to \infty}    \underset{t \in [0,T]}{\sup} \big| Y^{\ell,n}_t
  \n  - \n   \cY^\ell_t  \big| = 0$, \pas ~ It follows from   \eqref{eq:a281}   that \pas
  \bea \label{eq:a165}
   \cY^\ell_t = \lmtu{n \to \infty} Y^{\ell,n}_t
   = \lmtu{n \to \infty} \b1_{A_\ell} Y^n_t = \b1_{A_\ell} Y_t , \q  \fa  t \in [\nu,\tau_\ell] ,
  \eea
    which together with the continuity of $\cY^\ell$ shows that
  \bea   \label{eq:a205}
  \big\{  \b1_{A_\ell}  Y_{\nu \vee (\tau_\ell \land t)} \big\}_{t \in [0,T]} \hb{ is a continuous process.   }
  \eea
  On the other hand, the strong limit $\wt{\cZ}^\ell$ and the weak limit $\cZ^\ell$   of
  $ \big\{ \b1_{\{\nu < t \le \tau_\ell  \}} Z^n_t \big\}_{t \in [0,T]} $,
  $n \n  \in \n  \hN$  must coincide, i.e.
  $  \wt{\cZ}^\ell_t =    \cZ^\ell_t  $,   \dtp ,
    which together with  \eqref{eq:a161},  \eqref{eq:a281} and \eqref{eq:a165} and \eqref{eq:a353} shows   that
   \bea    \label{eq:a341}
      \lmt{n \to \infty}  \hE \[ \b1_{A_\ell} \underset{t \in [\nu,\tau_\ell]}{\sup}  \big|  Y^n_t
       \n - \n   Y_t \big|^2  \]
  \n   + \n   \lmt{n \to \infty}   \hE \n  \int_\nu^{\tau_\ell    } \n
   |Z^n_t \n  - \n  \cZ^\ell_t|^2 dt = 0  .
   \eea

  \ss \no {\bf (3)} By \eqref{eq:a527} and \eqref{eq:a205},
   $       Y_\nu \n- \n  Y_{\nu \vee (\tau_\ell \land t)}  \n = \n
      \b1_{A_\ell} \n \(  Y_\nu \n- \n  Y_{\nu \vee (\tau_\ell \land t)}  \) $, $t  \n \in \n  [0,T]$
      is an $\bF-$adapted continuous process, then so is
   \bea    \label{eq:xax014}
   \cK^\ell_t :=  Y_\nu \n - \n Y_{\nu \vee (\tau_\ell \land t)}
  \n   - \n   \int_0^t   \b1_{\{ \nu < s \le  \tau_\ell   \}}       g(s,  Y_s, \cZ^\ell_s) ds
  \n   + \n   \int_0^t   \b1_{\{ \nu < s \le  \tau_\ell   \}}   \cZ^\ell_s  d B_s    , \q
   t  \n \in \n  [0,T] .
  \eea
     One     can     deduce from \eqref{eq:a341}       that
   \bhe  \bea   \label{eq:a179}
   \lmt{n \to \infty}  \hE \[ \underset{t \in [\nu,\tau_\ell]}{\sup} \big| K^n_t
   - \cK^\ell_t \big|^2 \]   = 0 .
   \eea  \ehe
  So up to a subsequence of $\{K^n\}_{n \in \hN}$,  
  it holds \pas ~ that
  \bea \label{eq:a195}
     \lmt{n \to \infty} \, \underset{t \in [\nu,\tau_\ell]}{\sup} \big| K^n_t   - \cK^\ell_t \big| = 0 \q
    \hb{and thus} \q    \cK^\ell_t = \lmt{n \to \infty} K^n_t   , \q  \fa   t \in [\nu,\tau_\ell] ,
  \eea
  which   together with the monotonicity of $K^n$'s show that for \pas ~ $\o \in \O$, the
  path $ \cK^\ell_\cd (\o)   $ is increasing over  period $[\nu(\o),\tau_\ell(\o)]$.
    \if{0}
     {\it Define $g_\ell (t,\o,y,z) \n : = \n \b1_{\{t \le \tau_\ell (\o) \}} g(t,\o,y,z) $,
    $\fa (t,\o,y,z) \n \in \n [0,T] \n  \times \n   \O  \n  \times \n   \hR  \n  \times \n   \hR^d $.
     We now show that the processes $ \big\{ \(Y_{\tau_\ell \land t},
     \cZ^\ell_t, \cK^\ell_{\tau_\ell \land t} \) \big\}_{t \in [0,T]}  $  solves
    RBSDE$\(Y_{\tau_\ell}, g_\ell, \{ L_{\tau_\ell \land t} \}_{t \in [0,T]}  \)$, i.e.
    it holds \pas~ that
   \bea  \label{eq:a185}
   \begin{cases}
   \dis      L_{\tau_\ell \land t} \le  Y_{\tau_\ell \land t}
   =   Y_{\tau_\ell} + \int_{\tau_\ell \land t}^{\tau_\ell} g(s,Y_s, \cZ^\ell_s   )  ds
   + \cK^\ell_{\tau_\ell} - \cK^\ell_{\tau_\ell \land t}
   - \int_{\tau_\ell \land t}^{\tau_\ell} \cZ^\ell_s d B_s    , \q    t \in [0,T]  ,          \vspace{1mm}     \\
   \dis \int_0^{\tau_\ell}  (Y_t-L_t) d \cK^\ell_t = 0 .
   \end{cases}
    \eea  }
    \fi
  One can also deduce from \eqref{eq:a195} that  for \pas ~ $\o \in \O$,
  the  measure $d K^n_t (\o) $ converges weakly to
  the  measure $d \cK^\ell_t (\o) $ on    period  $   \[ \nu(\o) ,   \tau_\ell(\o) \] $. It then follows that \pas
   \bhe \bea \label{eq:a199}
    \int_t^{\tau_\ell}  (Y_s-L_s) d \cK^\ell_s = 0 ,  \q \fa t \in [\nu, \tau_\ell] .
   \eea \ehe

    \ss \no {\bf (4)} {\it Setting $\tau_0 := \nu $,   we next show that process $Y$ together with processes
 \bea \label{eq:a171}
    Z_t \n : = \n \sum_{\ell \in \hN} \b1_{\{ \tau_{\ell-1} < t \le \tau_\ell \}} \cZ^\ell_t
 \q \hb{and} \q     K_t  \n : =  \n   \sum_{\ell \in \hN}
 \big( \cK^\ell_{\tau_\ell \land t} \n  - \n   \cK^\ell_{\tau_{\ell-1} \land t} \big) , \q  t  \n \in \n  [0, T]
 \eea
 solves \eqref{RBSDEL_tau}.  }

 \ss   As  $ \big\{ \b1_{\{ \tau_{\ell-1} < t \le \tau_\ell \}} \big\}_{t \in [0,T]}$
 is an $\bF-$adapted c\`agl\`ad  process (thus    $\bF-$predictable) for each $\ell \in \hN$,
 the process $Z$ is   $\bF-$predictable.
 On the other hand, it is clear that   $K$ is an $\bF-$adapted process with $K_0 = 0$.

 Let  $\cN_3$ be the $\hP-$null set such that   for any $\o \n \in \n   \cN^c_3 $
  and $\ell  \n \in \n  \hN$,
  \beas
  \hb{ $\int_0^T \n |\cZ^\ell_t (\o)|^2 dt \n < \n \infty $
   and    the path $ \big\{\cK^\ell_t (\o) \big\}_{t \in [0,T]}$  is   continuous
   and increasing over period $[\nu(\o),\tau_\ell(\o)]$. }
  \eeas
 Given  $\o \in \(\cN_1 \cup \cN_3 \)^c $,   both sums in \eqref{eq:a171} are finite sums:
    \bea  \label{eq:a311}
   ~ \;  Z_t (\o) \n  : = \n \sum^{N_\o}_{\ell = 1} \b1_{\{ \tau_{\ell-1} (\o) < t \le \tau_\ell (\o) \}} \cZ^\ell_t (\o)
 ~ \hb{and} ~    K_t (\o) \n : =  \n   \sum^{N_\o}_{\ell = 1}
 \(   \cK^\ell \( \tau_\ell (\o) \land t, \o \) \n  - \n
    \cK^\ell \( \tau_{\ell-1} (\o) \land t,  \o \) \) , ~  t  \n \in \n  [0, T] .
  \eea
 The former implies that $ \dis \int_0^T  \n  |Z_t (\o)|^2 dt \n = \dn \int_0^{\tau(\o)}  \n  |Z_t (\o)|^2 dt
 \n  = \n \sum^{N_\o}_{\ell = 1}  \n
 \int_{\tau_{\ell-1} (\o)}^{\tau_\ell (\o)}  \n
 |\cZ^\ell_t (\o)|^2 dt \n \le \n  \sum^{N_\o}_{\ell = 1} \n  \int_0^T  \n  |\cZ^\ell_t (\o)|^2 dt  \n < \n  \infty $,
 so  $Z  \n \in \n  \hH^{2,0}$. We see from
   the latter of \eqref{eq:a311}   that the path $\{K_t (\o)\}_{t \in [0,T]}$ is equal to $0$ over
 period $[0,\nu(\o)]$, is
 a connection of continuous increasing pieces
    from $ \cK^\ell \( \tau_{\ell-1} (\o)  ,  \o \)$ to $\cK^\ell \( \tau_\ell (\o)  , \o \)$,
    $\ell  \n = \n 1, \cds  \n  , N_\o$ over period $[\nu(\o),\tau(\o)]$,
    and then remains constant over period $\big[\tau(\o),T\big]$. Thus,
    $\{K_t (\o)\}_{t \in [0,T]}$ is a continuous increasing path, which shows
    $ K   \n \in \n   \hK^0 $.

    Let $ \ell \in \hN $.  One can deduce that
     \bea
 K_t & \tn \dn = & \tn \dn  \sum_{i \in \hN}
 \big( \cK^i_{\tau_i \land   t} \n  - \n   \cK^i_{\tau_{i-1} \land   t} \big)
   \n = \n  \sum^\ell_{i =1}
 \big( \cK^i_{\tau_i   \land t} \n  - \n   \cK^i_{\tau_{i-1}   \land t} \big)
  \n = \n   \sum^\ell_{i =1} \n
 \(  \n -      Y_{\tau_i   \land t} \n  + \n   Y_{\tau_{i-1}   \land t}
  \n -  \dn    \int_{\tau_{i-1} \land t}^{\tau_i   \land t}  \n  g(s, Y_s, \cZ^i_s) ds
  \n + \dn  \int_{\tau_{i-1} \land t}^{\tau_i   \land t} \n   \cZ^i_s d B_s  \n  \) \nonumber \\
   & \tn \dn =  & \tn \dn  \sum^\ell_{i =1}  \n \(  \n -     Y_{\tau_i   \land t}  \n + \n  Y_{\tau_{i-1}   \land t}
   \n - \dn  \int_{\tau_{i-1}   \land t}^{\tau_i   \land t}  \n   g(s, Y_s,  Z_s) ds
   \n + \dn  \int_{\tau_{i-1}   \land t}^{\tau_i   \land t}  \n   Z_s d B_s  \n \)
  \n = \n  -    Y_{\tau_\ell   \land t} \n + \n Y_{\nu \land t}
   \n - \n  \int_{\nu \land t}^{\tau_\ell  \land t} \n  g(s, Y_s,  Z_s) ds
  \n + \n  \int_{\nu \land t}^{\tau_\ell   \land t} \n  Z_s d B_s  \nonumber \\
    & \tn \dn  =   & \tn \dn   -    Y_t \n + \n Y_\nu
   \n - \dn  \int_\nu^t  g(s, Y_s,  Z_s) ds
   \n + \dn  \int_\nu^t Z_s d B_s , ~ \fa  t \n  \in \n  [\nu, \tau_\ell].   \label{eq:b361}
 \eea
 It follows that
 \bea      \label{eq:a201}
  Y_t  \n = \n   Y_{\tau_\ell  }   \n + \n  \int_t^{\tau_\ell }  g(s, Y_s,  Z_s) ds
   \n  + \n  K_{\tau_\ell}  \n - \n K_t
  \n - \n  \int_t^{\tau_\ell }  Z_s d B_s , \q \fa  t \in [\nu, \tau_\ell] .
 \eea
 Since the increment of $K$ over $[\tau_{i-1}, \tau_i]$ is that of $\cK^i$
 over $[\tau_{i-1}, \tau_i]$ for any $i \in \hN$,
 \eqref{eq:a199}   implies   that
 \bea      \label{eq:a203}
 \int_\nu^{\tau_\ell} (Y_t -L_t) d K_t = \sum^\ell_{i =1} \int_{\tau_{i-1}}^{\tau_i} (Y_t -L_t)  d K_t
  = \sum^\ell_{i =1} \int_{\tau_{i-1}}^{\tau_i} (Y_t -L_t) d \cK^i_t  = 0 , \q \pas
 \eea

  Because of $\hP \{ Y_\tau \n \ge \n L_\tau \} \n = \n 1$, \eqref{RBSDEL_tau} clearly holds \pas ~
  on the set $ \{\nu  \n = \n  \tau \} $,
  and $ \big\{ ( Y_{\nu \vee t} )(\o)   \n \equiv \n  (Y_\nu) (\o) \big\}_{  t   \in    [0,T] } $
  is a constant path for any  $ \o  \n \in \n  \{\nu  \n = \n  \tau \} $.
  Let  $\cN_4$ be the $\hP-$null set such that
  for any $\o \n \in \n \{\nu \n < \n \tau\}    \cap    \cN^c_2     \cap    \cN^c_4 $
  and $\ell  \n \in \n  \hN$,
  \beas
  \hb{ \eqref{eq:a181} and \eqref{eq:a203}  hold  on  scenario $\o$,
  and   $ \big\{(Y_{\nu \vee (\tau_\ell \land t)}) (\o)\big\}_{t \in [0,T]}$
  is a continuous path \big(see \eqref{eq:a205}\big). }
  \eeas
  For any   $\o  \n \in \n   \{\nu  \n = \n  \tau \}  \n \cap \n
   \big(  \cN_1    \cup    \cN_2    \cup    \cN_4 \big)^c$,
  we can deduce from \eqref{eq:a201} that   \eqref{RBSDEL_tau} holds on  scenario $\o$ and
  $ \big\{(Y_{\nu \vee t}) (\o)= (Y_{\nu \vee (\tau \land t)}) (\o)\big\}_{t \in [0,T]}$ is a continuous path.  \qed

  \ss \no {\bf Proof of Proposition \ref{prop_RBSDE_comp}:}
 The  flat-off condition of reflected BSDEs  implies that  \pas
  \beas
  0 \n  \le \n  \int_t^s  \n  \b1_{\{ Y^1_r > Y^2_r \}}      d K^1_r
    \n   =   \n    \int_t^s  \n  \b1_{\{ L^1_r = Y^1_r > Y^2_r   \}}    d  K^1_r
   \n  \le \n   \int_t^s  \n    \b1_{\{ L^1_r  > L^2_r  \}} d K^1_r = 0 , \q \fa 0 \le t < s \le T .
  \eeas
  It follows that \pas
  \beas
   \int_t^s  \n  \b1_{\{ Y^1_r > Y^2_r \}}     ( d K^1_r- dK^2_r)
   = -  \int_t^s  \n  \b1_{\{ Y^1_r > Y^2_r \}}       dK^2_r \le 0,   \q \fa 0 \le t < s \le T .
  \eeas
 Then we can apply Proposition \ref{prop_BSDE_comp_basic} over period $[0,T]$ with $V^i \n = \n  K^i$,
 $i \n = \n 1,2$ to get the conclusion. \qed

 \ss \no {\bf Proof of Theorem \ref{thm_RBSDE_exist}:}
  {\bf (1) (existence)}
   For any $n \in \hN$, we define function $g_n$ as in \eqref{def_fn}, which satisfies (H1)$-$(H5) since
     $L \in \hS^1_+$.
 In light of Proposition \ref{prop_BSDE_exist}, the BSDE$(\xi,g_n)$ admits a unique solution
 $ (Y^n, Z^n) \in   \underset{p \in (0,1)}{\cap} ( \hS^p \times  \hH^{2,p} ) $ such that $Y^n$ is of class $($D$)$.
   Also,  Proposition \ref{prop_BSDE_comp} shows that
   for any $\o \in \O$ except on a $\hP-$null set $\cN$
 \bea \label{eq:a043}
  Y^n_t (\o) \le  Y^{n+1}_t (\o)  , \q    \fa       t \in [0,T] ,~ \fa n \in \hN.
 \eea
  We can let \eqref{eq:a043} hold for any $\o \in \O$ by setting
  $ Y^n_t (\o) : = \b1_{\{\o \in \cN^c \}} Y^n_t (\o)   $,
  $ (t,\o) \in [0,T] \times \O  $, $n \in \hN$ \Big(each modified $Y^n$
  still belongs to $\underset{p \in (0,1)}{\cap}   \hS^p$, of class (D) and satisfies BSDE$(\xi,g_n)$ with $Z^n$\Big).

  Applying Proposition \ref{prop_monotone_result_0} with $(Y^n,Z^n,J^n) = (Y^n,Z^n,0)$, $n \in \hN$ shows that
  the limit process $ Y_t := \lmtu{n \to \infty} Y^n_t $, $t \in [0,T]$ is
  an $\bF-$predictable process of class $($D$)$    satisfying
    $ \hE \n \[ \underset{t \in [0,T]}{\sup} |Y_t|^p \]  \n < \n \infty$, $\fa p  \n \in \n  (0,1)$.
   It   follows that $ \underset{t \in [0,T]}{\sup} \((Y^1_t)^-  \n + \n  Y^+_t\)
     \le Y^{1}_* + Y_* < \infty $, \pas ~
  As    $   Y_T = \lmtu{n \to \infty} Y^n_T =  \xi \ge L_T  $, \pas,
     applying   Proposition \ref{prop_monotone_result} with $(\nu,\tau) = (0, T)$ yields   that
     $ Y \in \underset{p \in (0,1)}{\cap}   \hS^p $ solves RBSDE$(\xi,g,L)$
     with some $(Z,K) \in \hH^{2,0} \times \hK^0 $.
     Moreover,  applying Lemma \ref{lem_RBSDE_estimate} with   $ (\nu,\tau) = (0, T) $
   and using H\"older's inequality  show  that
      \beas
    \hE \[ \bigg( \int_0^T |Z_s|^2 ds \bigg)^{p/2} \]  +  \hE \[   K^p_T     \]
    \le C_p \, \hE \[ (Y_*)^p \] + C_p \( \hE   \int_0^T  h_t dt  \)^p < \infty , \q \fa p \in (0,1) .
    \eeas
  Namely,  $(Z,K)   \in   \underset{p \in (0,1)}{\cap} (  \hH^{2,p} \times \hK^p )$.

\ss \no {\bf (2) (uniqueness)}
 Let $(Y^1,Z^1,K^1), (Y^2,Z^2,K^2)
 \in   \underset{p \in (0,1)}{\cap} ( \hS^p \times  \hH^{2,p} \times \hK^p ) $
 be two solutions of RBSDE$(\xi,g,L)$ such that $Y^1$, $Y^2$ is of class  (D).
 We know from  Proposition \ref{prop_RBSDE_comp} that
 $ \hP \{ Y^1_t = Y^2_t, \; \fa t \in [0,T]  \} = 1 $, so it holds \pas ~   that
 \beas
     \int_t^T g(s,Y^1_s, Z^1_s   )  ds + K^1_T - K^1_t  -\int_t^T Z^1_s d B_s
   =  \int_t^T g(s,Y^2_s, Z^2_s   )  ds + K^2_T - K^2_t  -\int_t^T Z^2_s d B_s    , \q    t \in [0,T] .
 \eeas
 Comparing martingale parts on both side shows that  $Z^1_t=Z^2_t$, \dtp ~ Then it follows that \pas
  \beas
  K^1_t \n   =   \n   Y^1_0  \n  - \n   Y^1_t  \n  - \n   \int_0^t g(s,Y^1_s, Z^1_s   )  ds
   \n  + \n    \int_0^t Z^1_s d B_s
   \n  = \n   Y^2_0  \n  - \n   Y^2_t  \n  - \n   \int_0^t g(s,Y^2_s, Z^2_s   )  ds
   \n  + \n    \int_0^t Z^2_s d B_s  \n  = \n   K^2_t , ~  t \in [0,T].
  \eeas

\ss \no {\bf (3) (proof of \eqref{eq:b661} and \eqref{eq:b437})}
 Fix $\nu \n \in \n \cT $ and $\ga   \n \in \n  \cT_{\nu,T}$.
 We will simply denote $\tau_\sharp(\nu)$ by $\wh{\tau}$.
 The uniform integrability of $\{Y_\ga\}_{\ga \in \cT}$ implies that $Y_\ga \n \in \n  L^1(\cF_\ga)$,
 so we see from \eqref{eq:b411}  that \pas
 \bea  \label{eq:b559}
      Y^{\ga, Y_\ga}_t    =   Y_\ga
        + \int_t^\ga  g \big(s,  Y^{\ga, Y_\ga}_s , Z^{\ga, Y_\ga}_s     \big)   ds
       - \int_t^\ga  Z^{\ga, Y_\ga}_s dB_s   , \q \fa t   \in    [\nu, \ga ] .
  \eea
  Since  it holds \pas ~ that
  \beas
      Y_t   =  Y_\ga       + \int_t^\ga  g \big(s,  Y_s , Z_s     \big)   ds
     +   K_\ga  - K_t - \int_t^\ga  Z_s dB_s   , \q \fa t   \in    [\nu, \ga ] ,
  \eeas
  applying Proposition \ref{prop_BSDE_comp_basic}
          with $(Y^1, Z^1, V^1) \n = \n \big(Y^{\ga, Y_\ga}, Z^{\ga, Y_\ga} ,0\big)$
     and $(Y^2, Z^2, V^2)   \n = \n  (Y,Z,K) $ yields that \pas,
$     Y^{\ga, Y_\ga}_t  \n \le \n  Y_t                $ for any $ t  \n \in \n  [\nu,\ga] $.
    In particular,
       \bea   \label{eq:b663}
       \cE^g_{\nu, \ga}[Y_\ga]    \n = \n  Y^{\ga, Y_\ga}_\nu  \n \le \n  Y_{\nu}     ,  \q   \pas
       \eea
       As  $  Y_\ga   \n \ge   \n    \b1_{\{\ga < T\}}  L_\ga   \n + \n \b1_{\{\ga = T\}} \xi
   \n  =  \n   \cR_\ga  $, \pas,  we see from the monotonicity of $g-$evaluation that
   \bea    \label{eq:b665}
     Y_{\nu} \ge \cE^g_{\nu, \ga}[Y_\ga] \ge \cE^g_{\nu, \ga}[\cR_\ga]   , \q  \pas
   \eea

 \ss   Since it holds \pas\;that
 $   Y_t \dn > \dn  \cR_t  \n = \dn  L_t $ for any $  t  \n \in \dn   [\nu,\wh{\tau} ) $,
     the flat-off condition in RBSDE\;$(\xi,g,L)$ implies  that \pas ~
    $  K_t  \n = \n  K_\nu   $ for any $ t  \n \in  \n  [ \nu , \wh{\tau}   ]  $.
     Then it holds \pas ~ that
        \beas
     Y_t  \n =  \n  Y_{\wh{\tau} \land \ga}  \n  + \dn  \int_t^{\wh{\tau} \land \ga} \n   g \big(s,  Y_s,  Z_s  \big)   ds
    \n + \n K_{\wh{\tau} \land \ga}  \n - \n  K_t      \n   -  \dn  \int_t^{\wh{\tau} \land \ga} \n  Z_s dB_s
     \n = \n  Y_{\wh{\tau} \land \ga}    \n + \dn  \int_t^{\wh{\tau} \land \ga} \n   g \big(s,  Y_s,  Z_s     \big)   ds
                    \n - \n   \int_t^{\wh{\tau} \land \ga} \n  Z_s dB_s   ,
                    \q \fa t  \n \in \n  \big[\nu,\wh{\tau}  \n \land \n  \ga\big] .
       \eeas
    Similar to \eqref{eq:b559}, one has that \pas
     \beas
    Y^{\wh{\tau} \land \ga, Y_{\wh{\tau} \land \ga}}_t    =   Y_{\wh{\tau} \land \ga}
        + \int_t^{\wh{\tau} \land \ga}  g \big(s,  Y^{\wh{\tau} \land \ga, Y_{\wh{\tau} \land \ga}}_s ,
        Z^{\wh{\tau} \land \ga, Y_{\wh{\tau} \land \ga}}_s     \big)   ds
       - \int_t^{\wh{\tau} \land \ga}  Z^{\wh{\tau} \land \ga, Y_{\wh{\tau} \land \ga}}_s dB_s   , \q \fa t   \in    \big[\nu, \wh{\tau} \land \ga \big] .
     \eeas
 Applying Proposition \ref{prop_BSDE_comp_basic} again yields that \pas,
 $Y_t \n = \n  Y^{\wh{\tau} \land \ga, Y_{\wh{\tau} \land \ga}}_t  $
 for any $t  \n \in \n  \big[\nu, \wh{\tau}  \n \land \n  \ga \big]$.
 It thus follows that
        \bea \label{eq:b667}
                   Y_\nu = Y^{\wh{\tau} \land \ga, Y_{\wh{\tau} \land \ga}}_\nu
         = \cE^g_{\nu,  \wh{\tau} \land \ga}    \[ Y_{\wh{\tau} \land \ga} \] , \q  \pas,
        \eea
         which together with \eqref{eq:b663} proves \eqref{eq:b661}.

      As $Y_T \n = \n \xi \n = \n \cR_T$, \pas,
  we can deduce from  the continuity of process $Y$ and the right-continuity of process $\cR$    that
 $Y_{\wh{\tau}}  \n = \n  \cR_{\wh{\tau}}$, \pas ~
 So taking $\ga \n = \n T$ in \eqref{eq:b667} yields that
 $Y_\nu  \n = \n  \cE^g_{\nu,  \wh{\tau}  }    \[ Y_{\wh{\tau}  } \]
  \n = \n  \cE^g_{\nu,  \wh{\tau}  }    \[ \cR_{\wh{\tau}  } \]$, \pas,
 which together with \eqref{eq:b665} implies \eqref{eq:b437}.       \qed

  \subsection{Proof of Theorem \ref{thm_DRBSDE_exist} }

  {\bf (1)  (existence)}
  We shall follow   \cite{Hamadene_Hassani_2005}'s approach by pasting local solutions
  to construct a global solution of DRBSDE $(\xi,g,L,U)$, see our introduction for a  synopsis.

 \ss \no   {\bf (1a)} {\it  \big(increasing penalization scheme\big)}

  For $n \in \hN$, we  define function $g_n$ as in \eqref{def_fn}
  which satisfies (H1)$-$(H5) since  $L \in \hS^1_+$.
  Theorem \ref{thm_RBSDE_exist} and Remark \ref{rem_RBSDEU2} show that
 the following reflected BSDE   with   
 generator $g_n$ and  upper obstacle $U$
 \bea
 \label{RBSDEU2}
   \begin{cases}
   \dis      U_t \ge  Y_t=    \xi + \int_t^T g_n ( s,Y_s, Z_s   )  ds - J_T + J_t
   -\int_t^T Z_s d B_s    , \q    t \in [0,T]  ,         \vspace{1mm}     \\
   \dis \int_0^T  (U_t - Y_t   ) d J_t = 0 .
   \end{cases}
    \eea
 admits    a unique solution
 $ (Y^n, Z^n, J^n) \in   \underset{p \in (0,1)}{\cap} ( \hS^p \times  \hH^{2,p} \times \hK^p ) $
 such that $Y^n$ is of class (D).  In light of Proposition \ref{prop_RBSDE_comp} and Remark \ref{rem_RBSDEU2},
 it holds   for any $\o \in \O$ except on a $\hP-$null set $\cN$ that
 \bea \label{eq:b155}
  Y^n_t (\o) \le  Y^{n+1}_t (\o)   , \q    \fa       t \in [0,T] ,~ \fa n \in \hN.
 \eea
  We can let \eqref{eq:b155} hold for any $\o \in \O$ by setting
  $ Y^n_t (\o) : = \b1_{\{\o \in \cN^c \}} Y^n_t (\o)   $,
  $ (t,\o) \in [0,T] \times \O  $, $n \in \hN$ \Big(each modified $Y^n$
  still belongs to $\underset{p \in (0,1)}{\cap}   \hS^p$, of class (D)
  and satisfies \eqref{RBSDEU2} with $(Z^n,J^n)$\Big).
  By Proposition \ref{prop_monotone_result_0},
  the limit process $ Y_t := \lmtu{n \to \infty} Y^n_t $, $t \in [0,T]$
  is an $\bF-$predictable process of class $($D$)$ that  satisfies
    \bea \label{eq:b311}
     \hE \n \[ \underset{t \in [0,T]}{\sup} |Y_t|^p \]  \n < \n \infty , \q \fa p  \n \in \n  (0,1) .
     \eea

    Let $\nu \n \in \n  \cT$. For any $n  \n \in \n  \hN$, define a  stopping time
  $\ga^n_\nu  \n := \n  \inf\{ t  \n \in \n  [\nu, T] \n : Y^n_t  \n = \n  U_t  \}  \n \land \n  T  \n \in \n  \cT$.
  As it holds \pas ~ that $Y^n_t  \n < \n  U_t$ for any $t \n \in \n \big[\nu, \ga^n_\nu\big)$,
  we can deduce from the flat-off condition in \eqref{RBSDEU2}   that
   $ \hP \big\{ J^n_t  \n = \n  J^n_\nu   , ~ \fa t  \n \in  \n  [ \nu ,\ga^n_\nu   ]   \big\}  \n = \n  1 $.
    It then follows that \pas
    \bea \label{eq:b117}
  0 & \tn = & \tn  J^n_{\ga^n_\nu} - J^n_t
  =  Y^n_{\ga^n_\nu} - Y^n_t
  + \int_t^{\ga^n_\nu} g_n   ( s,Y^n_s, Z^n_s   )  ds
  - \int_t^{\ga^n_\nu} Z^n_s d B_s .
  \eea
  Clearly, $\ga^n_\nu$     is decreasing in $n $,
  and their limit $ \ga_\nu : = \lmtd{n \to \infty} \ga^n_\nu \ge \nu $ is still a stopping time thanks to
  the right continuity of filtration $\bF$. We claim that
  \bea   \label{eq:b151}
     Y_{\ga_\nu} = \b1_{\{\ga_\nu = T\}}  \xi  +   \b1_{\{\ga_\nu < T \}} U_{\ga_\nu } , \q \pas
   \eea
 (which will be shown in the appendix).
 So  $ Y_{\ga_\nu} 
  \n \ge \n  \b1_{\{\ga_\nu =  T\}} L_T  \n + \n  \b1_{\{\ga_\nu < T \}}  L_{\ga_\nu  }  \n  = \n  L_{\ga_\nu}  $,   \pas ~
   Since $ \hE   \bigg[ \underset{t \in [0,T]}{\sup} |Y^1_t|^p \n  + \n \underset{t \in [0,T]}{\sup} |Y_t|^p \bigg]
   \n < \n \infty$, $\fa p  \n \in \n  (0,1)$
  and since it holds \pas ~ that
  \bea \label{eq:b157}
  Y^n_t =  Y^n_{\ga_\nu}   + \int_t^{\ga_\nu} g_n   ( s,Y^n_s, Z^n_s   )  ds
  - \int_t^{\ga_\nu} Z^n_s d B_s , \q   \fa  t \in  [ \nu ,\ga_\nu   ]
  \eea
  for any $n \in \hN$ by \eqref{eq:b117}, applying Proposition \ref{prop_monotone_result}
  to $\big\{ (Y^n, Z^n  )\big\}_{n \in \hN}$ yields that
  process $\big\{ Y_{\nu \vee (\ga_\nu \land t)} \big\}_{t \in [0,T]}$ has \pas ~ continuous paths
  and there exist $(Z^\nu,K^\nu)  \n \in \n  \hH^{2,0}  \n \times \n  \hK^0 $
  such that \pas
  \bea  \label{eq:b333}
  \begin{cases}
  \dis  L_t \le  Y_t
  = Y_{\ga_\nu}  + \int_t^{\ga_\nu}  g(s, Y_s, Z^\nu_s) ds
  + K^\nu_{\ga_\nu}  - K^\nu_t
  - \int_t^{\ga_\nu}  Z^\nu_s dB_s, \q \fa t \in [\nu, \ga_\nu  ] ,   \vspace{1mm}  \\
  \dis \int_\nu^{\ga_\nu}  (Y_t-L_t) d K^\nu_t = 0 .
  \end{cases}
  \eea
  Since $\hE[|Y_\nu|] \n < \n \infty$ by the uniform integrability of $\{Y_\z \}_{\z \in \cT}$,
  Lemma \ref{lem_RBSDE_estimate}, H\"older's inequality and \eqref{eq:b311} show that
  \bea   \label{eq:b315}
    \hE \[  \(\int_\nu^{\tau  }     |Z_t|^2 dt \)^{p/2} \]
   \le  C_p \, \hE \bigg[ \underset{t \in [\nu, \tau]}{\sup} | Y_t |^p \bigg]
   + C_p    \( \hE \n \int_\nu^\tau h_t dt \)^p    < \infty , \q \fa   p \in (0, 1) .
   \eea

 \no  {\bf (1b)} {\it   \big(decreasing penalization scheme\big)}

 Similar to $g_L$ discussed in Remark \ref{rem_assum} (4),
 $ g_U (t,\o,y) := \(y-U_t(\o)\)^+   $, $ ( t,\o, y )  \in [0,T] \times \O \times \hR $
 is clearly a  $\sP  \n  \otimes \n   \sB(\hR) / \sB(\hR) -$measurable function satisfying (H2)$-$(H4).
 For any $n \in \hN$,   we see from Remark \ref{rem_assum}   (3) that
   \beas
     \wt{g}_n (t,\o, y,z) \n :=  \n    g(t,\o, y,z)  \n - \n  n \(y  \n - \n  U_t(\o)\)^+ ,
     \q \fa (t,\o, y,z)  \n \in \n  [0,T]  \n \times \n  \O  \n \times \n  \hR  \n \times \n  \hR^d
    \eeas
  defines a  generator,
 and  Theorem \ref{thm_RBSDE_exist} shows that      RBSDE$\big(\xi, \wt{g}_n, L\big)$
 admits    a unique solution
 $ \big( \wt{Y}^n, \wt{Z}^n, \wt{K}^n \big) \in   \underset{p \in (0,1)}{\cap} ( \hS^p \times  \hH^{2,p} \times \hK^p ) $
 such that $\wt{Y}^n$ is of class (D). Since $\wt{g}_n$ is decreasing in $n$,
     Proposition \ref{prop_RBSDE_comp} shows that \pas
 \bea  \label{eq:b155b}
  \wt{Y}^n_t   \ge  \wt{Y}^{n+1}_t    , \q    \fa       t \in [0,T] ,~ \fa n \in \hN.
 \eea
 As in \eqref{eq:b155}, we can assume that \eqref{eq:b155b} holds  everywhere on $  \O$.

 \ss Set $\big(\wt{L}, \wt{U}\big): = ( -   U , - L) \in \hS^1_+ \times \hS^1_- $.
 For any $n \in \hN$, $ \big(\wh{Y}^n, \wh{Z}^n, \wh{J}^n \big) :  = \big( -\wt{Y}^n, -\wt{Z}^n, -\wt{K}^n \big) $
 satisfies that \pas
 \bea
   \wt{U}_t & \tn = & \tn  - L_t \ge \wh{Y}^n_t  =    - \xi - \int_t^T g \big(s, \wt{Y}^n_s, \wt{Z}^n_s \big) ds
 + n \int_t^T \big(    \wt{Y}^n_s  \n - \n  U_s \big)^+ ds
 - \wt{K}^n_T + \wt{K}^n_t  + \int_t^T \wt{Z}^n_s d B_s \nonumber \\
  & \tn = & \tn \wh{Y}^n_T + \int_t^T g_- \big(s, \wh{Y}^n_s, \wh{Z}^n_s \big) ds
  + n \int_t^T \big(   \wh{Y}^n_s  \n - \n  \wt{L}_s \big)^- ds
 + \wh{J}^n_T - \wh{J}^n_t  - \int_t^T \wh{Z}^n_s d B_s   , \q t \in [0,T] .   \label{eq:b159}
 \eea
 Since $g_-$ is a  generator by Remark \ref{rem_assum}   (1),
  applying  Proposition \ref{prop_monotone_result_0} to
  $\Big\{\big(\wh{Y}^n, \wh{Z}^n, \wh{J}^n \big)\Big\}_{n \in \hN}$
  yields that   $ \wh{Y}_t  \n := \n  \lmtu{n \to \infty} \wh{Y}^n_t $, $t  \n \in \n  [0,T]$
  is an $\bF-$predictable process of class $($D$)$ that  satisfies
  $   \hE \n \[ \underset{t \in [0,T]}{\sup} \big| \wh{Y}_t \big|^p \]  \n < \n \infty $, $  \fa p  \n \in \n  (0,1) $.

 \ss  Let $\nu \n \in \n  \cT$.  The   stopping times
  $\tau^n_\nu  \n := \n  \inf\big\{ t  \n \in \n  [\nu, T] \n : \wh{Y}^n_t  \n = \n  \wt{U}_t  \big\}  \n \land \n  T
   \n = \n  \inf\big\{ t  \n \in \n  [\nu, T] \n : \wt{Y}^n_t  \n = \n  L_t  \big\}  \n \land \n  T  \n \in \n  \cT$
  is decreasing in $n $.    Analogous to \eqref{eq:b151},
     $ \tau_\nu : = \lmtd{n \to \infty} \tau^n_\nu \ge \nu $ is still a stopping time that satisfies
  \bea   \label{eq:b161}
     \wh{Y}_{\tau_\nu} = - \b1_{\{\tau_\nu = T\}}  \xi +   \b1_{\{\tau_\nu < T \}} \wt{U}_{\tau_\nu }
     \ge - \b1_{\{\tau_\nu = T\}}  U_T -   \b1_{\{\tau_\nu < T \}} L_{\tau_\nu }
     \ge - U_{\tau_\nu }   \ge \wt{L}_{\tau_\nu } , \q \pas
   \eea
   For any $n \in \hN$, similar to \eqref{eq:b157}, we can deduce from \eqref{eq:b159}  that \pas
   \beas
  \wh{Y}^n_t =  \wh{Y}^n_{\tau_\nu}   + \int_t^{\tau_\nu} g_n   \big(s,\wh{Y}^n_s, \wh{Z}^n_s   \big)  ds
  + n \int_t^{\tau_\nu} \big(   \wh{Y}^n_s  \n - \n  \wt{L}_s \big)^- ds
  - \int_t^{\tau_\nu} \wh{Z}^n_s d B_s , \q   \fa  t \in  [ \nu ,\tau_\nu   ] .
  \eeas
  As $ \hE \n \[ \underset{t \in [0,T]}{\sup} \big|\wh{Y}^1_t\big|^p
  + \underset{t \in [0,T]}{\sup} \big|\wh{Y}_t\big|^p \]
   \n < \n \infty$, $\fa p  \n \in \n  (0,1)$, using \eqref{eq:b161} and
    applying Proposition \ref{prop_monotone_result} yield that
  process $\big\{ \wh{Y}_{\nu \vee (\tau_\nu \land t)} \big\}_{t \in [0,T]}$ has \pas ~ continuous paths
  and there exist $\big(\wh{Z}^\nu,\wh{K}^\nu\big)  \n \in \n  \hH^{2,0}  \n \times \n  \hK^0 $
  such that \pas
  \bea  \label{eq:b319}
  \begin{cases}
  \dis  \wt{L}_t \le  \wh{Y}_t
  = \wh{Y}_{\tau_\nu}  + \int_t^{\tau_\nu}  g_-  (s, \wh{Y}_s, \wh{Z}^\nu_s) ds
  + \wh{K}^\nu_{\tau_\nu}  - \wh{K}^\nu_t
  - \int_t^{\tau_\nu}  \wh{Z}^\nu_s dB_s, \q \fa t \in [\nu, \tau_\nu  ] ,   \vspace{1mm}  \\
  \dis \int_\nu^{\tau_\nu}  (\wh{Y}_t-\wt{L}_t) d \wh{K}^\nu_t = 0 .
  \end{cases}
  \eea
  Since $g_-$ satisfies (H4) and (H5) with the same function $h$ as $g$, an analogy to \eqref{eq:b315} shows that
    \bea   \label{eq:b317}
    \hE \[  \(\int_\nu^{\tau  }     |\wh{Z}_t|^2 dt \)^{p/2} \]
   \le  C_p \, \hE \bigg[ \underset{t \in [\nu, \tau]}{\sup} \big| \wh{Y}_t \big|^p \bigg]
   + C_p    \( \hE \n \int_\nu^\tau h_t dt \)^p    < \infty , \q \fa   p \in (0, 1) .
   \eea

 Set $\big(\wt{Y}, \wt{Z}^\nu,\wt{J}^\nu\big) = \big( - \wh{Y}, \n - \n \wh{Z}^\nu, -\wh{K}^\nu\big)$,
 it   follows from \eqref{eq:b319} that \pas
 \bea  \label{eq:b331}
  \begin{cases}
  \dis  U_t \ge  \wt{Y}_t
  = \wt{Y}_{\tau_\nu}  + \int_t^{\tau_\nu}  g(s, \wt{Y}_s, \wt{Z}^\nu_s) ds
  - \wt{J}^\nu_{\tau_\nu}  + \wt{J}^\nu_t
  - \int_t^{\tau_\nu}  \wt{Z}^\nu_s dB_s, \q \fa t \in [\nu, \tau_\nu  ] ,   \vspace{1mm}  \\
  \dis \int_\nu^{\tau_\nu}  ( U_t - \wt{Y}_t   ) d \wt{J}^\nu_t = 0 .
  \end{cases}
  \eea

 \no  {\bf (1c)} {\it  Next, we show that except on a $\hP-$null set $\cN_1$
 \bea  \label{eq:b339}
 L_t \le   \wt{Y}_t = Y_t   \le U_t , \q t \in [0,T] .
 \eea }

 Given $n \in \hN$, we set $V^n_t \n := \n  n \int_0^t (Y^n_s  \n - \n  L_s)^- ds  \n - \n  J^n_t$
 and $ \wt{V}^n_t  \n := \n  - n \int_0^t (\wt{Y}^n_s  \n - \n  U_s)^+ ds  \n + \n  K^n_t $, $t  \n \in \n  [0,T]$.
 As $  (Y^n, Z^n, J^n) $ solves \eqref{RBSDEU2} and
 $\big( \wt{Y}^n, \wt{Z}^n, \wt{K}^n \big)$ solves RBSDE$\big(\xi, \wt{g}_n, L\big)$, it holds \pas ~ that
  \bea \label{eq:b335}
  Y^n_t \n \le \n  U_t \q \hb{and} \q \wt{Y}^n_t  \n \ge \n  L_t , \q \fa t  \n \in \n  [0,T]  .
  \eea
 We can then deduce that   \pas
 \beas
  && \hspace{-1.2cm} \int_t^s \b1_{\{Y^n_r > \wt{Y}^n_r\}}  \big( d V^n_r \n  - \n  d \wt{V}^n_r \big)
   \n \le  \n   n  \n  \int_t^s \b1_{\{Y^n_r > \wt{Y}^n_r\}}  \n
   \( (Y^n_r  \n - \n  L_r)^-  \n + \n  (\wt{Y}^n_r  \n - \n  U_r)^+ \)  dr
   \n = \n  n \n  \int_t^s \n  \Big( \b1_{\{L_r \ge Y^n_r > \wt{Y}^n_r\}} (Y^n_r  \n - \n  L_r)^-  \\
    &    &    +   \b1_{\{Y^n_r > \wt{Y}^n_r \ge U_r\}} (\wt{Y}^n_r  \n - \n  U_r)^+ \Big) dr
     \n \le \n  n \n  \int_t^s  \n  \Big( \b1_{\{L_r   > \wt{Y}^n_r\}} (Y^n_r  \n - \n  L_r)^-
     \n + \n    \b1_{\{Y^n_r >   U_r\}} (\wt{Y}^n_r  \n - \n  U_r)^+ \Big) dr  \n = \n  0  ,
    \q \fa 0  \n \le \n  t  \n < \n  s  \n \le \n  T     .
 \eeas
 Since $Y^n_T  \n = \n  \wt{Y}^n_T  \n = \n  \xi$, \pas,
 applying Proposition \ref{prop_BSDE_comp_basic} over period $[0,T]$
 with $g^1 \n = \n g^2 \n = \n g$, $(Y^1,Z^1,V^1)   \n = \n  (Y^n, Z^n, V^n  ) $ and
  $(Y^2,Z^2,V^2)   \n = \n \big(\wt{Y}^n, \wt{Z}^n, \wt{V}^n\big) $ yields that
  $\hP \big\{ Y^n_t \n \le \n \wt{Y}^n_t , ~ \fa t  \n \in \n  [0,T]  \big\}  \n = \n  1 $.
  It follows that \pas
  \bea  \label{eq:b337}
  Y_t = \lmtu{n \to \infty} Y^n_t \le \lmtd{n \to \infty} \wt{Y}^n_t = \wt{Y}_t , \q t \in [0,T] .
  \eea

 On the other hand, let   $\nu \in \cT$.  By \eqref{eq:b161},
 \beas
 \wt{Y}_{\tau_\nu \land \ga_\nu}
 & \dn  \dn = &  \dn  \dn   \b1_{\{\tau_\nu > \ga_\nu \}} \wt{Y}_{\ga_\nu}
 \n + \n  \b1_{\{\tau_\nu \le \ga_\nu     \}} \wt{Y}_{\tau_\nu}
  \n = \n  \b1_{\{\tau_\nu > \ga_\nu \}} \wt{Y}_{\ga_\nu}
   \n + \n  \b1_{\{\tau_\nu \le \ga_\nu, \tau_\nu < T \}} L_{\tau_\nu}
   \n + \n  \b1_{\{\tau_\nu = \ga_\nu  = T   \}} \xi \\
  &  \dn  \dn   \le &  \dn  \dn   \b1_{\{\tau_\nu > \ga_\nu \}} U_{\ga_\nu}
   \n + \n  \b1_{\{\tau_\nu \le \ga_\nu, \tau_\nu < T \}} Y_{\tau_\nu}
   \n + \n  \b1_{\{\tau_\nu = \ga_\nu  = T   \}} \xi
   \n = \n  \b1_{\{\tau_\nu > \ga_\nu \}} Y_{\ga_\nu}  \n + \n  \b1_{\{\tau_\nu \le \ga_\nu     \}} Y_{\tau_\nu}
   \n = \n  Y_{\tau_\nu \land \ga_\nu}  , \q \pas
 \eeas
 Also, we see from \eqref{eq:b331} and \eqref{eq:b333} that \pas
 \beas
  \wt{Y}_t
 & \tn =  & \tn  \wt{Y}_{\tau_\nu \land \ga_\nu}  + \int_t^{\tau_\nu \land \ga_\nu}  g(s, \wt{Y}_s, \wt{Z}^\nu_s) ds
  - \wt{J}^\nu_{\tau_\nu \land \ga_\nu}  + \wt{J}^\nu_t
  - \int_t^{\tau_\nu \land \ga_\nu}  \wt{Z}^\nu_s dB_s,    \\  \hb{and} \q
 Y_t    & \tn = & \tn  Y_{\tau_\nu \land \ga_\nu}  + \int_t^{\tau_\nu \land \ga_\nu}  g(s, Y_s, Z^\nu_s) ds
  + K^\nu_{\tau_\nu \land \ga_\nu}  - K^\nu_t
  - \int_t^{\tau_\nu \land \ga_\nu}  Z^\nu_s dB_s, \q \fa t \in [\nu, \tau_\nu \land \ga_\nu  ] .
 \eeas
 Since both $Y$ and $\wt{Y}$ are of class (D), using \eqref{eq:b315}, \eqref{eq:b317} and
  applying Proposition \ref{prop_BSDE_comp_basic} over stochastic interval
  $[\n [ \nu, \tau_\nu   \land   \ga_\nu ]\n ]$ with
  $(Y^1,Z^1,V^1)   \n = \n \big(\wt{Y}, \wt{Z}^\nu, - \wt{J}^\nu \big) $ and
  $(Y^2,Z^2,V^2)   \n = \n (Y,Z^\nu,K^\nu) $ yield that \pas,
  $ \wt{Y}_t  \n \le \n  Y_t $, $\fa t  \n \in \n  [\nu, \tau_\nu \land \ga_\nu]$. In particular,
 one has  $ \wt{Y}_\nu \le Y_\nu $, \pas ~
 As $\nu   $ varies over $\cT$, the cross-section theorem (see Theorem IV.86 of \cite{Proba_Potential_1})
 and \eqref{eq:b335} imply that  \pas ~
 \beas
 L_t \le \lmtd{n \to \infty} \wt{Y}^n_t = \wt{Y}_t \le Y_t = \lmtu{n \to \infty} Y^n_t \le U_t , \q t \in [0,T] ,
 \eeas
 which together with \eqref{eq:b337} proves \eqref{eq:b339}.
 In particular, we see from \eqref{eq:b333} and \eqref{eq:b331} that
 $(Y, Z^\nu,K^\nu,0)$ locally solves the doubly reflected BSDE over the stochastic interval $[\n[ \nu, \ga_\nu ]\n]$
 and $(Y, \wt{Z}^\nu,0, \wt{J}^\nu ) = (\wt{Y}, \wt{Z}^\nu,0, \wt{J}^\nu ) $
 locally solves the doubly reflected BSDE over the stochastic interval $[\n[ \nu, \tau_\nu ]\n]$.

 \ss  \no  {\bf (1d) } {\it \big(construction of a   solution via pasting\big)}

 For any $n \in \hN$ and $t \in [0,T]$, set $\cI^n_t := [(t-2^{-n}) \vee 0, (t+2^{-n}) \land T]$.
 Similar to \eqref{eq:b343}, we can deduce from  the continuity of $Y^n$'s, $\wt{Y}^n$'s and \eqref{eq:b339}  that \pas
  \beas
   \hspace{-3mm}   \lmtu{n \to \infty} \n \underset{s \in \cI^n_t }{\inf} \, Y_s
   & \tn \dn  = & \tn  \dn    \lmtu{n \to \infty} \n  \underset{s \in \cI^n_t}{\inf}   \lmtu{m \to \infty}   Y^m_s
     \n   \ge  \n   \lmtu{m \to \infty}  \n  \lmtu{n \to \infty}  \n   \underset{s \in \cI^n_t}{\inf}   Y^m_s
       \n = \n    \lmtu{m \to \infty} Y^m_t    \n = \n  Y_t \n  = \n  \wt{Y}_t  \n = \n   \lmtd{m \to \infty} \wt{Y}^m_t
    \n  = \n  \lmtd{m \to \infty} \lmtd{n \to \infty}  \n
      \underset{s \in \cI^n_t}{\sup}   \wt{Y}^m_s   \\
    & \tn  \dn    \ge  & \tn  \dn   \lmtd{n \to \infty}  \underset{s \in \cI^n_t}{\sup} \lmtd{m \to \infty}   \wt{Y}^m_s
      =  \lmtd{n \to \infty}  \underset{s \in \cI^n_t}{\sup} \wt{Y}_s
      =  \lmtd{n \to \infty}  \underset{s \in \cI^n_t}{\sup} Y_s
      \ge \lmtu{n \to \infty} \underset{s \in \cI^n_t }{\inf} \, Y_s , \q    \fa       t \in  [0, T ] ,
  \eeas
 which shows  that $Y$ is a continuous process. So   $Y \in \underset{p \in (0,1)}{\cap}   \hS^p  $ by \eqref{eq:b311}.

 Let $\nu_1 := 0$, we recursively set stopping times $\nu'_\ell := \ga_{\nu_\ell}   $,
  $\nu_{\ell+1} := \tau_{\nu'_\ell}   $, $\ell \in \hN$, and define processes
 \bea  \label{eq:b353}
     Z_t \n := \n  \sum_{\ell \in \hN} \b1_{\{\nu_\ell < t \le \nu'_\ell \}} Z^{\nu_\ell}_t
  \n + \n  \b1_{\{\nu'_\ell < t \le \nu_{\ell+1} \}} \wt{Z}^{\nu'_\ell}_t , ~
  K_t  \n := \n   \sum_{\ell \in \hN} \(  K^{\nu_\ell}_{\nu'_\ell \land t}
   \n - \n  K^{\nu_\ell}_{\nu_\ell \land t} \) ,   ~
  J_t  \n := \n   \sum_{\ell \in \hN} \(  \wt{J}^{\,\nu'_\ell}_{\nu_{\ell+1} \land t}
   \n - \n  \wt{J}^{\,\nu'_\ell}_{\nu'_\ell \land t} \),
  ~ t  \n \in \n  [0,T] . \q
 \eea
 Since $\big\{\b1_{\{\nu_\ell < t \le \nu'_\ell \}} \big\}_{t \in [0,T]}$ and
 $\big\{\b1_{\{\nu'_\ell < t \le \nu_{\ell+1} \}} \big\}_{t \in [0,T]}$ are
 $\bF-$adapted c\`agl\`ad  processes (thus    $\bF-$predictable) for each $\ell \in \hN$,
 the process $Z$ is   $\bF-$predictable.
 Also, it is clear that    $K$ and $J$ are $\bF-$adapted processes with  $ K_0 \n = \n J_0 = \n  0 $.

 Let $\cN_2$ be the $\hP-$null set such that for any $\o \n \in \n  \cN^c_2$, the paths $L_\cd (\o)$, $U_\cd (\o)$
 $Y_\cd (\o)$  are continuous and $L_t(\o)  \n < \n  U_t(\o)$ for any $  t  \n \in \n  [0,T]$.
 By \eqref{eq:b151} and \eqref{eq:b161}, it holds except on a $\hP-$null set $\cN_3$ that
 \bea  \label{eq:b341}
 \b1_{\{ \nu'_\ell  < T \}} Y_{ \nu'_\ell } =   \b1_{\{ \nu'_\ell  < T \}} U_{ \nu'_\ell  }
 \q \hb{and} \q \b1_{\{ \nu_{\ell+1}  < T \}} \wt{Y}_{\nu_{\ell+1}}
 =  \b1_{\{\nu_{\ell+1} < T \}} L_{\nu_{\ell+1} } , \q \fa \ell  \n \in \n  \hN .
 \eea
 We claim that $\{\nu_n \}_{n \in \hN}$ is stationary: more precisely,
 for any $\o \n \in \n ( \cN_1 \cup \cN_2 \cup \cN_3 )^c $
  \bea  \label{eq:b351}
 T = \nu_{N_\o} (\o) \hb{ for some } N_\o \in \hN .
 \eea
 Assume not, then it holds for some $ \o  \n \in \n  ( \cN_1 \cup \cN_2 \cup \cN_3 )^c $
 that $\nu_n (\o ) < T $ for each $n \in \hN$.
 Given $n \in \hN$, as $ \nu_n (\o ) \le \nu'_n (\o ) \le  \nu_{n+1} (\o ) < T $, \eqref{eq:b341} shows that
 \bea    \label{eq:b347}
 \(Y_{\nu'_n}\) (\o ) \n = \n  \(U_{\nu'_n}\)  (\o )
 \q \hb{and} \q    \(Y_{\nu_{n+1}}\)  (\o)  \n = \n  \big(\wt{Y}_{\nu_{n+1}}\big)  (\o)
  \n = \n  \(L_{\nu_{n+1}}\)  (\o)  .
 \eea
 Let $t_* = t_* (\o) = \lmtu{n \to \infty} \nu_n (\o) = \lmtu{n \to \infty} \nu'_n (\o) \in [0, T] $.
 As $n \to \infty$ in \eqref{eq:b347}, we see from the continuity of paths $L_\cd (\o)$, $U_\cd (\o)$ and
 $Y_\cd (\o)$ that
 \beas
 L_{t_*} (\o)  \n = \n  \lmt{n \to \infty} \(L_{\nu_{n+1}}\)  (\o)
  \n = \n  \lmt{n \to \infty} \big(\wt{Y}_{\nu_{n+1}}\big)  (\o)
  \n = \n  \wt{Y}_{t_*}  (\o)  \n = \n  Y_{t_*}  (\o)  \n = \n  \lmt{n \to \infty} \(Y_{\nu'_n}\) (\o )
  \n = \n  \lmt{n \to \infty} \(U_{\nu'_n}\)  (\o )  \n = \n  U_{t_*}  (\o ) .
 \eeas
 A contradiction appears, so \eqref{eq:b351} holds. Then the three sums in \eqref{eq:b353} are finite sums.
 An analogous   discussion to the one below \eqref{eq:a311}
 shows that $Z \n \in \n  \hH^{2,0}$ and $K,J  \n \in \n  \hK^0$.

 Let $\ell \in \hN$ with $\ell \ge 2$. Similar to \eqref{eq:b361},
 we can deduce from \eqref{eq:b331}, \eqref{eq:b333} and \eqref{eq:b339} that \pas
 \beas
 && \hspace{-9mm}  K_t  \n - \n   J_t
  = \sum^{\ell-1}_{i = 1}  \(  K^{\nu_i}_{\nu'_i \land t} - K^{\nu_i}_{\nu_i \land t} \)
  - \sum^{\ell-1}_{i = 1} \(  \wt{J}^{\,\nu'_i}_{\nu_{i+1} \land t} - \wt{J}^{\,\nu'_i}_{\nu'_i \land t} \) \\
 && \hspace{-4mm} = \sum^{\ell-1}_{i = 1} \( - \n Y_{\nu'_i \land t}  \n + \n  Y_{\nu_i \land t}
  \n - \dn  \int_{\nu_i \land t}^{\nu'_i \land t}  \n  g \big(s, Y_s, Z^{\nu_i}_s\big) ds
  \n + \dn  \int_{\nu_i \land t}^{\nu'_i \land t}  \n  Z^{\nu_i}_s d B_s
   \n - \n  \wt{Y}_{\nu_{i+1} \land t}  \n + \n  \wt{Y}_{\nu'_i \land t}
  \n - \dn  \int_{\nu'_i \land t}^{\nu_{i+1}  \land t} \n  g \big(s, \wt{Y}_s, \wt{Z}^{\nu'_i}_s\big) ds
  \n + \dn  \int_{\nu'_i \land t}^{\nu_{i+1}  \land t} \n  \wt{Z}^{\nu'_i}_s d B_s \) \\
&& \hspace{-4mm}  = \sum^{\ell-1}_{i = 1} \( - \n  Y_{\nu_{i+1} \land t}  \n + \n  Y_{\nu_i \land t}
  \n - \n  \int_{\nu_i \land t}^{\nu_{i+1} \land t} \n  g \big(s, Y_s, Z_s\big) ds
  \n + \n  \int_{\nu_i \land t}^{\nu_{i+1} \land t} \n  Z_s d B_s \)
  \n = \n  - \n  Y_t  \n + \n  Y_0   \n - \n  \int_0^t \n  g \big(s, Y_s, Z_s\big) ds
  \n + \n  \int_0^t \n  Z_s d B_s , ~ \fa t  \n \in \n  [0,\nu_\ell] .
 \eeas
 It follows that \pas
 \bea \label{eq:b371}
 Y_t = Y_{\nu_\ell}  \n + \n  \int_t^{\nu_\ell} \n  g \big(s, Y_s, Z_s\big) ds
 \n + \n K_{\nu_\ell} \n - \n K_t \n - \n   J_{\nu_\ell}  \n + \n  J_t
  \n - \n  \int_t^{\nu_\ell} \n  Z_s d B_s , \q \fa t  \n \in \n  [0,\nu_\ell] .
 \eea
 Since the increment of $K$ over $[\nu_i, \nu'_i]$ is that of $K^{\nu_i}$ over $[\nu_i, \nu'_i]$
 ($K$ is constant over $[\nu'_i, \nu_{i+1}]$)
 and since the increment of $J$ over $[\nu'_i, \nu_{i+1}]$ is that of $J^{\nu'_i}$ over $[\nu'_i, \nu_{i+1}]$
 ($J$ is constant over $[\nu_i, \nu'_i]$), \eqref{eq:b331}, \eqref{eq:b333} and \eqref{eq:b339} again imply that
 \bea
 \int_0^{\nu_\ell} \n  (Y_t \n - \n  L_t) dK_t
 & \dn \dn  =& \dn \dn   \sum^{\ell-1}_{i = 1} \int_{\nu_i}^{\nu'_i} \n  (Y_t  \n - \n  L_t) d K_t
  \n = \n  \sum^{\ell-1}_{i = 1} \int_{\nu_i}^{\nu'_i} \n  (Y_t  \n - \n  L_t) d K^{\nu_i}_t  \n = \n  0 , \label{eq:b373} \\
\hb{and} ~ \;
 \int_0^{\nu_\ell} \n  (U_t \n - \n Y_t ) d J_t
 & \dn \dn =& \dn \dn   \int_0^{\nu_\ell}  \n \big(U_t \n - \n \wt{Y}_t \big) d J_t
   \n = \n  \sum^{\ell-1}_{i = 1} \int_{\nu'_i}^{\nu_{i+1}} \n  \big(U_t \n - \n \wt{Y}_t\big) d J_t
  \n = \n  \sum^{\ell-1}_{i = 1} \int_{\nu'_i}^{\nu_{i+1}} \n  \big(U_t \n - \n \wt{Y}_t\big) d J^{\nu'_i}_t  \n = \n   0 ,
  \q \pas \qq   \label{eq:b375}
 \eea
   Clearly, $ Y_T = \lmtu{n \to \infty} Y^n_T = \xi 
   $, \pas ~ Letting $\ell \to \infty$ in \eqref{eq:b371}, \eqref{eq:b373} and \eqref{eq:b375}, we see from
   \eqref{eq:b351} and \eqref{eq:b339} that $(Y,Z, K, J)$ solves DRBSDE$(\xi,g,L,U)$.

    \ss \no  {\bf (2) (proof of \eqref{eq:b671}$-$\eqref{eq:b537})}
 Fix $\nu \n \in \n \cT $.
 We will simply    denote $\tau^*_\nu$ by $\wh{\tau}$ and $\ga^*_\nu$ by $\wh{\ga}$.
 Since it holds \pas ~ that
 \beas
    Y_t \n > \n  L_t, \q \fa    t  \n \in \n  \big[\nu,\wh{\tau}\big)
    \q \hb{and} \q Y_t \n < \n  U_t, \q \fa    t  \n \in \n  \big[\nu,\wh{\ga}\big) ,
    \eeas
     the flat-off conditions in DRBSDE$(\xi,g,L,U)$ implies  that \pas
    \bea   \label{eq:b611}
     K_t  \n = \n  K_\nu   , \q \fa t  \n \in  \n  [ \nu , \wh{\tau}   ]
      \q \hb{and} \q   J_t  \n = \n  J_\nu, \q \fa    t  \n \in \n  \big[\nu,\wh{\ga}\big]  .
     \eea

 Let $\tau, \ga   \n \in \n  \cT_{\nu,T}$, we see from \eqref{eq:b611}  that \pas
     \beas
Y_t = Y_{  \wh{\tau} \land \ga} + \int_t^{  \wh{\tau} \land \ga}
 g (s, Y_s, Z_s) ds  - J_{  \wh{\tau} \land \ga} + J_t
 - \int_t^{  \wh{\tau} \land \ga} Z_s d B_s , \q \fa t \in \big[\nu, \wh{\tau} \land \ga\big] .
\eeas
 As $Y_{  \wh{\tau} \land \ga} \n \in \n  L^1(\cF_{\wh{\tau} \land \ga})$ by
   the uniform integrability of $ \{ Y_{\ga'}  \}_{\ga' \in \cT}$,
  \eqref{eq:b411} shows that \pas
   \bea  \label{eq:b615}
    Y^{\wh{\tau} \land \ga,Y_{\wh{\tau} \land \ga}}_t    = Y_{\wh{\tau} \land \ga}
        + \int_t^{\wh{\tau} \land \ga}  g \Big(s,  Y^{\wh{\tau} \land \ga,Y_{\wh{\tau} \land \ga}}_s ,
        Z^{\wh{\tau} \land \ga,Y_{\wh{\tau} \land \ga}}_s     \Big)   ds
       - \int_t^{\wh{\tau} \land \ga}  Z^{\wh{\tau} \land \ga,Y_{\wh{\tau} \land \ga}}_s dB_s   ,
       \q \fa t   \in    \big[\nu, \wh{\tau} \land \ga \big] .
     \eea
 Applying Proposition \ref{prop_BSDE_comp_basic}
          with $(Y^1, Z^1, V^1) \n = \n (Y,Z,-J) $
     and $(Y^2, Z^2, V^2)
      \n = \n \Big(Y^{\wh{\tau} \land \ga,Y_{\wh{\tau} \land \ga}}, Z^{\wh{\tau} \land \ga,Y_{\wh{\tau} \land \ga}} ,0\Big)
        $ yields that \pas, $Y_t \n \le  \n  Y^{\wh{\tau} \land \ga,Y_{\wh{\tau} \land \ga}}_t$
        for any $t   \n  \in  \n    \big[\nu, \wh{\tau}  \n \land \n  \ga \big]$. It follows that
  \bea  \label{eq:b673}
  Y_\nu  \n \le \n   Y^{\wh{\tau} \land \ga,Y_{\wh{\tau} \land \ga}}_\nu  \n = \n  \cE^g_{\nu, \wh{\tau} \land \ga}
 \big[Y_{\wh{\tau} \land \ga}\big] , \q  \pas
 \eea
  Similarly, we can deduce that
  \bea \label{eq:c673}
  Y_\nu  \n \ge \n \cE^g_{\nu, \tau \land \wh{\ga}  }
 \big[Y_{\tau \land \wh{\ga}  }\big]  , \q  \pas,
 \eea
  proving \eqref{eq:b671}.

 The continuity of processes $Y$, $L$ and $U$ implies that
 $ \b1_{ \{\wh{\tau} < T\}}  Y_{\wh{\tau}} \n = \n \b1_{ \{\wh{\tau} < T\}}  L_{\wh{\tau}}  $
 and $\b1_{ \{\wh{\ga} < T\}}  Y_{\wh{\ga}} \n = \n \b1_{ \{\wh{\ga} < T\}}  U_{\wh{\ga}}$, \pas ~
 It follows that \pas
 \bea
    R(\wh{\tau}, \ga) & \tn  \tn = & \tn  \tn    \b1_{\{\wh{\tau} < \ga\} }  L_{\wh{\tau}}
    \n +  \n  \b1_{\{ \ga \le \wh{\tau} \}\cap \{\ga < T\}}   U_\ga  \n +  \n    \b1_{\{\wh{\tau} = \ga =  T \}} \xi
     \ge \b1_{\{\wh{\tau} < \ga\} }  Y_{\wh{\tau}}
    \n +  \n  \b1_{\{ \ga \le \wh{\tau} \}\cap \{\ga < T\}}   Y_\ga  \n +   \n   \b1_{\{\wh{\tau} = \ga =  T \}} Y_T
      \n = \n  Y_{  \wh{\tau} \land \ga} , \qq \qq \label{eq:b675} \\
 \hb{and }   R(\tau, \wh{\ga})  & \tn  \tn  =  & \tn  \tn   \b1_{\{\tau < \wh{\ga}\} }  L_{\tau}
    \n +  \n  \b1_{\{ \wh{\ga} \le \tau \}\cap \{\wh{\ga} < T\}}   U_{\wh{\ga}}
     \n +  \n    \b1_{\{\tau = \wh{\ga} =  T \}} \xi
     \le \b1_{\{\tau < \wh{\ga}\} }  Y_{\tau}
    \n +  \n  \b1_{\{ \wh{\ga} \le \tau \}\cap \{\wh{\ga} < T\}}   Y_{\wh{\ga}}
     \n +  \n    \b1_{\{\tau = \wh{\ga} =  T \}} Y_T
     \n  = \n  Y_{  \tau \land \wh{\ga}} . \label{eq:c675}
 \eea
Then   \eqref{eq:b673}, \eqref{eq:c673} and the monotonicity of $g-$evaluation show that
 \beas
 \cE^g_{\nu, \tau \land \wh{\ga}} \big[R(\tau, \wh{\ga})\big]  \n \le \n
  \cE^g_{\nu, \tau \land \wh{\ga}} \big[Y_{\tau \land \wh{\ga}}\big]   \n \le \n
 Y_\nu \n \le  \n  \cE^g_{\nu, \wh{\tau} \land \ga} \big[Y_{\wh{\tau} \land \ga}\big]
  \n \le \n  \cE^g_{\nu, \wh{\tau} \land \ga} \big[R(\wh{\tau}, \ga)\big] , \q  \pas
  \eeas
  Taking essential supremum over $\tau \in \cT_{\nu,T}$
  and essential infimum over $\ga \in \cT_{\nu,T}$ respectively yields that
 \bea
 \esup{\tau \in \cT_{\nu,T}} \einf{\ga \in \cT_{\nu,T}} \cE^g_{\nu, \tau \land \ga} \big[R(\tau, \ga)\big] & \tn \le  & \tn
\einf{\ga \in \cT_{\nu,T}} \esup{\tau \in \cT_{\nu,T}} \cE^g_{\nu, \tau \land \ga} \big[R(\tau, \ga)\big] \n \le \n
\esup{\tau \in \cT_{\nu,T}} \cE^g_{\nu,  \tau  \land \wh{\ga}} \big[R( \tau , \wh{\ga})\big] \nonumber \\
  & \tn \le & \tn  Y_\nu
   \n \le \n  \einf{\ga \in \cT_{\nu,T}} \cE^g_{\nu, \wh{\tau} \land \ga} \big[R(\wh{\tau}, \ga)\big]
   \n \le \n  \esup{\tau \in \cT_{\nu,T}} \einf{\ga \in \cT_{\nu,T}}
   \cE^g_{\nu, \tau \land \ga} \big[R(\tau, \ga)\big] , \q \pas \qq \qq   \label{eq:b623}
 \eea
 By \eqref{eq:b611} again, it holds \pas ~ that
 \beas
Y_t = Y_{  \wh{\tau} \land \wh{\ga}} + \int_t^{  \wh{\tau} \land \wh{\ga}}  g (s, Y_s, Z_s) ds
 - \int_t^{  \wh{\tau} \land \wh{\ga}} Z_s d B_s , \q \fa t \in \big[\nu, \wh{\tau} \land \wh{\ga}\big] .
\eeas
Comparing it to \eqref{eq:b615} with $\ga  \n = \n  \wh{\ga}$, we can deduce from
applying Proposition \ref{prop_BSDE_comp_basic}   that
\pas, $Y_t \n =  \n  Y^{\wh{\tau} \land \wh{\ga},Y_{ \wh{\tau} \land  \wh{\ga}}}_t$
        for any $t   \n  \in  \dn    \big[\nu, \wh{\tau}  \n \land \n  \wh{\ga} \big]$.
        Taking $ \ga  \n = \n  \wh{\ga} $ in \eqref{eq:b675} and
        $ \tau  \n = \n  \wh{\tau} $ in \eqref{eq:c675} yields that
  $Y_\nu  \n = \n   Y^{\wh{\tau} \land \wh{\ga},Y_{ \wh{\tau} \land  \wh{\ga}}}_\nu
   \n = \n  \cE^g_{\nu, \wh{\tau} \land \wh{\ga}} \big[Y_{ \wh{\tau} \land  \wh{\ga}}\big]
 = \cE^g_{\nu, \wh{\tau} \land \wh{\ga}} \big[ R( \wh{\tau} ,  \wh{\ga}) \big] $, \pas,
 which together with \eqref{eq:b623}
  proves \eqref{eq:b811} and \eqref{eq:b537}.

   \ss \no  {\bf (3) (uniqueness)}  Let $(\sY,\sZ, \sK, \sJ)
     \n \in \n   \( \underset{p \in (0,1)}{\cap}  \hS^p \n \) \n
      \n \times \n   \hH^{2,0}  \n \times \n  \hK^0  \n \times \n  \hK^0   $ be another solution of
      DRBSDE\,$(\xi,g,L,U)$    such that $\sY$  is of class $($D$)$. Since $\sY$ also satisfies
      \eqref{eq:b537}, it holds for any $t \in [0,T]$ that
      \bea  \label{eq:b631}
            \sY_t \n =  \n   \esup{  \tau \in  \cT_{t,T}}  \einf{  \ga \in  \cT_{t,T}}
        \cE^g_{t, \tau \land \ga}\big[R(\tau,\ga) \big]
         \n =  \n   \einf{  \ga \in  \cT_{t,T}}  \esup{  \tau \in  \cT_{t,T}}
        \cE^g_{t, \tau \land \ga}\big[R(\tau,\ga) \big] \n =  \n Y_t , \q \pas
      \eea
      The continuity of $Y$ and $\sY$ then shows that \pas
      \beas
  & &     \xi  + \int_t^T g \( s,Y_s, Z_s \)  ds    +    K_T   -   K_t
   -     J_T   +   J_t
   - \int_t^T Z_s d B_s = Y_t = \sY_t \\
   & = &  \xi  + \int_t^T g \( s,\sY_s, \sZ_s \)  ds    +    \sK_T   -   \sK_t
   -   \n  \sJ_T   +  \n   \sJ_t
   - \int_t^T \sZ_s d B_s   , \q    t \in [0,T] .
      \eeas
   Comparing the martingale parts on both sides shows that $Z_t = \sZ_t$, \dtp, and it follows that \pas
 \bea   \label{eq:b633}
 K_t - J_t = \sK_t - \sJ_t , \q    t \in [0,T] .
 \eea
 The flat-off conditions in DRBSDE$(\xi,g,L,U)$ implies  that \pas
 \bea \label{eq:b635}
  K_t \n = \dn  \int_0^t \n  \b1_{\{Y_s = L_s\}} d K_s , ~ \sK_t  \n = \dn  \int_0^t \n  \b1_{\{\sY_s = L_s\}} d \sK_s, ~
  J_t  \n = \dn  \int_0^t \n  \b1_{\{Y_s = U_s\}} d J_s , ~
  \sJ_t  \n = \dn  \int_0^t \n  \b1_{\{\sY_s = U_s\}} d \sJ_s , \q t  \n \in \n  [0,T] .
 \eea
 As $\hP\{L_t  \n < \n  U_t, \,  \fa t    \in \n  [0,T] \}  \n = \n  1$, we can deduce that \pas
 \beas
 \int_0^t \n  \b1_{\{Y_s = U_s\}} d K_s
 \n = \dn  \int_0^t \n  \b1_{\{Y_s = U_s\}} \b1_{\{Y_s = L_s\}} d K_s  \n = \n  0
 ~ \; \hb{and}  ~ \;   \int_0^t \n  \b1_{\{\sY_s = U_s\}} d \sK_s
 \n = \dn  \int_0^t \n  \b1_{\{\sY_s = U_s\}} \b1_{\{\sY_s = L_s\}} d \sK_s  \n = \n  0 , \q t  \n \in \n  [0,T] ,
 \eeas
 which together with \eqref{eq:b635}, \eqref{eq:b631} and \eqref{eq:b633} leads to that \pas
 \beas
 J_t  & \tn =  & \tn  \int_0^t \b1_{\{Y_s = U_s\}} d J_s + \int_0^t \b1_{\{\sY_s = U_s\}} d \sK_s
 =  \int_0^t \b1_{\{Y_s = U_s\}} d J_s + \int_0^t \b1_{\{Y_s = U_s\}} d \sK_s  \\
  & \tn =  & \tn  \int_0^t \b1_{\{Y_s = U_s\}} d \sJ_s + \int_0^t \b1_{\{Y_s = U_s\}} d K_s
 =  \int_0^t \b1_{\{\sY_s = U_s\}} d \sJ_s + \int_0^t \b1_{\{Y_s = U_s\}} d K_s = \sJ_t , \q t  \n \in \n  [0,T] .
 \eeas
 Then it easily follows from \eqref{eq:b633} that \pas, $K_t = \sK_t$, $\fa t  \n \in \n  [0,T]$. \qed

\appendix
\renewcommand{\thesection}{A}
\refstepcounter{section}
\makeatletter
\renewcommand{\theequation}{\thesection.\@arabic\c@equation}
\makeatother

\section{Appendix}

  \begin{lemm} Given $ \xi  \n \in \n  L^1(\cF_\tau) $, let $ \xi  \n \in \n  L^1(\cF_\tau) $
  and let  $g (t,\o, y , z) = \b1_{\{t \le \tau(\o)\}} g (t,\o, y , z) $,
  $ (t,\o,y,z)  \n \in \n  [0,T]  \n \times \n  \O  \n \times \n  \hR  \n \times \n  \hR^d $ be a generator.
  Then one has
  \bea \label{eq:b409}
     \hP \big\{ Y^{\tau, \xi}_t \n = \n Y^{\tau, \xi}_{\tau \land t}, ~  \fa t  \n \in \n  [0,T] \big\}  \n = \n 1
 \q \hb{and} \q        Z^{\tau, \xi}_t =\b1_{\{t \le \tau \}}Z^{\tau, \xi}_t  ,   \q      \dtp~
  \eea
   \big(see \eqref{eq:xax011} for the notation $(Y^{\tau, \xi},Z^{\tau, \xi})$\big).
       In particular, it holds \pas ~ that
         \bea \label{eq:b411}
         Y^{\tau, \xi}_t = \xi + \int_t^\tau g \big(s, Y^{\tau, \xi}_s, Z^{\tau, \xi}_s \big) ds
         - \int_t^\tau Z^{\tau, \xi}_s d B_s , \q \fa t \in [0,\tau].
         \eea

\end{lemm}

\ss \no {\bf Proof:}  Given $n \n \in \n \hN$,  we define a stopping time
 \bea \label{eq:b407}
 \ga_n \n  := \n \inf\left\{ t  \n \in \n  [0,T] \n :
 \int_0^t \n  \big|Z^{\tau, \xi}_s\big|^2 ds  \n > \n  n \right\}  \n \land \n  T  \n \in \n  \cT  .
 \eea
  Since
 $ Y^{\tau, \xi}_{\tau \land \ga_n} \n = \n Y^{\tau, \xi}_{  \ga_n}
  \n + \n  \int_{\tau \land \ga_n}^{  \ga_n} \n  \b1_{\{s \le \tau\}}
 g \big(s, Y^{\tau, \xi}_s , Z^{\tau, \xi}_s \big) ds
  \n - \n  \int_{\tau \land \ga_n}^{  \ga_n} \n  Z^{\tau, \xi}_s d B_s
  \n = \n  Y^{\tau, \xi}_{  \ga_n}    \n - \n  \int_{\tau \land \ga_n}^{  \ga_n} \n  Z^{\tau, \xi}_s d B_s $, \pas,
  taking conditional expectation $\hE \big[\cd |\cF_{\tau \land \ga_n}\big]$ yields that \pas
  \bea \label{eq:b417}
  \q   Y^{\tau, \xi}_{\tau \land \ga_n} \n = \n  \hE \big[  Y^{\tau, \xi}_{  \ga_n}  \big|\cF_{\tau \land \ga_n}\big]
    \n = \n  \b1_{\{\tau \le  \ga_n\}} \hE \big[  Y^{\tau, \xi}_{  \ga_n}   \big|\cF_{\tau  }\big]
    \n + \n  \b1_{\{\tau >  \ga_n\}}  \hE \big[  Y^{\tau, \xi}_{  \ga_n}   \big|\cF_{\ga_n  }\big]
    \n = \n  \b1_{\{\tau \le  \ga_n\}} \hE \big[  Y^{\tau, \xi}_{  \ga_n}   \big|\cF_{\tau  }\big]
    \n + \n  \b1_{\{\tau >  \ga_n\}}     Y^{\tau, \xi}_{  \ga_n} .
  \eea
   As  $Z^{\tau, \xi} \in \underset{p \in (0,1)}{\cap} \hH^{2,p} \subset  \hH^{2,0}$, $\{\ga_n\}_{n \in \hN}$ is stationary.
 Letting $n \to \infty$, we can deduce from the uniform integrability of $\big\{Y^{\tau, \xi}_\ga\big\}_{\ga \in \cT}$ that
 \beas
 Y^{\tau, \xi}_\tau  = \b1_{\{\tau \le  T\}} \hE \big[  Y^{\tau, \xi}_T   \big|\cF_{\tau  }\big]
   + \b1_{\{\tau >  T\}}     Y^{\tau, \xi}_T = \hE \big[  Y^{\tau, \xi}_T   \big|\cF_{\tau  }\big]
   = \hE \big[  \xi   \big|\cF_{\tau  }\big] = \xi , \q \pas
 \eeas
  Then it follows that \pas
\bea
Y^{\tau, \xi}_{\tau \land t} & \tn = & \tn  Y^{\tau, \xi}_\tau
  \n + \n  \int_{\tau \land t}^\tau \n  \b1_{\{s \le \tau\}}
 g \big(s, Y^{\tau, \xi}_s , Z^{\tau, \xi}_s \big) ds
  \n - \n  \int_{\tau \land t}^\tau \n  Z^{\tau, \xi}_s d B_s \label{eq:b415} \\
 & \tn  = & \tn  \xi \n + \n  \int_t^T \n  \b1_{\{s \le \tau\}}
 g \big(s, Y^{\tau, \xi}_{\tau \land s} , \b1_{\{s \le \tau\}} Z^{\tau, \xi}_s \big) ds
  \n - \n  \int_t^T \n  \b1_{\{s \le \tau\}}  Z^{\tau, \xi}_s d B_s   , \q t \in [0,T] . \nonumber
\eea
which shows that $\big\{ \big(Y^{\tau, \xi}_{\tau \land t}, \b1_{\{t \le \tau \}}Z^{\tau, \xi}_t \big) \big\}_{t \in [0,T]}$
also solves BSDE$(\xi,g_\tau)$. Clearly, $\big\{  Y^{\tau, \xi}_{\tau \land t}  \big\}_{t \in [0,T]}$
is an $\bF-$adapted continuous process such that
$\hE \n \[ \underset{t \in [0,T]}{\sup} \big|Y^{\tau, \xi}_{\tau \land t}\big|^p \]
\le \hE \n \[ \underset{t \in [0,T]}{\sup} \big|Y^{\tau, \xi}_t \big|^p \] < \infty $ for any $p \in (0,1)$
and that $\big\{  Y^{\tau, \xi}_\ga \big\}_{\ga \in \cT_{0,\tau}}$ is uniformly integrable.
As $\big\{\b1_{\{t \le \tau\}} \big\}_{t \in [0,T]}$
 is an $\bF-$adapted c\`agl\`ad process
 \big(and thus $\bF-$predictable\big), we see that $\big\{  \b1_{\{t \le \tau \}}Z^{\tau, \xi}_t   \big\}_{t \in [0,T]}$
 is an $\bF-$predictable process satisfying
 $ \hE \n \[ \( \int_0^T \b1_{\{t \le \tau \}} |Z^{\tau, \xi}_t|^2 \)^{p/2} \]
 \le \hE \n \[ \( \int_0^T   |Z^{\tau, \xi}_t|^2 \)^{p/2} \] < \infty  $ for any $p \in (0,1)$.
 Hence,  by  the uniqueness of solution of BSDE$(\xi,g_\tau)$, \eqref{eq:b409} holds.

 Moreover, \eqref{eq:b415} can be alternatively expressed as: \pas
\beas
  Y^{\tau, \xi}_{\tau \land t}   =      \xi   +    \int_{\tau \land t}^\tau
 g \big(s, Y^{\tau, \xi}_s , Z^{\tau, \xi}_s \big) ds
    -    \int_{\tau \land t}^\tau    Z^{\tau, \xi}_s d B_s,   \q t \in [0,T]   ,
\eeas
 which leads to \eqref{eq:b411}.  \qed

  \begin{lemm}   \label{lem_RBSDE_estimate}
 Let  $g \n : [0,T]  \n \times \n  \O  \n \times \n  \hR  \n \times \n  \hR^d \to \hR$ be a
   $\sP     \n  \otimes  \n     \sB(\hR)    \n   \otimes  \n     \sB(\hR^d)/\sB(\hR) -$measurable function
   satisfying   {\rm (H1)} and  {\rm(H4)}. Given $\nu, \tau \n  \in \n  \cT$ with $\nu  \n \le \n  \tau$,
 let $(Y,Z,K)  \n \in \n  \hS^0  \n \times \n  \wt{\hH}^{2,0}  \n \times \n  \hK^0 $
 satisfy   that \pas
 \bea \label{eq:a251}
 Y_t  = Y_\tau + \int_t^\tau g(s, Y_s, Z_s) ds + K_\tau - K_t
 - \int_t^\tau   Z_s  d B_s , \q \fa t \in [\nu,\tau] .
 \eea
 If $ \hE \big[ | Y_\nu| \big] < \infty  $, then for any $p \in (0, \infty)$,
  $  \hE \[  \(\int_\nu^{\tau  }     |Z_t|^2 dt \)^{p/2} \] + \hE \big[ ( K_\tau  \n - \n K_\nu)^p \big]
   \le  C_p \, \hE \bigg[ \underset{t \in [\nu, \tau]}{\sup} | Y_t |^p \bigg]
   + C_p \hE \Big[ \( \int_\nu^\tau h_t dt \)^p  \Big]   $.

\end{lemm}

  \ss \no {\bf Proof:} Let $ \hE \big[ | Y_\nu| \big] < \infty  $ and fix $p \n \in \n  (0, \infty)$.
  By the Burkholder-Davis-Gundy  inequality, there exists
  $\fc_p > 0$ such that for any continuous local martingale $M$
  \bea  \label{eq:a211}
  \hE \[ (M_*)^p \] \le \fc_p \hE \[ \lan M \ran^{p/2}_T \] \q \hb{and} \q
  \hE \[ (M_*)^{p/2} \] \le \fc_p \hE \[ \lan M \ran^{p/4}_T \] .
  \eea

  Set $\Psi \n := \n \underset{t \in [\nu, \tau]}{\sup} | Y_t | $
 and  suppose $ \hE \[ \Psi^p \]  \n < \n  \infty$,
  otherwise the result trivially holds.
 We let $n  \n  \in \n   \hN$  and define a stopping time
 $  \tau_n   \n   := \n \inf\{ t  \n \in \n  [\nu,\tau] \n :
 \int_\nu^t     |Z_s|^2 ds  \n > \n  n \}
  \n \land \n  \tau  \n \in \n \cT  $. It is clear that  $ \nu \n \le \n  \tau_n  \n \le \n  \tau $.
 Since (H1), (H4) and H\"older's inequality imply that
 \beas
 K_{\tau_n } \dn - \n  K_\nu & \tn \dn  =  & \tn  \dn    Y_\nu  \n  - \n   Y_{\tau_n }  \n   - \n
     \int_\nu^{\tau_n }  \n    g(t,  Y_t, Z_t) dt
     \n  + \n   \int_\nu^{\tau_n }  \n      Z_t  d B_t
     \n  \le  \n  2 \Psi \n   + \n
     \int_\nu^{\tau_n }  \n    \( h_t  \n  + \n   \k |Y_t|  \n  + \n   \k | Z_t | \)  dt
     \n  + \n  \bigg| \int_0^T  \n  \b1_{\{\nu \le s \le \tau_n\}}    Z_t  d B_t \bigg| \\
  & \tn  \dn  \le  & \tn  \dn  ( 2 \n + \n \k T ) \Psi \n + \n  \int_\nu^\tau h_t dt
   \n + \n  \k \sqrt{T} \( \int_\nu^{\tau_n }  \n  |Z_t|^2 dt \)^{1/2}
  \n  + \n  \underset{t \in [0,T]}{\sup} \bigg| \int_0^t  \n
  \b1_{\{\nu \le s \le \tau_n\}} Z_s  d B_s \bigg| , \q \pas ,
 \eeas
 taking the expectation of $p-$th power, we can deduce from \eqref{eqn-d011}  and \eqref{eq:a211} that
 \if{0}
 \beas
 \hE \[ K^p_{\tau_n } \] & \tn \dn \le & \tn  \dn   ( 2 \n + \n \k T )^p \, \hE \[  \Psi^p \]
  \n + \n  \hE \[ \bigg( \int_\nu^\tau h_t dt \bigg)^p \, \]
 \dn + \n  \k^p T^{p/2} \, \hE \[ \bigg( \int_\nu^{\tau_n } |Z_t|^2 dt \bigg)^{p/2} \]
 \dn + \n \hE \[ \underset{t \in [0,T]}{\sup}  \bigg| \int_\nu^{\tau_n \land t }  \n  Z_s  d B_s \bigg|^p \]  \\
 & \tn  \dn  \le & \tn  \dn   ( 2 \n + \n \k T )^p \,   \hE \[  \Psi^p \] + \hE \[ \( \int_\nu^\tau h_t dt \)^p \, \]
 \n + \n  ( \k^p T^{p/2}  \n + \n  \fc_p ) \, \hE \[ \( \int_\nu^{\tau_n } |Z_t|^2 dt \)^{p/2} \]  .
 \eeas
 \fi
 \bea   \label{eq:a215}
 \hE \Big[ \big(K_{\tau_n } \dn - \n  K_\nu \big)^p \Big]
  \n \le    \n (1  \n \vee \n  4^{p-1}) \Bigg\{ ( 2 \n + \n \k T )^p \, \hE \[  \Psi^p \]
  \n + \n  \hE \[ \( \int_\nu^\tau h_t dt \)^p \, \]
 \dn + \n  \( \k^p T^{p/2}  \n + \n  \fc_p \)  \, \hE \[ \bigg( \int_\nu^{\tau_n } |Z_t|^2 dt \bigg)^{p/2} \]
 \Bigg\} .   \q
 \eea

  As $ \hE \big[ | Y_\nu| \big]  \n < \n  \infty  $, Corollary \ref{cor_martingale} implies that
 there exists a unique $\wt{Z}  \n \in \n  \underset{p \in (0,1)}{\cap}   \hH^{2,p} $
 such that $ \hP \big\{ \hE[Y_\nu|\cF_t] \n = \n  \hE[Y_\nu]
  \n + \n  \int_0^t \n  \wt{Z}_s dB_s, ~ \fa t  \n \in \n  [0,T] \big\}  \n = \n 1 $.
 Similar to \eqref{eq:b241},    \eqref{eq:a251} shows  that \pas
 \bea
 \wt{Y}_t & \dn \dn: =& \dn \dn \hE[Y_\nu|\cF_{\nu \land t}]  \n  + \n  Y_{\nu \vee (\tau \land t)}  \n - \n  Y_\nu
 \n = \n \hE[Y_\nu]
 \n - \n  \int_0^t \n  \b1_{\{\nu < s \le \tau\}}     g(s, Y_s, Z_s) ds
  \n - \n  \int_0^t \n  \b1_{\{\nu < s \le \tau\}}    d K_s \nonumber \\
  & \dn \dn    & \dn \dn
    +    \int_0^t  \n  \( \b1_{\{s \le \nu \}} \wt{Z}_s
  \n + \n \b1_{\{\nu < s \le \tau\}}  Z_s \) d B_s   , \q t \in [0,T] . \label{eq:a261}
 \eea
 So $\wt{Y}$ is an $\bF-$adapted continuous process, i.e. $\wt{Y} \in \hS^0$.

   Set $a : = 2 ( \k + \k^2 )  $ and
    $ \d : = \[ 3 (1 \vee 4^{p/2-1}) (1 \vee 4^{p-1}) \( \k^p T^{p/2}  \n + \n  \fc_p \)\]^{-2/p}   $.
    Applying It\^o's formula to process   $ \big\{ e^{at}  | \wt{Y}_t   |^2 \big\}_{t \in [0,T]}  $,
    we can deduce from \eqref{eq:a261}   that \pas
   \beas
  e^{at} \big|\wt{Y}_t\big|^2 & \tn = & \tn  \big(\hE[Y_\nu]\big)^2 \n + \n a \int_0^t e^{as} \big| \wt{Y}_s \big|^2 ds
    \n - \n 2 \int_0^t \n
    \b1_{\{\nu < s \le \tau\}}  e^{as}  \wt{Y}_s  g(s, Y_s, Z_s) ds
    \n +\n \int_0^t e^{as} \n \( \b1_{\{s \le \nu \}} |\wt{Z}_s|^2
  \n + \n  \b1_{\{\nu < s \le \tau\}}  |Z_s|^2 \) ds  \\
  & \tn  & \tn  \n -   2 \int_0^t \n  \b1_{\{\nu < s \le \tau\}}  e^{as}  \wt{Y}_s  d K_s
  \n + \n 2 \int_0^t  \n e^{as}  \wt{Y}_s \( \b1_{\{s \le \nu \}} \wt{Z}_s
  \n + \n  \b1_{\{\nu < s \le \tau\}}  Z_s \) d B_s    , \q t \in [0,T] .
   \eeas
   An analogy to \eqref{eq:b221} shows that $\wt{Y}_t = Y_t$, $\fa t \in [\nu, \tau]$. Hence,  it holds \pas ~ that
   \bea
     e^{a \tau} |Y_\tau|^2  &   \dn \dn = &   \dn \dn    e^{a \tau}   \big|\wt{Y}_\tau\big|^2
     =  e^{at}  \big|\wt{Y}_t\big|^2
   \n + \n a \int_t^\tau e^{as} \big| \wt{Y}_s \big|^2 ds
    \n - \n 2 \int_t^\tau \n
    \b1_{\{\nu < s \le \tau\}}  e^{as}  \wt{Y}_s  g(s, Y_s, Z_s) ds
   \n - \n   2 \int_t^\tau \n  \b1_{\{\nu < s \le \tau\}}  e^{as}  \wt{Y}_s  d K_s   \nonumber \\
  &   \dn \dn   & \dn  \dn
    +   2 \int_t^\tau  \n e^{as}  \wt{Y}_s \( \b1_{\{s \le \nu \}} \wt{Z}_s
  \n + \n  \b1_{\{\nu < s \le \tau\}}  Z_s \) d B_s
  \n +\n \int_t^\tau e^{as} \n \( \b1_{\{s \le \nu \}} |\wt{Z}_s|^2
  \n + \n  \b1_{\{\nu < s \le \tau\}}  |Z_s|^2 \) ds     \nonumber \\
   &   \dn \dn   =  &   \dn \dn   e^{at}  | Y_t |^2
    \n + \dn   \int_t^\tau \n e^{as} \( a | Y_s |^2 \n + \n |Z_s|^2 \n - \n 2  Y_s  g(s, Y_s, Z_s) \)   ds
   \n -   2 \n  \int_t^\tau \n     e^{as}  Y_s  d K_s
  \n + \n 2 \n  \int_t^\tau  \n  e^{as}  Y_s   Z_s   d B_s  , ~ \fa t \n \in \n [\nu, \tau] . \qq \q  \label{eq:a237}
   \eea
   Then    (H1) and  (H4) imply  that \pas
  \beas
 && \hspace{-1.5cm}  e^{a \nu} \big| Y_\nu \big|^2
   \n +  \n    \int_\nu^{\tau_n }   \n  e^{as}  \( a | Y_s |^2 \n + \n |Z_s|^2 \) ds
   \n  =   \n    e^{a \tau_n }  \big| Y_{\tau_n } \big|^2
   \n  +  \n    2   \n   \int_\nu^{\tau_n }  \n   e^{as}    Y_s g(s,  Y_s, Z_s) ds
   \n  +  \n    2   \n   \int_\nu^{\tau_n }  \n   e^{as}   Y_s d K_s
  \n        -  \n  2   \n   \int_\nu^{\tau_n }  \n   e^{as}   Y_s Z_s d B_s  \\
 &&     \le  \n   e^{aT} \Psi^2
    \n + \n   2  \int_\nu^{\tau_n } \n e^{as} \n \( | Y_s | h_s  \n  +  \n  \k | Y_s|^2  \n  +  \n  \k | Y_s| |Z_s| \) ds
    \n + \n   2  e^{aT}    \Psi ( K_{\tau_n } \dn - \n K_\nu)
    \n + \n   2  \bigg|  \int_0^T  \n  \b1_{\{\nu \le s \le \tau_n\}} e^{as}  Y_s Z_s d B_s \bigg|  \\
 &&    \le  \n  \Big( 1  \n + \n \frac{2}{\d}  \Big) e^{2 aT} \Psi^2
   \n + \n   2 e^{aT} \Psi  \n \cd \n  \int_\nu^\tau h_s ds \n+\n   2 ( \k   \n + \n    \k^2 )
   \int_\nu^{\tau_n }   \n  e^{as}  |Y_s|^2  ds
    \n + \n   \frac12  \int_\nu^{\tau_n }   \n  e^{as}  |Z_s|^2  ds
    \n  +  \n    \frac{\d}{2}  ( K_{\tau_n } \dn - \n K_\nu)^2 \\
    &&  \q   \n   +       2  \bigg|  \int_0^T   \n \b1_{\{\nu \le s \le \tau_n\}}  e^{as}   Y_s Z_s d B_s \bigg| .
  \eeas
  It follows that \pas
  \beas
  \q \int_\nu^{\tau_n }  \n    |Z_t|^2 dt  \n \le  \n  \int_\nu^{\tau_n }  \n e^{as}  |Z_t|^2 dt
    \n  \le   \n  \Big( 4    \n + \n        \frac{4}{\d}  \Big)   e^{2 aT} \Psi^2
     \n  +  \n  2  \( \int_\nu^\tau h_t dt \)^2   \n +  \n    \d  ( K_{\tau_n } \dn - \n K_\nu)^2
   \n + \n   4 \underset{t \in [0,T]}{\sup}
   \bigg| \int_0^t   \n \b1_{\{\nu \le s \le \tau_n\}} e^{as}  Y_s Z_s d B_s \bigg| .
  \eeas
 Taking the expectation of $p/2-$th power, we can deduce from \eqref{eqn-d011} and  \eqref{eq:a215}  that
  \beas
 && \hspace{-1.2cm} \hE \[  \(\int_\nu^{\tau_n }  \n   |Z_t|^2 dt\)^{p/2} \]
  \n \le   \n  (1 \vee 4^{p/2-1}) \Bigg\{ \n \Big( 4  \n + \n   \frac{4}{\d}   \Big)^{p/2}
  e^{apT} \hE \[  \Psi^p \]  \n  + \n  2^{p/2} \hE \[ \(  \int_\nu^\tau h_t dt  \)^p \, \]  \\
  && \q    +      \d^{p/2}  \hE \[ ( K_{\tau_n } \dn - \n K_\nu)^p \]
  \n + \n  4^{p/2} \fc_p \hE\[ \( \int_\nu^{\tau_n } \n e^{2at} |Y_t|^2  | Z_t |^2  d t \)^{p/4}\] \Bigg\} \\
  && \n     \le   \n  C_p  \hE \[  \Psi^p \] + C_p   \hE \[ \( \int_\nu^\tau h_t dt \)^p \, \]
    \n  +  \n  \frac13  \hE \[ \( \int_\nu^{\tau_n }  |Z_t|^2  dt \)^{p/2} \]
    \n  + \n  C_p \hE\[ (\Psi)^{p/2} \( \int_\nu^{\tau_n   }  \n   | Z_t |^2  d t \)^{p/4}\] \\
     && \n     \le   \n  C_p  \hE \[  \Psi^p \] + C_p \hE \[ \( \int_\nu^\tau h_t dt \)^p \, \]
    \n  +  \n  \frac23  \hE \[ \( \int_\nu^{\tau_n }  \n   |Z_t|^2  dt \)^{p/2} \]  .
  \eeas
  So   $   \hE \[  \(\int_\nu^{\tau_n }    |Z_t|^2 dt \)^{p/2} \]
   \le  C_p \hE \[ \Psi^p \] + C_p \hE \[ \( \int_\nu^\tau h_t dt \)^p \, \]   $, which together with
   \eqref{eq:a215} shows that
  \bea  \label{eq:a221}
   \hE \[  \(\int_\nu^{\tau_n }     |Z_t|^2 dt \)^{p/2} \]
   + \hE \Big[ \big(K_{\tau_n } \dn - \n  K_\nu \big)^p \Big]
   \le  C_p \hE \[ \Psi^p \] + C_p \hE \[ \( \int_\nu^\tau h_t dt \)^p \, \]   .
  \eea
  As $Z \in \wt{\hH}^{2,0}$, it holds for \pas ~ $\o \in \O$ that $\tau(\o) = \tau_{N_\o}(\o)$
  for some $N_\o \in \hN$.  Then letting
   $n \to \infty$ in \eqref{eq:a221}, we can apply the monotone convergence theorem to obtain the conclusion. \qed

 \begin{lemm}   \label{lem_increasing}
 Let $X$ be an $\bF-$optional process with \pas ~ right upper semi-continuous paths \big(i.e.,
 for any $\o \in \O$ except a $\hP-$null set $\cN_X$,
 $    X_t \n  \ge \n   \underset{s \searrow t}{\limsup} \, X_s $,
   $ \fa  t  \n  \in \n   [0,T) $\big).
  If  $X_\nu  \n  \le \n   X_{\wt{\nu}} $, \pas ~ for any  $\nu, \wt{\nu}  \n  \in \n   \cT$
  with $\nu  \n  \le \n   \wt{\nu}$, \pas, then
 $X$ is an increasing process.
 \end{lemm}

 \ss \no {\bf Proof:}
   Set  $\cD_k  := \big\{ t^k_i :=  \frac{i}{2^k} \land T   \big\}^{\lceil 2^k T \rceil}_{i = 0}  $, $\fa k \in \hN$
 and $\cD := \underset{k \in \hN}{\cup} \cD_k$.    Given $t \in [0,T)$,   we define
    $  \ul{X}_t \n  := \n    \lmtu{n \to \infty} \underset{s \in \Th^n_t }{\inf} X_s $,
   where $ \Th^n_t : =   \cD  \cap  \( t, (t+2^{-n}) \land T \] $.
    Clearly,
    \bea   \label{eq:a125}
     \Th^n_t = \underset{k > n}{\cup} \Th^{n,k}_t  ,
   \q \hb{where } \, \Th^{n,k}_t  := \cD_k \cap  \( t, (t+2^{-n}) \land T \] .
    \eea

    For any $m,n \in \hN$ with $m < n  $,  since $\Th^n_t$ is a countable subset of $ \( t, (t+2^{-n}) \land T \] $,
  the random variable $   \underset{s \in \Th^n_t }{\inf} X_s  $ is clearly $\cF_{(t+2^{-n}) \land T}-$measurable.
 So $ \ul{X}_t
 =  \lmtu{\substack{n \to \infty \\ n > m}} \underset{s \in \Th^n_t }{\inf} X_s \in \cF_{(t+2^{-m}) \land T} $.
 As $m \to \infty $,    the right-continuity of the filtration $\bF$ shows  that
 \bea    \label{eq:a113}
 \ul{X}_t \in \underset{m \in \hN}{\cap} \cF_{(t+2^{-m}) \land T}  = \cF_{t+} = \cF_t .
 \eea

 \ss \no {\bf (1)} { \it Additionally setting $\ul{X}_T := X_T \in \cF_T$, we first show
 the process $ \ul{X} $ is $\bF-$progressively measurable.}

 \ss  For any $t \n  \in \n   [0,T)$, $c  \n  \in \n   \hR$ and $n ,k  \n  \in \n   \hN$ with $k  \n  > \n   n $,
 since it holds  for   $i  \n  = \n   0, \cds , \lfloor 2^k t \rfloor $
 and   any $ s  \n  \in \n   [t^k_i, t^k_{i+1})  \n  \cap \n   [0,t] $   that
    \beas
     \Th^{n,k }_i : = \Th^{n,k }_{t^k_i} =  \{ t^k_j: j =i+1, \cds, i+ 2^{k-n} \}
    =  \Th^{n,k}_s \subset \( s, (s+2^{-n}) \land T \] \subset \( 0, (t+2^{-n}) \land T \] ,
    \eeas
    we can deduce that
     \bea
  && \hspace{-0.8cm}  \Big\{(s,\o) \n  \in \n   [0,t]  \n  \times \n   \O \n  :
   \underset{r \in \Th^{n,k}_s }{\min} X_r (\o)  \n  \ge \n    c   \Big\}
    \n  = \n   \underset{i = 0 }{\overset{\lfloor 2^k t \rfloor}{\cup}}
   \Big\{(s,\o)  \n  \in \n    \( [t^k_i, t^k_{i+1})  \n  \cap \n   [0,t] \)
   \n  \times  \n    \O \n  :  \underset{r \in \Th^{n,k}_s }{\min} X_r (\o)   \n  \ge \n    c   \Big\} \nonumber  \\
  && =  \n   \underset{i = 0 }{\overset{\lfloor 2^k t \rfloor}{\cup}}
   \Big\{(s,\o)  \n  \in \n    \( [t^k_i, t^k_{i+1})  \n  \cap \n   [0,t] \)
   \n  \times  \n    \O \n  :  \underset{r \in \Th^{n,k }_i }{\min} X_r (\o)   \n  \ge  \n   c   \Big\}
   \n   =  \n   \underset{i = 0 }{\overset{\lfloor 2^k t \rfloor}{\cup}} \,
   \underset{r \in \Th^{n,k}_i}{\cap}
   \Big\{(s,\o)  \n  \in \n    \( [t^k_i, t^k_{i+1})  \n  \cap \n   [0,t] \)
   \n  \times  \n    \O \n  :    X_r (\o)   \n  \ge \n    c   \Big\}   \nonumber  \\
 && = \n   \underset{i = 0 }{\overset{\lfloor 2^k t \rfloor}{\cup}} \,
   \underset{r \in \Th^{n,k}_i}{\cap}
  \( [t^k_i, t^k_{i+1})  \n  \cap \n   [0,t] \)
   \n  \times  \n    \{    X_r    \ge  c  \}
   \n  \in \n   \sB([0,t])  \n  \otimes \n   \cF_{(t+2^{-n}) \land T} .   \label{eq:a111}
 \eea

 Now, let $\wt{t} \in [0,T]$ and $\wt{c} \in \hR$. If $\wt{t} = 0$, then \eqref{eq:a113} shows that
  $  \big\{ (s,\o) \n  \in \n   \[0,\wt{t} \;  \]  \n  \times \n   \O \n  : \ul{X}_s(\o)  \n  > \n  \wt{c} \big\}
   \n  = \n   \{0\}  \n  \times \n   \{ \ul{X}_0   \n  > \n  \wt{c} \}
    \n  \in \n   \sB(\{0\})  \n  \otimes \n   \cF_0 $; if $\wt{t} > 0$,
    for any  $ m >  m_0 :=  \left\lceil - \frac{ \ln \wt{t} }{ \ln 2 } \right\rceil  $,
    we can deduce from  \eqref{eq:a111} and \eqref{eq:a125}  that
 \beas
 && \hspace{-1cm}  \big\{ (s,\o) \n  \in \n   \[0,\wt{t}-2^{-m}\]  \n  \times \n   \O \n  :
 \ul{X}_s(\o)  \n  > \n   \wt{c}  \big\}
  \n  = \n       \Big\{ (s,\o)  \n  \in \n   \[0,\wt{t} \n  - \n  2^{-m}\]  \n  \times \n   \O \n  :
 \lmtu{\substack{  n \to \infty \\ n > m }} \, \underset{r \in \Th^n_s }{\inf}  X_r(\o)  \n  > \n   \wt{c}  \Big\}   \\
 && =   \n   \underset{n > m}{\cup}  \Big\{(s,\o)  \n  \in \n   \[0,\wt{t} \n  - \n  2^{-m}\]
  \n  \times \n   \O \n  :
  \underset{r \in \Th^n_s }{\inf}  X_r(\o)  \n  > \n   \wt{c}  \Big\}
  \n  = \n     \underset{n > m}{\cup} \,
   \underset{\ell \in \hN}{\cup}  \Big\{(s,\o)  \n  \in \n   \[0,\wt{t} \n  - \n  2^{-m}\]
    \n  \times \n   \O \n  :
  \underset{r \in \Th^n_s }{\inf}  X_r(\o)  \n  \ge \n   \wt{c}  \n  + \n   1 / \ell  \Big\} \\
  &&    =  \n    \underset{n > m}{\cup} \,
   \underset{\ell \in \hN}{\cup}  \, \underset{k > n}{\cap} \Big\{(s,\o)
    \n  \in \n   \[0,\wt{t} \n  - \n  2^{-m}\]  \n  \times \n   \O \n  :
  \underset{r \in \Th^{n,k}_s }{\min}  X_r(\o)  \n  \ge \n   \wt{c}  \n  + \n   1 / \ell  \Big\}
  \n  \in \n   \sB\(\[0, \wt{t} \n  - \n  2^{-m}\]\)  \n  \otimes \n   \cF_{ \wt{t}  } ,
 \eeas
 which together with \eqref{eq:a113} shows that
 \beas
  \big\{ (s,\o) \n  \in \n   \[0,\wt{t} \, \]  \n  \times \n   \O \n  :
 \ul{X}_s(\o)  \n  > \n   \wt{c}  \big\}  =
  \Big\{ (s,\o) \n  \in \n  \( \underset{m > m_0}{\cup} \[0,\wt{t}-2^{-m}\] \) \n  \times \n   \O \n  :
 \ul{X}_s(\o)  \n  > \n   \wt{c}  \Big\}   \cup
 \big\{ (s,\o) \n  \in \n   \big\{ \wt{t} \, \big\}  \n  \times \n   \O \n  :
 \ul{X}_s(\o)  \n  > \n   \wt{c}  \big\} \\
 = \( \underset{m > m_0}{\cup}
  \big\{ (s,\o) \n  \in \n   \[0,\wt{t}-2^{-m}\]  \n  \times \n   \O \n  :
 \ul{X}_s(\o)  \n  > \n   \wt{c}  \big\} \) \cup
 \(   \big\{ \wt{t} \, \big\}  \times \big\{ \ul{X}_{\wt{t}}  \n  > \n   \wt{c} \big\} \)
 \n   \in  \n   \sB \(\[0, \wt{t} \, \]\)  \n  \otimes \n   \cF_{ \wt{t}  } .  \qq \q
 \eeas
 So $\L := \Big\{ \cE \subset \hR :   \big\{ (s,\o) \n  \in \n   \[0,\wt{t} \, \]  \n  \times \n   \O \n  :
 \ul{X}_s(\o)   \in \cE  \big\}
  \n   \in  \n   \sB \(\[0, \wt{t} \, \]\)  \n  \otimes \n   \cF_{ \wt{t}  }   \Big\}$
  contains all open sets of form $\( \wt{c}, \infty  \)$, which generates $\sB(\hR)$. Clearly, $\L$ is a
  $\si-$field of $\hR$. It follows that $ \sB(\hR) \subset \L $, i.e.
  $ \big\{ (s,\o) \n  \in \n   \[0,\wt{t} \, \]  \n  \times \n   \O \n  :
 \ul{X}_s(\o)   \in \cE  \big\}
  \n   \in  \n   \sB \(\[0, \wt{t} \, \]\)  \n  \otimes \n   \cF_{ \wt{t}  }  $
  for any $\cE \in  \sB(\hR) $. Hence, $\ul{X}$ is $\bF-$progressively measurable.

 \ss \no {\bf (2)} Fix $\ell \in \hN$. Since both $X$ and $\ul{X}$ are $\bF-$progressively measurable,
  the Debut theorem shows that
  \beas
  \tau_\ell : = \inf\{t \in [0,T] : \ul{X}_t \le X_t - 1/\ell \} \land  T .
  \eeas
  defines a stopping time, i.e. $\tau_\ell \in \cT$.
  We claim that $A_\ell := \{\tau_\ell  < T \} \in \cF_T $ is a $\hP-$null set:
  {\it Assume not}, so   $A_\ell \backslash \cN_X$ is not empty.
   Let $\o \in A_\ell \backslash \cN_X $ and set $s := \tau_\ell(\o) $.
  there exists  $\{s_i\}_{i \in \hN} \subset [s, T)$ with $\lmtd{i \to \infty} s_i = s  $
  such that
  \bea   \label{eq:a117}
  \ul{X}_{s_i} (\o) \le X_{s_i} (\o) - 1/\ell , \q \fa  i \in \hN  .
  \eea
   Given $m \in \hN$, we can find some
  $\wh{i} = \wh{i} (m) \in \hN$ and $\wh{n} = \wh{n}(m) \ge m $
  such that for any $ i \ge \wh{i} $ and $n \ge \wh{n}$,
   $\(s_i, (s_i+2^{-n}) \land T \] \subset \(s , (s +2^{-m}) \land T \]$ and thus
   \beas
   \Th^n_{s_i} = \( \underset{k > n}{\cup} \cD_k \) \cap \(s_i, (s_i+2^{-n}) \land T \]
   \subset \( \underset{k > m}{\cup} \cD_k \) \cap \(s , (s +2^{-m}) \land T \] = \Th^m_s .
   \eeas
  It follows that $\underset{r \in \Th^m_s}{\inf} X_r (\o) \le \underset{r \in  \Th^n_{s_i}}{\inf} X_r (\o) $.
  Letting $n \to \infty$, we see that  $\underset{r \in \Th^m_s}{\inf} X_r (\o) \le \ul{X}_{s_i} (\o) $.
  As $i \to \infty$,   \eqref{eq:a117} and the  right upper semi-continuity of $X_\cd (\o)$ imply that
    \beas
    \underset{r \in \Th^m_s}{\inf} X_r (\o)  \le \underset{i \to \infty}{\liminf} \ul{X}_{s_i} (\o)
    \le \underset{i \to \infty}{\limsup} X_{s_i} (\o) - 1/\ell \le
    \underset{r \searrow s}{\limsup} X_r (\o) - 1/\ell \le X_s (\o) - 1/\ell .
    \eeas
    Now, letting $m \to \infty$ yields that
    $ \ul{X}_s (\o) \le X_s (\o) - 1/\ell $, which shows that
    \bea \label{eq:a123}
      \ul{X}_{\tau_\ell} \le X_{\tau_\ell} - 1/\ell   \q   \hb{on } \,  A_\ell \backslash \cN_X .
    \eea

  The $\bF-$optional measurability of $X$ implies that of the stopped process
  $\big\{ X_{ \tau_\ell \land t } \big\}_{t \in [0,T]}$ (see e.g. Corollary 3.24 of \cite{HWY_1992}), so
  $ \cX^\ell_t \n  : =  \n    \b1_{\{ X_{  \tau_\ell \land t  } \le X_t \}}$, $t  \n  \in \n   [0,T]$ is also an
  $\bF-$optional  process.
  Since $ \cX^\ell_\nu  \n  = \n   \b1_{\{ X_{  \tau_\ell \land \nu  } \le X_\nu \}}  \n  = \n   1 $, \pas ~
  for any $\nu  \n  \in \n   \cT$,
  the cross-section theorem (see Theorem IV.86 of \cite{Proba_Potential_1}) shows that
  for any $\o  \n  \in \n   \O$ except on a $\hP-$null set  $\cN_\ell$,
  \bea \label{eq:a121}
 \cX^\ell_t (\o) = 1 \q \hb{or} \q \( X_{ \tau_\ell \land t  } \) (\o) \le X_t (\o),  \q \fa t \in [0,T]  .
  \eea
  Let $\o \in A_\ell \backslash (\cN_X \cup \cN_\ell) $.
  As $ X \(\tau_\ell (\o),\o \) \le X(t,\o) $, $\fa t \in \[ \tau_\ell (\o), T \]$ by \eqref{eq:a121},
  we can deduce from \eqref{eq:a123} that
  \beas
  X \(\tau_\ell (\o),\o \) \le  \ul{X} \(\tau_\ell (\o),\o \)
  \le  X \( \tau_\ell (\o), \o \) - 1/\ell  .
  \eeas
  An contradiction appears, so $ 0 = \hP(A_\ell) = \hP \{ \ul{X}_t \le  X_t - 1/\ell  \hb{ for some } t \in [0,T) \}  $.
   Letting $\ell \to \infty$  yields that
  $ \hP \{  \ul{X}_t <   X_t , \hb{ for some } t \in [0,T)  \} = \lmtu{\ell \to \infty}
  \hP \{  \ul{X}_t \le  X_t -1/\ell  \hb{ for some } t \in [0,T)  \} = 0 $, which together with
   the right upper semi-continuity of $X$ shows that except on a $\hP-$null set $\cN$
 \beas
  \ul{X}_t   \n    \ge   \n      X_t
  \n  \ge \n   \underset{s \searrow t}{\limsup} \, X_s  \n  =  \n
     \lmtd{n \to \infty} \underset{s \in   (t, (t+2^{-n})\land  T]}{\sup} X_s
    \n  \ge  \n      \lmtd{n \to \infty} \underset{s \in \Th^n_t }{\sup} X_s
   \n   \ge  \n   \lmtu{n \to \infty} \underset{s \in \Th^n_t }{\inf} X_s
   \n   = \n   \ul{X}_t,   \q \fa t  \n  \in \n   [0,T) .
 \eeas
  To wit, it holds for any $\o \in \cN^c$   that
  \bea  \label{eq:a127}
   X_t (\o) = \lmt{   \substack{  s  \searrow t  \\   s \in \cD \cap (t,T]  }  }  X_s (\o)    , \q \fa t \in [0,T).
  \eea

  Set $\wt{\cN} : = \cN \cup \( \underset{\substack{s,s' \in \cD,  s<s'}}{\cup}
  \big\{ X_s > X_{s'} \big\} \)$, which is also a $\hP-$null set. Given $\o \in \wt{\cN}^c$ and $t,t' \in [0,T]$
  with $t<t'$, let $\{s_n\}_{n \in \hN} \subset \cD \cap (t,t')$ with
  $\lmtd{n \to \infty} s_n = t$    and let $\{s'_n\}_{n \in \hN} \subset \cD \cap \( (t',T) \cup \{T\} \)$
  with $\lmtd{n \to \infty} s'_n = t'$. We can deduce from \eqref{eq:a127} that
  $   X_t (\o) = \lmt{n \to \infty} X_{s_n}(\o) \le \lmt{n \to \infty} X_{s'_n}(\o) = X_{t'} (\o) $.
  Therefore,   $X$ is an increasing process.   \qed

           \ss \no {\bf Proof of \eqref{eq:a187}:}
  {\bf (1)}   The continuity of $Y^n$'s   implies that for \pas ~ $\o \in \O$
  \bea
    \underset{s \searrow t}{\liminf} \, Y_s (\o)
   & \tn   = & \tn   \lmtu{n \to \infty} \underset{s \in (t, (t+2^{-n}) \land T]}{\inf} \, Y_s (\o)
   \n  = \n   \lmtu{n \to \infty} \underset{s \in (t, (t+2^{-n}) \land T]}{\inf}   \lmtu{m \to \infty} \, Y^m_s (\o)
     \n   \ge  \n   \lmtu{m \to \infty} \lmtu{n \to \infty}
      \underset{s \in (t, (t+2^{-n}) \land T]}{\inf} \, Y^m_s (\o) \nonumber \\
    & \tn     =  & \tn    \lmtu{m \to \infty} \underset{s \searrow t}{\liminf} \, Y^m_s (\o)
     =  \lmtu{m \to \infty} Y^m_t (\o) = Y_t (\o) , \q    \fa       t \in \big[\nu(\o), \tau(\o)\big) , \label{eq:b343}
  \eea
 which shows  that the process $\big\{Y_{\nu \vee (\tau_\ell \land t)}\big\}_{t \in [0,T]}$
 has \pas ~ right lower semi-continuous paths.
 It then  follows from \eqref{eq:a133}   that
 $\wt{K}^\ell$ has \pas ~ right upper semi-continuous paths.

 \ss  \no {\bf (2)}  {\it  We next show that
    $\wt{K}^\ell_{  \ga }$ is a weak limit of $\big\{K^n_{\tau_\ell \land \ga} \big\}_{n \in \hN}$
    in $L^2(\cF_T)$ for any $\ga \in \cT$. }

  Let $\chi \in L^2(\cF_T) $. In virtue of martingale representation theorem, there exists a unique $Z^\chi \in
  \hH^{2,2}$ such that \pas
  \beas
  M^\chi_t := \hE [\chi|\cF_t] = \hE [\chi] + \int_0^t Z^\chi_s d B_s ,   \q   \fa t \in [0,T] .
  \eeas
  Set $\z \n  = \n  \z^\ell \n  : = \n  \nu  \n \vee \n  (\tau_\ell  \n \land \n  \ga)  \n \in \n  \cT $
  and let $n  \n \in \n  \hN$.    We   define
  $ \U^{\ell,n}_t     \n :=     \n   K^n_{\nu  \vee ( \z \land t )}  \n + \n  Y^{\ell,n}_{\nu  \vee ( \z \land t )}
   \n - \n  Y^{\ell,n}_\nu   \n - \n  \( \wt{K}^\ell_{\nu  \vee ( \z \land t )}
   \n + \n  Y_{\nu  \vee ( \z \land t )}  \n - \n  Y_\nu \) $, $t  \n \in \n  [0,T]$.
  As $K^n_\nu  \n = \n  0$ by \eqref{eq:a245}, one can deduce from   \eqref{eq:a333} that \pas
  \beas
  \U^{\ell,n}_t   & \tn \n  =& \tn \n   - \int_\nu^{\nu  \vee ( \z \land t )} \n\(  g(s, Y^{\ell,n}_s, Z^n_s)
  -    g(s,  Y_s, 0) \n  - \n  \wt{h}^\ell_s \) ds
   \n + \n  \int_\nu^{\nu  \vee ( \z \land t )}  \n    \( Z^n_s  \n- \n
      \cZ^\ell_s \) d B_s  \\
    & \tn \n   =& \tn  \n - \dn  \int_0^t \n \b1_{\{ \nu < s \le  \z   \}}
     \(   g(s, Y^{\ell,n}_s, Z^n_s) \n - \n
     g(s,  Y_s, 0)   \n  - \n  \wt{h}^\ell_s  \) ds
     \n + \n  \int_0^t \n \b1_{\{ \nu < s \le  \z   \}}
      \(   Z^n_s \n - \n \cZ^\ell_s\) dB_s , ~ \;  t \in [0,T] ,
  \eeas
  thus $\U^{\ell,n}$ is   an   $\bF-$adapted continuous process.
  Since \eqref{eq:a281}, \eqref{eq:a527} and \eqref{eq:a095} shows that
    $  | \U^{\ell,n}_t | \n  \le  \n  4 \ell  \n +  \n
  K^n_{\nu  \vee ( \z \land t )}  \n + \n  \big| \wt{K}^\ell_{\nu  \vee ( \z \land t )} \big|  $,
  $ \fa  t  \n \in \n  [0,T] $,       \eqref{eqn-d011},     \eqref{eq:a093}  and \eqref{eq:a135} imply   that
 \bea  \label{eq:a153}
 \q    \hE \[ \(\U^{\ell,n}_*\)^2   \]  \n   \le   \n
 3 \hE  \[ 16 \ell^2 \n + \n \( K^n_{\tau_\ell   } \)^2
 \n  + \n    \big( \wt{K}^{\ell}_*  \big)^2  \]
     \n  \le  \n     C_0   \ell^2
      \n  + \n   C_0 \hE \int_\nu^{\tau_\ell} \( | \wt{h}^\ell_t |^2  \n  + \n   | \cZ^\ell_t |^2 \) dt < \infty ,
 \eea
 \if{0}
  This together with \eqref{eq:a157} leads to that
   \beas
   \hE \int_0^T \( |M^\chi_t|^2 \n   + \n   |\U^{\ell,n}_t|^2 \) dt
    \n \le \n T \cd \hE \[ \underset{t \in [0,T]}{\sup}  |M^\chi_t|^2
     \n  + \n  \underset{t \in [0,T]}{\sup}   |\U^{\ell,n}_t|^2 \]
     \n  \le  \n     C_0  \ell^2   \n  + \n   C_0 (\hE [\chi] )^2
    \n  + \n   C_0 \hE \int_0^T \n \( |Z^\chi_t|^2  \n  + \n   | \wt{h}^\ell_t |^2
     \n  + \n   | \cZ^\ell_t |^2 \) dt ,
   \eeas
namely, $ M^\chi, \U^{\ell,n} \in \hH^{2,2} $.
 \fi
 which shows that $ \U^{\ell,n} \in \hS^2 $.  Like  \eqref{eq:a245}, one has
  \bea    \label{eq:a357}
  \wt{K}^\ell_t = 0, \q \fa t \in [0,\nu] .
  \eea
  So $ \U^{\ell,n}_\nu \n = \n  K^n_\nu  \n - \n  \wt{K}^\ell_\nu  \n = \n  0  $.
   Integrating by parts the process $ M^\chi \U^{\ell,n} $   yields that \pas
  \bea
 \chi \U^{\ell,n}_T & \tn   = & \tn   M^\chi_T \U^{\ell,n}_T
 =  M^\chi_t \U^{\ell,n}_t + \int_t^T M^\chi_s d \U^{\ell,n}_s
  + \int_t^T \U^{\ell,n}_s d M^\chi_s \n + \n
  \lan  M^\chi, \U \ran_T \n - \n \lan  M^\chi, \U \ran_t  \nonumber \\
  & \tn  = & \tn   -  \n   \int_t^T  \n   \b1_{\{ \nu < s \le  \z \}} M^\chi_s
    \(   g(s, Y^{\ell,n}_s, Z^n_s)   \n  - \n   g(s,  Y_s, 0)   \n  - \n  \wt{h}^\ell_s \) ds
  \n   +  \n   \int_t^T  \n   \b1_{\{ \nu < s \le  \z \}}  M^\chi_s
   \(    Z^n_s  \n  - \n    \cZ^\ell_s \)   d B_s    \nonumber \\
  &  \tn   &  \tn
        + \n   \int_t^T   \n     \U^{\ell,n}_s Z^\chi_s dB_s
     \n   + \n   \int_t^T  \n   \b1_{\{  \nu < s \le  \z  \}}
    Z^\chi_s   \(    Z^n_s  \n  - \n    \cZ^\ell_s \) ds , \q t \in [0,T].   \label{eq:a155}
  \eea

  Since  Doob's  martingale inequality shows that
  $   \hE \n \[ \(  M^{\chi}_* \)^2  \] \n \le \n 4 \hE \n  \[  | M^\chi_T|^2 \]
   \n = \n  4 \hE \big[ |\chi|^2 \big]   \n < \n  \infty  $
  (i.e. $M^\chi \n \in \n  \hS^2  \n \subset \n  \hH^{2,2} $),
    applying the Burkholder-Davis-Gundy inequality and H\"older's inequality,
  we see from   \eqref{eq:a153} that
  \beas
 \q && \hspace{-1.2cm} \hE \[ \underset{t \in [0,T]}{\sup} \bigg|
  \int_0^t  \n   \b1_{\{ \nu < s \le  \z \}} M^\chi_s
   \(    Z^n_s  \n  - \n    \cZ^\ell_s \)   d B_s \bigg|
   +   \underset{t \in [0,T]}{\sup} \bigg|  \int_0^t  \n    \U^{\ell,n}_s Z^\chi_s dB_s \bigg| \]
   \n  \le  \n  C_0 \hE \[ M^{\chi}_*   \bigg(  \int_\nu^\z  \n
   \big|    Z^n_s  \n  - \n    \cZ^\ell_s \big|^2   d  s \bigg)^{1/2} \] \\
 && \hspace{-0.5cm} +    C_0 \hE \n \[   \U^{\ell,n}_*   \bigg(  \n  \int_0^T  \n  |  Z^\chi_s |^2  d s \bigg)^{1/2} \]
   \dn  \le \n  C_0 \Bigg\{ \hE \n  \[   \(  M^{\chi}_* \)^2  \]  \n \cd \n
   \hE  \n  \int_\nu^\z  \n  \big|    Z^n_s  \n  - \n    \cZ^\ell_s \big|^2   d  s \Bigg\}^{1/2}
    \dn + \n  C_0 \Bigg\{ \hE  \n  \[  \(\U^{\ell,n}_*\)^2  \]  \n \cd \n
   \hE \n  \int_0^T  \n  |  Z^\chi_s |^2   d  s \Bigg\}^{1/2}  \n < \n  \infty   \,   .
  \eeas
  Namely,    $\big\{ \n  \int_0^t  \n    \b1_{\{ \nu < s \le  \z \}} M^\chi_s
   \(    Z^n_s  \n  - \n    \cZ^\ell_s \)  d B_s \big\}_{t \in [0,T]}$ and
   $ \big\{ \n \int_0^t  \n  \b1_{\{s \ge \nu\}}   \U^{\ell,n}_s Z^\chi_s dB_s \big\}_{t \in [0,T]}$
   are uniformly integrable martingales.
   Then taking expectation in \eqref{eq:a155} for $ t = 0 $ yields that
  \beas
     \hE \[ \chi \n \(  \n  K^n_\z  \n  - \n  \wt{K}^\ell_\z \) \]
  & \tn  \dn  =  & \tn  \dn   \hE \[ \chi   \(  - Y^{\ell,n}_\z  \n  + \n   Y_\z
     \n  + \n  Y^{\ell,n}_\nu  \n  - \n   Y_\nu \) \]  \n  + \n   \hE [\chi \U^{\ell,n}_T] \nonumber  \\
   & \tn  = &  \tn       \hE \[ \chi   \(  - Y^{\ell,n}_\z  \n  + \n   Y_\z
     \n  + \n  Y^{\ell,n}_\nu  \n  - \n   Y_\nu \) \]
   - \hE    \int_0^T \b1_{\{ \nu < s \le  \z \}} M^\chi_s
   \Big(  g(s, Y^{\ell,n}_s, Z^n_s) \n  - \n   g(s,  Y^{\ell,n}_s, 0)
     \n  - \n  \wt{h}^\ell_s \Big) ds    \nonumber  \\
  & \tn \dn &  \tn  \dn   - \hE  \n  \int_0^T  \n  \b1_{\{ \nu < s \le  \z \}} M^\chi_s
   \(     g(s,  Y^{\ell,n}_s, 0)  \n  - \n   g(s,  Y_s, 0) \) ds
    \n + \n  \hE  \n  \int_0^T \n  \b1_{\{ \nu < s \le  \z \}}
   Z^\chi_s   \big(    Z^n_s  \n  - \n    \cZ^\ell_s \big)
      ds  \n : =  \n  I^n_1  \n - \n  I^n_2  \n - \n  I^n_3  \n + \n  I^n_4   .      \qq
  \eeas

  As $M^\chi   , Z^\chi \n \in \n  \hH^{2,2} \n $,
   the weak  convergence of $\big\{\b1_{\{\nu < s \le \tau_\ell\}}
    (g(s,Y^{\ell,n}_s,Z^n_s)    - \n  g(s,Y^{\ell,n}_s,0)) \big\}_{s \in [0,T]}$,
  $n \n  \in \n   \hN$ to $ \big\{\b1_{\{\nu < s \le \tau_\ell\}} \wt{h}^\ell_s \big\}_{s \in [0,T]} $ and
   that of $\big\{\b1_{\{\nu < s \le \tau_\ell\}} Z^n_s \big\}_{s \in [0,T]}$,
   $n \n  \in  \n    \hN$ to   $\big\{\b1_{\{\nu < s \le \tau_\ell\}} \cZ^\ell \big\}_{s \in [0,T]}  $
   by \eqref{eq:a353}  show that
   $ \lmt{n \to \infty} I^n_2  \n = \n  \lmt{n \to \infty} I^n_4  \n = \n  0 $.
     Since
      $ \Big| \chi  \big(  - Y^{\ell,n}_\z  \n  + \n   Y_\z
     \n  + \n  Y^{\ell,n}_\nu  \n  - \n   Y_\nu \big)  \Big|  \n \le \n  4 \ell |\chi| $
     by  \eqref{eq:a281}, \eqref{eq:a527},   \eqref{eq:a095} and since
     $  \hE[|\chi|] \n \le \n  1 \n + \n  \hE[|\chi|^2]  \n < \n  \infty $ by
       \eqref{eqn-d011b},   the dominated convergence theorem imply
 that $\lmt{n \to \infty} I^n_1
      \n = \n  0$.   Moreover, (H3)   shows that
     \beas
        \lmt{n \to \infty}  \b1_{\{ \nu < s \le  \z \}}
      \big(     g(s,  Y^{\ell,n}_s, 0) \n    -  \n   g(s,  Y_s, 0) \big)  \n = \n
        \lmt{n \to \infty}  \b1_{A_\ell \cap \{ \nu < s \le  \z \}}
      \big(     g(s,  Y^n_s, 0) \n    -  \n   g(s,  Y_s, 0) \big)  \n = \n 0  ,  \q    \dsp ,
      \eeas
     while  (H4), \eqref{eq:a281} and \eqref{eq:a095} imply that \dsp
     \beas
        \Big| \b1_{\{ \nu < s \le  \z \}}  M^\chi_s
   \big(     g(s,  Y^{\ell,n}_s, 0)  \n  - \n   g(s,  Y_s, 0) \big) \Big|
    \n \le \n    \b1_{\{ \nu < s \le  \z \}} | M^\chi_s |
     \big( 2 h_s    +   \k \big| Y^{\ell,n}_s \big| + \k |Y_s|  \big)
    \n \le \n    \b1_{\{ \nu < s \le  \z \}} | M^\chi_s |
     \big( 2 h_s    +    2 \k \ell   \big)  .
     \eeas
   As \eqref{eq:a273} and  H\"older's inequality   show that
     \beas
  \q    \hE \int_0^T \n \b1_{\{ \nu < s \le  \z \}} | M^\chi_s | \( 2 h_s  \n + \n  2 \k \ell   \)  ds
    \n \le \n   2 \ell ( 1  \n + \n  \k   T ) \, \hE \big[   M^{\chi}_*   \big]
     \n \le  \n  2 \ell ( 1  \n + \n  \k   T )
     \Big\{ \hE \[ \(  M^{\chi}_* \)^2    \] \Big\}^{1/2}
      \n < \n  \infty ,
     \eeas
   we can  apply  the dominated convergence theorem again to obtain  $\lmt{n \to \infty} I^n_3 \n = \n 0  $.
   Hence  $ \lmt{n \to \infty}  \hE \Big[ \chi  \Big( K^n_\z  \n  - \n  \wt{K}^\ell_\z \Big) \Big] = 0 $.

   Since \eqref{eq:a245} and \eqref{eq:a357} imply that for any $ n \in \hN $
       $$
    K^n_{\tau_\ell \land \ga} \n  - \n  \wt{K}^\ell_\ga
    \n = \n    K^n_{\tau_\ell \land \ga}  \n  - \n  \wt{K}^\ell_{\tau_\ell \land \ga}
    \n = \n    K^n_{ \tau_\ell \land \ga }  \n  - \n K^n_{\nu \land (\tau_\ell \land \ga)}
     \n  - \n \Big( \wt{K}^\ell_{ \tau_\ell \land \ga } \n - \n \wt{K}^\ell_{\nu \land (\tau_\ell \land \ga)} \Big)
    \n = \n    K^n_{ \nu \vee (\tau_\ell \land \ga) }  \n  - \n K^n_\nu
     \n  - \n \Big( \wt{K}^\ell_{ \nu \vee (\tau_\ell \land \ga) } \n - \n \wt{K}^\ell_\nu \Big)
     \n = \n   K^n_\z  \n  - \n  \wt{K}^\ell_\z ,
   $$
   one gets   $ \lmt{n \to \infty}  \hE \Big[ \chi  \big( K^n_{\tau_\ell \land \ga}
    - \n   \wt{K}^\ell_{  \ga} \big) \Big] = 0 $, which shows that
  $\big\{K^n_{\tau_\ell \land \ga} \big\}_{n \in \hN}$ converges weakly to
  $\wt{K}^\ell_{  \ga}$ in $L^2(\cF_T)$.

   \ss \no {\bf (3)}  Now,    let $\ga, \wt{\ga} \n  \in \n   \cT$ such that
  $\ga  \n  \le \n   \wt{\ga}$, \pas ~ For any $n  \n  \in \n   \hN$, since  $K^n$
  is an increasing process,  it holds \pas ~ that
   \bea  \label{eq:a107}
   K^n_{\tau_\ell \land \ga }  \n  \le \n    K^n_{\tau_\ell \land \wt{\ga} } \, .
   \eea
  Then we must have   $ \wt{K}^\ell_{ \ga }  \n  \le \n   \wt{K}^\ell_{ \wt{\ga} } $, \pas:
  {\it Assume not}, i.e.
  the $\hP-$measure of set $A  \n  := \n   \big\{ \wt{K}^\ell_{ \ga }
   \n  > \n   \wt{K}^\ell_{ \wt{\ga} } \big\}
   \n  \in \n   \cF_T$ is  strictly   larger than $0$,   it would follow that
  $\hE \[\b1_A \wt{K}^\ell_{ \ga } \]
    \n  > \n   \hE \[\b1_A \wt{K}^\ell_{ \wt{\ga} } \] $.
   However, we know from part (2) and \eqref{eq:a107} that
   \beas
   \hE \[\b1_A \wt{K}^\ell_{ \ga } \]
   = \lmt{n \to \infty} \hE \[\b1_A K^n_{\tau_\ell \land \ga } \]
    \le  \lmt{n \to \infty} \hE \[\b1_A K^n_{\tau_\ell \land \wt{\ga} } \]
    =  \hE \[\b1_A \wt{K}^\ell_{ \wt{\ga} } \] .
   \eeas
  An contradiction appears. Therefore, $ \wt{K}^\ell_{ \ga } \le \wt{K}^\ell_{ \wt{\ga} } $, \pas ~
  Then Lemma \ref{lem_increasing} shows that $   \wt{K}^\ell $ is an increasing process. \qed

  \ss \no {\bf Proof of \eqref{eq:a189}:}
 Set $\fra \n : = \n  2(\l^+  \n + \n  \k^2)$ and Fix $m, n   \n \in \n  \hN$ with $m  \n > \n  n$.
 We   define processes  $\Xi^{m,n}_t  \n := \n  \Xi^m_t  \n - \n  \Xi^n_t $,
 $t  \n \in \n  [0,T]$ for $\Xi = Y, Y^\ell, Z$.
 Similar to   \eqref{eq:a237},  we can deduce from \eqref{eq:a519}   that   \pas
  \bea
 \q && \hspace{-1.5cm} e^{\fra t  }
  \big|  Y^{\ell,m,n}_{t} \big|^2
 \n  +  \dn   \int_{ t}^{\tau_\ell  }  \n
  e^{\fra s  } \( \fra | Y^{\ell,m,n}_s|^2 \n + \n | Z^{m,n}_s|^2 \)  d s
     \n   =   \n    e^{\fra \tau_\ell      }
    \big| Y^{\ell,m,n}_{ \tau_\ell    } \big|^2
   \n   +  \n   2  \n   \int_{ t}^{\tau_\ell  } \n
    e^{\fra s  }  Y^{\ell,m,n}_s  \( g(s, Y^{\ell,m}_s, Z^m_s)  \n  - \dn   g(s, Y^{\ell,n}_s, Z^n_s) \)   ds  \nonumber    \\
 &  \dn   \dn   &  \dn \n  +    2  \n  \int_{ t}^{\tau_\ell  } \n
     e^{\fra s  } \,  Y^{\ell,m,n}_s  d K^m_s
   \n   -    2 \n   \int_{ t}^{\tau_\ell  } \n
  e^{\fra s  } \, Y^{\ell,m,n}_s  d K^n_s
   - \n   2 \int_{ t}^{\tau_\ell  } \n
    e^{\fra s  } \, Y^{\ell,m,n}_s     Z^{m,n}_s d B_s , \q \fa t \in [\nu, \tau_\ell] .    \label{eq:a145}
   \eea

   By (H1) and (H2), it holds \dsp ~ that
   \bea
  \hspace{-5mm}  Y^{\ell,m,n}_s \n \( g(s, Y^{\ell,m}_s, Z^m_s)  \n  - \dn   g(s, Y^{\ell,n}_s, Z^n_s) \)
    & \tn \dn = & \tn  \dn   Y^{\ell,m,n}_s  \n  \( g(s, Y^{\ell,m}_s, Z^m_s)  \n  - \dn   g(s, Y^{\ell,n}_s, Z^m_s) \)
     \n + \dn  Y^{\ell,m,n}_s  \n  \( g(s,  Y^{\ell,n}_s, Z^m_s)  \n  - \dn   g(s, Y^{\ell,n}_s, Z^n_s) \) \nonumber \\
   & \tn  \dn  \le & \tn  \dn   \l |Y^{\ell,m,n}_s|^2  \n  + \n   \k |Y^{\ell,m,n}_s| |Z^{m,n}_s|
   \n  \le \n  (\l^+  \n  + \n   \k^2) |Y^{\ell,m,n}_s|^2  \n  + \n   \frac14   |Z^{m,n}_s|^2  .  \label{eq:a141}
   \eea
  Also, one can deduce from the definition of process $K^m$ that
   \bea
    \int_t^{\tau_\ell  } \n   e^{\fra s  } \,  Y^{\ell,m,n}_s  d K^m_s
     & \tn \dn  = &  \tn  \dn   \b1_{A_\ell} \n \int_t^{\tau_\ell  } \n   e^{\fra s  } \,  Y^{m,n}_s  d K^m_s
    \n  = \n  \b1_{A_\ell} \n  \int_t^{\tau_\ell  } \n  \b1_{\{Y^m_s < L_s\}}
     e^{\fra s  } \,  Y^{m,n}_s  d K^m_s
      \n   \le  \n   \b1_{A_\ell} \n  \int_t^{\tau_\ell  } \n  \b1_{\{Y^m_s < L_s\}}
     e^{\fra s  } \, ( L_s  \n  - \n   Y^n_s )  d K^m_s    \nonumber \\
       &  \tn   \dn      \le     &  \tn   \dn        e^{\fra T  } \b1_{A_\ell} \n \int_{ \nu  }^{\tau_\ell  } \n
       (    Y^n_s  \n  - \n   L_s )^-  d K^m_s
      \n  \le \n    e^{\fra T  }  \b1_{A_\ell} \n
      \( \underset{s \in [\nu,\tau_\ell]}{\sup} (  Y^n_s  \n  - \n   L_s )^- \)
     K^m_{\tau_\ell  }     , \q \fa t \in [\nu, \tau_\ell] .     \label{eq:a143}
   \eea
   Similarly,
    \beas
   - \n  \int_t^{\tau_\ell  } \n   e^{\fra s  } \,  Y^{\ell,m,n}_s  d K^n_s
      & \tn   \le    &  \tn   \b1_{A_\ell} \n    \int_t^{\tau_\ell  } \n  \b1_{\{Y^n_s < L_s\}}
     e^{\fra s  } \, ( L_s \n  - \n   Y^m_s ) \, d K^n_s
     \n   \le  \n   e^{\fra T  }  \b1_{A_\ell} \n
     \( \underset{s \in [\nu,\tau_\ell]}{\sup} (    Y^m_s  \n  - \n   L_s )^- \)
       K^n_{\tau_\ell  }   \\
        &  \tn     \le    &  \tn      e^{\fra T  }  \b1_{A_\ell} \n
       \( \, \underset{s \in [\nu,\tau_\ell]}{\sup} ( Y^n_s  \n  - \n   L_s )^- \)  K^n_{\tau_\ell  } ,
        \q \fa t \in [\nu, \tau_\ell] .
   \eeas
 Plugging this and \eqref{eq:a141}, \eqref{eq:a143} back into \eqref{eq:a145} shows   that \pas
   \beas
   e^{\fra t  }  \big|  Y^{\ell,m,n}_t \big|^2
 \n  +  \n  \frac12 \n    \int_t^{\tau_\ell  }  \n
  e^{\fra s  }  | Z^{m,n}_s|^2 d s     \le    \eta
    - \n   2 \int_t^{\tau_\ell  } \n
    e^{\fra s  } \, Y^{\ell,m,n}_s     Z^{m,n}_s d B_s , \q \fa t \in [\nu, \tau_\ell] ,
   \eeas
   where $\eta  :=   e^{\fra \tau_\ell      }
    \big| Y^{\ell,m,n}_{ \tau_\ell    } \big|^2
    \n  + \n   2  e^{\fra T  }  \b1_{A_\ell} \n  \( \, \underset{s \in [\nu,\tau_\ell]}{\sup}
    (  Y^n_s  \n  - \n   L_s )^- \)
    \(  K^m_{\tau_\ell  }
     \n  + \n   K^n_{\tau_\ell  }   \)  $.

 Taking expectation for $t=\nu $, we   see from  H\"older's inequality, \eqref{eq:a281} and \eqref{eq:a093} that
 \bea
 \qq && \hspace{-1.5cm} \hE \int_0^T \b1_{\{\nu < t \le \tau_\ell  \}} | Z^{m,n}_t  |^2 dt  \n \le  \n
 \hE  \int_\nu^{\tau_\ell  }  \n
  e^{\fra s  }  | Z^{m,n}_s|^2 d s \n  \le  \n  2 \hE[\eta ] \n  \le  \n
 2 \hE \[    e^{\fra \tau_\ell      }   \big| Y^{\ell,m,n}_{ \tau_\ell    } \big|^2 \]
  \n  + \n   4 e^{\fra T  } \Bigg\{ \hE \[  \b1_{A_\ell}  \underset{s \in [\nu,\tau_\ell]}{\sup}
 \( (  Y^n_s  \n  - \n   L_s )^- \)^2
  \]  \times   \nonumber \\
  & &  \hE \[ \(  K^m_{\tau_\ell  }
     \n  + \n   K^n_{\tau_\ell  }   \)^2  \] \Bigg\}^{1/2}  \n \le \n  C_0 \hE \[
   \b1_{A_\ell} \big| Y_{ \tau_\ell    } \n  - \n  Y^n_{ \tau_\ell    } \big|^2 \]
 + C_0 \ell \left\{ \hE \[  \b1_{A_\ell}  \underset{s \in [\nu,\tau_\ell]}{\sup}
 \( (  Y^n_s  \n  - \n   L_s )^- \)^2
  \] \right\}^{1/2}   .   \label{eq:a149}
 \eea
 On the other hand, the Burkholder-Davis-Gundy inequality implies that
\beas
 \q  && \hspace{-1cm} \hE \[ \underset{t \in [\nu,\tau_\ell]}{\sup}
  \big|  Y^{\ell,m,n}_t \big|^2 \]
    \n   \le     \n
 \hE \[ \underset{t \in [\nu,\tau_\ell]}{\sup}   e^{\fra t  }
  \Big|  Y^{\ell,m,n}_t \Big|^2 \]  \le
  \hE  [ \eta   ] \n  + \n  2 \hE \[ \underset{t \in [0,T]}{\sup}
  \bigg| \int_t^T  \n  \b1_{\{\nu < s \le  \tau_\ell \}}
    e^{\fra s  } \, Y^{\ell,m,n}_s     Z^{m,n}_s d B_s \bigg| \] \\
   & &   \le    \n    \hE  [ \eta   ] \n  + \n  C_0 \hE
     \[ \( \underset{t \in [\nu,\tau_\ell]}{\sup} \big| Y^{\ell,m,n}_t \big| \) \n  \cd \n
     \(  \int_\nu^{\tau_\ell  } e^{\fra t  } \,  | Z^{m,n}_t  |^2 dt \)^{1/2} \]
     \n   \le   \n   \hE  [ \eta   ] \n  + \n
    \frac12 \hE \[ \underset{t \in [\nu,\tau_\ell]}{\sup} \big|  Y^{\ell,m,n}_t \big|^2 \]
    \n  + \n   C_0 \hE   \int_\nu^{\tau_\ell  } e^{\fra t  } \,  | Z^{m,n}_t  |^2 dt   .
\eeas
 As $ \hE \[ \underset{t \in [\nu,\tau_\ell]}{\sup}
 \big|   Y^{\ell,m,n}_t \big|^2 \] \le 4 \ell^2 $ by \eqref{eq:a281} and \eqref{eq:a095},
 it follows from \eqref{eq:a149} that
 \bea  \label{eq:a283}
 \hE \[ \underset{t \in [\nu,\tau_\ell]}{\sup}  \big|  Y^{\ell,m,n}_t \big|^2 \] \le 2  \hE  [ \eta   ]
 \n  + \n   C_0 \hE   \int_\nu^{\tau_\ell  } e^{\fra t  } \,  | Z^{m,n}_t  |^2 dt \le C_0  \hE  [ \eta   ]   .
 \eea

 Since Doob's martingale inequality and \eqref{eq:a281} show  that
  \beas
 \hE \[ \underset{t \in [0,\nu]}{\sup}   \Big|  Y^{\ell,m,n}_t \Big|^2 \]
  \n \le \n  \hE \[ \underset{t \in [0,T]}{\sup}    \Big| \hE \big[  \b1_{A_\ell}  Y^{m,n}_\nu \big| \cF_t  \big] \Big|^2 \]
  \n  \le  \n  
   4 \hE \[      \big|  \b1_{A_\ell}  Y^{m,n}_\nu   \big|^2 \]
  \n \le \n  4 \hE \[ \underset{t \in [\nu,\tau_\ell]}{\sup}   \big|  Y^{\ell, m,n}_t \big|^2 \] ,
 \eeas
 we see from \eqref{eq:a281} and \eqref{eq:a283} that
 \beas
 \hE \[ \underset{t \in [0,T]}{\sup}   \Big|  Y^{\ell,m,n}_t \Big|^2 \]
 \n \le \n  \hE  \n \[ \underset{t \in [0,\nu]}{\sup}   \Big|  Y^{\ell,m,n}_t \Big|^2   \n + \n
   \underset{t \in [\nu,\tau_\ell]}{\sup}   \Big|  Y^{\ell,m,n}_t \Big|^2 \]
  \n \le \n  5 \hE \[ \underset{t \in [\nu,\tau_\ell]}{\sup}   \big|  Y^{\ell,m,n}_t \big|^2 \]
  \n \le \n  C_0  \hE  [ \eta   ]     .
 \eeas
 This together with   \eqref{eq:a149} leads to that
  \beas
  \hspace{-6mm}   \underset{m > n}{\sup} \Bigg\{  \hE \n \[ \underset{t \in [0,T]}{\sup}  \Big|  Y^{\ell,m,n}_t \Big|^2 \]
   \dn  + \n  \hE  \dn  \int_0^T  \n  \b1_{\{\nu < t \le \tau_\ell  \}} | Z^{m,n}_t  |^2 dt  \Bigg\}
    \n \le \n C_0  \hE  [ \eta   ]     \n   \le  \n    C_0 \hE \n \[
        \b1_{A_\ell}    \big|  Y_{ \tau_\ell    } \dn  - \n  Y^n_{ \tau_\ell    }  \big|^2 \]
   \n   +   \n    C_0 \ell \left\{ \hE \n \[   \b1_{A_\ell}   \underset{t \in [\nu,\tau_\ell]}{\sup}
 \( (  Y^n_t  \dn  - \n   L_t )^- \)^2   \] \right\}^{1/2} .
  \eeas
  As $  \b1_{A_\ell}  \big|  Y_{ \tau_\ell   } \n  - \n  Y^n_{ \tau_\ell   } \big|
   \n \le \n  2 \ell $, $ \fa n  \n \in \n  \hN$ by \eqref{eq:a095}, letting $n  \n \to \n  \infty$,
   we see from   bounded convergence theorem and \eqref{eq:a151} that
  \beas
  \lmt{n \to \infty} \,  \underset{m > n}{\sup} \left\{  \hE \[ \underset{t \in [0,T]}{\sup}
  \big|  Y^{\ell,m,n}_t  \big|^2 \]
   \n  + \n     \hE  \int_0^T \b1_{\{\nu < t \le \tau_\ell  \}}     | Z^{m,n}_t  |^2 dt  \right\}   = 0 .
  \eeas
  Hence, $  \big\{ Y^{\ell,n} \big\}_{n   \in   \hN} $  is Cauchy sequence in
  $\hS^2$ and $ \big\{ \b1_{\{\nu < t \le \tau_\ell  \}} Z^n_t \big\}_{t \in [0,T]} $,
  $n \n \in \n \hN$ is Cauchy sequence in   $ \hH^{2,2} $.   \qed

  \ss \no {\bf Proof of \eqref{eq:a179}:}  As $K^n_\nu \n =  \n  \cK^\ell_\nu  \n = \n  0 $ by \eqref{eq:a245}
  and \eqref{eq:xax014},   one can deduce from  \eqref{eq:a333} that \pas
  \beas
  \hspace{-0.2cm}   K^n_t \n   -  \n   \cK^\ell_t
      \n = \n  ( K^n_t \n   -  \n K^n_\nu ) \n   -  \n ( \cK^\ell_t \n - \n \cK^\ell_\nu )
  &  \dn \dn  =  &  \dn \dn   Y^{\ell,n}_\nu \n - \n Y^{\ell,n}_t \n  - \n \b1_{A_\ell} \( Y_\nu  \n  - \n       Y_t \)
    \n  - \n   \int_\nu^t     \( g(s, Y^{\ell,n}_s, Z^n_s)  \n  - \n   g(s,  Y_s, \cZ^\ell_s) \) ds \\
   &  \dn \dn   &  \dn \dn     + \n  \int_\nu^t     (Z^n_s  \n  - \n   \cZ^\ell_s)  d B_s,
   \q  \fa t  \n  \in \n   [\nu,\tau_\ell] .
  \eeas
   Then \eqref{eq:a281} and (H1) show that \pas
  \beas
    \hspace{-0.3cm}   \big|  K^n_t \n   -  \n   \cK^\ell_t  \big|
   \n   \le   \n  \b1_{A_\ell}  |Y^n_\nu  \dn  - \n   Y_\nu|  \n  + \n  \b1_{A_\ell}  \big|   Y^n_t  \dn  - \n  Y_t \big|
   \n   + \dn   \int_\nu^t  \dn   \( \k |Z^n_s  \n  - \dn   \cZ^\ell_s | \n  +  \n
     | g(s, Y^{\ell,n}_s, \cZ^\ell_s)  \n  - \dn   g(s, Y_s, \cZ^\ell_s)  | \) \n ds
   \n     +   \n    \bigg| \int_\nu^t \dn
    \( Z^n_s  \n  - \dn   \cZ^\ell_s \) d B_s \bigg| , ~ \fa t  \n  \in \n   [\nu,\tau_\ell] .
  \eeas
  Since H\"older's inequality and \eqref{eqn-d011} imply  that
  \beas
   \big|  K^n_t \n   -  \n   \cK^\ell_t  \big|^2
   & \tn \dn   \le  & \tn  \dn     C_0 \b1_{A_\ell}  |Y^n_\nu  \n  - \n   Y_\nu|^2
    \n  + \n    C_0 \b1_{A_\ell}   \big|  Y^n_t  \n  - \n    Y_t \big|^2
   \n   + \n  C_0 \int_\nu^t \n   | Z^n_s \n  - \n  \cZ^\ell_s |^2 ds  \n  +  \n
   C_0  \( \b1_{A_\ell} \n \int_\nu^t \n
     | g(s, Y^{\ell,n}_s, \cZ^\ell_s)  \n  - \n   g(s, Y_s, \cZ^\ell_s)  |   ds \)^2    \\
     & \tn  \dn  & \dn  \dn  + C_0 \,  \underset{\wt{t} \in [0,T]}{\sup}
   \bigg| \int_0^{\wt{t}} \n \b1_{\{ \nu < s \le  \tau_\ell \}} \( Z^n_s  \n  - \n   \cZ^\ell_s \) d B_s \bigg|^2  ,
   \q  \fa  t  \n  \in \n   [\nu,\tau_\ell] ,
  \eeas
   we can deduce from  Doob's martingale inequality   that
  \beas
 \hE \[ \underset{t \in [\nu,\tau_\ell]}{\sup} \big| K^n_t - \cK^\ell_t \big|^2 \]
 & \tn   \le & \tn  C_0 \hE \[  \b1_{A_\ell}  |Y^n_\nu - Y_\nu |^2 \]  \n   + \n  C_0 \hE
 \[  \b1_{A_\ell}  \underset{t \in [\nu,\tau_\ell]}{\sup}  \big|  Y^n_t  -  Y_t \big|^2  \]
  \n   + \n  C_0 \hE \n  \int_\nu^{\tau_\ell    } \n
  |Z^n_t \n  - \n  \cZ^\ell_t|^2 dt      \\
   & \tn  & \tn    +  \,  C_0 \hE \[ \(     \int_\nu^{\tau_\ell  } \n    \b1_{A_\ell}
   \big| g(t, Y^n_t, \cZ^\ell_t)  \n  - \n   g(t, Y_t, \cZ^\ell_t)  \big|  dt \)^2 \]  .
 \eeas
   \if{0}
  \beas
  \hE \[ \underset{t \in [0,T]}{\sup} \big| Y^n_{\tau_\ell \land t} - \cY^\ell_{\tau_\ell \land t}  \big|^2 \]
  =  \hE \[ \underset{t \in [0,\tau_\ell]}{\sup} \big| Y^n_t - \cY^\ell_t  \big|^2 \]
  =  \hE \[ \underset{t \in [0,\tau_\ell]}{\sup} \big| Y^n_{\tau_\ell \land t} - \cY^\ell_t  \big|^2 \]
  \le  \hE \[ \underset{t \in [0,T]}{\sup} \big| Y^n_{\tau_\ell \land t} - \cY^\ell_t  \big|^2 \]
  \eeas

  \beas
  \hE  \int_0^{\tau_\ell  }  \n    |   Z^n_t - \cZ^\ell_t  |^2 dt
   =   \hE  \int_0^{\tau_\ell  }  \n    | \b1_{\{ t \le \tau_\ell \}} Z^n_t - \cZ^\ell_t  |^2 dt
   \le  \hE  \int_0^T  \n    | \b1_{\{ t \le \tau_\ell \}} Z^n_t - \cZ^\ell_t  |^2 dt
  \eeas
  \fi
   The bounded convergence theorem and \eqref{eq:a095} imply  that
   $ \lmtd{n \to \infty} \hE \[  \b1_{A_\ell}  |Y^n_\nu \n - \n Y_\nu|^2 \] \n = \n 0 $.
   Thanks to \eqref{eq:a341}, it remains to show that
  \bea  \label{eq:a167}
  \lmt{n \to \infty} \hE \[ \(    \int_\nu^{\tau_\ell  } \n \b1_{A_\ell}
   \big| g \(t,  Y^n_t  , \cZ^\ell_t   \)
   \n  - \n   g \(t,   Y_t , \cZ^\ell_t \)  \big|  dt \)^2 \]  = 0 .
  \eea

    By (H3), it holds \dtp ~ that
    $   \lmt{n \to \infty} \b1_{A_\ell \cap \{\nu < t \le \tau_\ell\}}
   \big| g(t, Y^n_t, \cZ^\ell_t)  \n  - \n   g(t, Y_t, \cZ^\ell_t)  \big| \n = \n 0 $.
   Also,  (H1), (H4)  and \eqref{eq:a095} imply that for any $n \in \hN$
  \bea
    && \hspace{-1.5cm}  \b1_{A_\ell \cap \{\nu < t \le \tau_\ell \}} \big| g(t, Y^n_t, \cZ^\ell_t)
     \n  - \dn   g(t, Y_t, \cZ^\ell_t)  \big|
    \n  \le     \n    \b1_{A_\ell \cap \{\nu < t \le \tau_\ell \}}
    \Big( \big| g(t, Y^n_t, 0) \big| \n  + \n
    \big| g(t, Y^n_t, \cZ^\ell_t) \n - \dn g(t, Y^n_t, 0) \big|
  \n   + \n  \big| g(t, Y_t, 0) \big|     \nonumber \\
    & &       +      \big| g(t, Y_t, \cZ^\ell_t) \n - \dn g(t, Y_t, 0) \big|   \Big)
   \n \le \n   \b1_{A_\ell \cap \{\nu < t \le \tau_\ell \}}  \(  2   h_t \n  +\n   2 \k \ell
   \n  +\n   2 \k   |\cZ^\ell_t| \) : = \fh^\ell_t   , ~ \dtp  \qq  \q  \label{eq:a163}
  \eea
  As    $  \hE \int_0^T  \fh^\ell_t \, dt  \n \le  \n   2\ell \n  + \n  2 \k \ell  T  \n  +\n
     2 \k T^{1/2} \( \hE \int_\nu^{\tau_\ell} |\cZ^\ell_t|^2 dt  \)^{1/2}  \n < \n  \infty $
     by \eqref{eq:a273} and  H\"older's inequality,
     applying   the dominated convergence theorem yields that
   \beas
   \lmt{n \to \infty} \hE \int_\nu^{\tau_\ell}  \n \b1_{A_\ell}
   \big|g(t, Y^n_t, \cZ^\ell_t) \n - \n g(t,Y_t,\cZ^\ell_t) \big| dt = 0 .
   \eeas
  So up to a subsequence of  $\{ Y^n \}_{n \in \hN} $, it holds \pas ~ that
    $ \lmt{n \to \infty} \n \int_\nu^{\tau_\ell}  \n \b1_{A_\ell}
     \big|g(t, Y^n_t, \cZ^\ell_t) \n - \dn  g(t,Y_t,\cZ^\ell_t) \big| dt
    \n = \n  0 $. Since  \eqref{eq:a163}  shows   that for any $n \n \in \n  \hN$,
  $   \(  \int_\nu^{\tau_\ell    } \n  \b1_{A_\ell}
   \big| g(t, Y^n_t, \cZ^\ell_t)  \n  - \dn   g(t, Y_t, \cZ^\ell_t)  \big|    dt \)^2  \n \le \n
  \(  \int_0^T \n  \fh^\ell_t  dt \)^2   $,    \pas ~ and since  H\"older's inequality implies that
   $  \hE \n \[ \(  \int_0^T \n  \fh^\ell_t  dt \)^2 \]  \dn \le \n
   \hE  \n  \[ \big(   2 \ell \n  + \n  2 \k \ell  T
   \n  + \n   2 \k   \int_\nu^{\tau_\ell}  \n  |\cZ^\ell_t|  dt \big)^2 \]
    \dn \le \n  C_0   \ell^2   \n + \n  C_0 \hE    \n    \int_\nu^{\tau_\ell} \n  |\cZ^\ell_t|^2  dt   \n < \n  \infty $,
  applying the dominated convergence theorem again yields \eqref{eq:a167}.    \qed

  \ss \no {\bf Proof of \eqref{eq:a199}:}
   \ss  For any $n \in \hN$, H\"older's inequality and  \eqref{eq:a093} imply that
  \beas
   \hE  \n  \int_\nu^{\tau_\ell} \dn \big| Y^n_t \dn - \n Y_t \big| d  K^n_t
 & \dn \dn   = & \dn \dn  \hE \[ \b1_{A_\ell} \n  \int_\nu^{\tau_\ell} \dn \big| Y^n_t \dn - \n Y_t \big| d  K^n_t \]
    \n  \le   \n   \hE \n \[  \b1_{A_\ell}
   \n \( \underset{t \in [\nu,\tau_\ell]}{\sup}
  \big| Y^n_t \dn - \n Y_t \big| \)  K^n_{\tau_\ell} \] \\
  & \dn \dn  \le & \dn \dn   \Bigg\{ \hE \n \[ \( K^n_{\tau_\ell}   \)^2 \]
   \hE \n \[  \b1_{A_\ell} \underset{t \in [\nu,\tau_\ell]}{\sup} \big| Y^n_t \n - \n Y_t \big|^2\] \dn \Bigg\}^{1/2}
    \n  \le  \n   C_0   \ell
    \Bigg\{ \hE \n \[ \b1_{A_\ell} \underset{t \in [\nu,\tau_\ell]}{\sup}
    \big| Y^n_t \n - \n Y_t \big|^2\] \dn \Bigg\}^{1/2} .
  \eeas
  As $n \to \infty$, \eqref{eq:a341} shows that
  $ \dis \lmt{n \to \infty} \hE \[ \bigg| \n  \int_\nu^{\tau_\ell} \dn \big( Y^n_t \dn
  - \n Y_t \big) d  K^n_t \bigg| \] = 0 $. So up to a subsequence of $\{Y^n\}_{n \in \hN}$,
    \bea   \label{eq:a197}
     \lmt{n \to \infty}    \int_\nu^{\tau_\ell} \dn \big( Y^n_t \n  - \n Y_t \big) d  K^n_t   = 0 , \q \pas
    \eea
  For \pas ~ $\o \in \O$ such that \eqref{eq:a195} holds and that  path
  $  Y_t (\o)  \n - \n  L_t (\o)  $ is continuous
  from $t \n = \n \nu (\o) $ to   $ t \n = \n   \tau_\ell(\o)   $ by \eqref{eq:a205},
  one can   deduce from \eqref{eq:a195} that  measure $d K^n_t (\o) $ converges weakly to
   measure $d \cK^\ell_t (\o) $ on    period  $   \[ \nu(\o) ,   \tau_\ell(\o) \] $, so
   \beas
   \lmt{n \to \infty} \int_{\nu(\o)}^{\tau_\ell (\o)} \dn \( Y_t (\o) \n  - \n L_t (\o) \) d  K^n_t (\o)
   = \int_{\nu(\o)}^{\tau_\ell (\o)} \dn \( Y_t (\o) \n  - \n L_t (\o) \) d  \cK^\ell_t (\o) .
   \eeas
   Adding this to   \eqref{eq:a197}, we see from  \eqref{eq:a181}  that  \pas
    \beas
    0 & \tn \le & \tn  \b1_{\{ \nu < \tau_\ell \}} \int_t^{\tau_\ell}  \n  (Y_s \n - \n L_s) d \cK^\ell_s
    \le  \b1_{\{ \nu < \tau_\ell \}} \int_\nu^{\tau_\ell}  \n  (Y_s \n - \n L_s) d \cK^\ell_s
       \n = \n  \int_\nu^{\tau_\ell}  \n  (Y_s \n - \n L_s) d \cK^\ell_s
       \n = \n  \lmt{n \to \infty} \n \int_\nu^{\tau_\ell} \n \big( Y^n_s \n  - \n L_s \big) d  K^n_s \\
        & \tn  =  & \tn   \lmt{n \to \infty}  \n  \int_\nu^{\tau_\ell}  \n  \b1_{\{Y^n_s < L_s \}}
        (Y^n_s \n - \n L_s) d  K^n_s  \n \le \n  0 , \q \fa t \in [\nu, \tau_\ell] ,
    \eeas
   proving \eqref{eq:a199}.   \qed

    \ss \no {\bf Proof of Claim \eqref{eq:b151}:} It is clear that
 $ Y_{\ga_\nu}  = \b1_{\{\ga_\nu = T\}}  Y_T + \b1_{\{\ga_\nu < T\}}  Y_{\ga_\nu}
  \le \b1_{\{\ga_\nu = T\}}  \xi + \b1_{\{\ga_\nu < T\}}  U_{\ga_\nu}   $,   \pas,
  so we only need to show the converse inequality.

 \ss   Fix $n \in \hN$. Clearly,
 $K^n_s  \n : = \n  n \int_0^s (Y^n_r  \n - \n  L_r)^- dr $, $s  \n \in \n  [0,T]$
  is   a process of $\hK^0$   satisfying that \pas
    \beas
  Y^n_t =   Y^n_{\ga^n_\nu}
  + \int_t^{\ga^n_\nu} g(s,Y^n_s, Z^n_s   )  ds + K^n_{\ga^n_\nu} - K^n_\nu
  - \int_t^{\ga^n_\nu} Z^n_s d B_s, \q   \fa  t \in  [ \nu ,\ga^n_\nu   ]
  \eeas
  by \eqref{eq:b117}.
  Since $\hE[|Y^n_\nu|] \n < \n \infty$ by the uniform integrability of $\{Y^n_\z \}_{\z \in \cT}$,
   applying Lemma \ref{lem_RBSDE_estimate} with $(Y,Z,K)  \n = \n  (Y^n,Z^n,K^n) $ and $ \tau  \n = \n  \ga^n_\nu $,
   we see from \eqref{eqn-d011b}   that  for any $p \n \in \n (0,1)$
  \bea  \label{eq:b137}
  \hE \n \[ \( \int_\nu^{\ga^n_\nu} \n  | Z^n_t |^2 dt \)^{p/2} \]
   \dn \le \n  C_p \hE \n  \[\underset{t \in [\nu, \ga^n_\nu]}{\sup} |Y^n_t|^p
   \n +  \n    \bigg(\int_\nu^{\ga^n_\nu} \n  h_t d t \bigg)^p\]
      \n \le \n  C_p \hE \n  \[ 1 \n + \n \underset{s \in [0, T]}{\sup} |Y^1_s|^p    \n + \n
   \underset{s \in [0, T]}{\sup} |Y_s|^p
    \n + \dn     \int_0^T \n  h_t d t \]  .
  \eea

 Let $j  \n \in \n  \hN$ and define a stopping time
  $\z^n_j  \n  := \n \inf \big\{ t  \n \in \n  [0,T] \n : \int_0^t \n  |Z^n_s|^2 ds
   \n > \n  j \big\}  \n \land \n  T  \n \in \n  \cT $.
   Since \eqref{eq:b117} shows that
     \beas
  Y_{\ga_\nu \land \z^n_j} & \tn \ge  & \tn  Y^n_{\ga_\nu \land \z^n_j} = Y^n_{\ga^n_\nu \land \z^n_j}
  + \int_{\ga_\nu \land \z^n_j}^{\ga^n_\nu \land \z^n_j} g(s,Y^n_s, Z^n_s   )  ds
  + K^n_{\ga^n_\nu \land \z^n_j} - K^n_{\ga_\nu \land \z^n_j}
  - \int_{\ga_\nu \land \z^n_j}^{\ga^n_\nu \land \z^n_j} Z^n_s d B_s \\
   & \tn  \ge & \tn  Y^n_{\ga^n_\nu \land \z^n_j}
   + \int_{\ga_\nu \land \z^n_j}^{\ga^n_\nu \land \z^n_j} g(s,Y^n_s, Z^n_s   )  ds
   - \int_{\ga_\nu \land \z^n_j}^{\ga^n_\nu \land \z^n_j} Z^n_s d B_s   , \q \pas ,
  \eeas
 taking conditional expectation $\hE\[ ~ \cd ~ \Big|\cF_{\ga_\nu \land \z^n_j}\]$ yields that \pas
 \bea
 \hspace{-2mm}  Y_{\ga_\nu \land \z^n_j} & \tn \dn \ge  & \tn  \dn  \hE \[   Y^n_{\ga^n_\nu \land \z^n_j}
   + \int_{\ga_\nu \land \z^n_j}^{\ga^n_\nu \land \z^n_j}  \n
   g(t,Y^n_t, Z^n_t   )  dt \bigg| \cF_{\ga_\nu \land \z^n_j} \]
   \n = \n  \b1_{\{\ga_\nu \ge \z^n_j\}} \hE \n \[   Y^n_{\ga^n_\nu \land \z^n_j}
   \n  + \dn  \int_{\ga_\nu \land \z^n_j}^{\ga^n_\nu \land \z^n_j}  \n
   g(t,Y^n_t, Z^n_t   )  dt \bigg| \cF_{  \z^n_j} \] \nonumber \\
  & \tn  \dn     & \tn  \dn \qq \q
  \n +      \b1_{\{\ga_\nu < \z^n_j\}} \hE \n  \[   Y^n_{\ga^n_\nu \land \z^n_j}
  \n + \dn  \int_{\ga_\nu \land \z^n_j}^{\ga^n_\nu \land \z^n_j}  \n
  g(t,Y^n_t, Z^n_t   )  dt \bigg| \cF_{\ga_\nu  } \] \n : = \n I^{n,j}_1 \n + \n I^{n,j}_2  . \label{eq:b131}
 \eea
 As $\{\ga_\nu \ge \z^n_j\} \subset \{\ga^n_\nu \ge \z^n_j\}$,  it holds \pas ~  that
  \bea
 I^{n,j}_1  \n  =    \n     \hE \[  \b1_{\{\ga_\nu \ge \z^n_j\}} Y^n_{\ga^n_\nu \land \z^n_j}
    \n + \n  \b1_{\{\ga_\nu \ge \z^n_j\}}  \n \int_{\ga_\nu \land \z^n_j}^{\ga^n_\nu \land \z^n_j} \n
    g(t,Y^n_t, Z^n_t   )  dt \bigg| \cF_{  \z^n_j} \]
            \n = \n         \hE \[  \b1_{\{\ga_\nu \ge \z^n_j\}} Y^n_{ \z^n_j}  \big| \cF_{  \z^n_j} \]
           \n = \n        \b1_{\{\ga_\nu \ge \z^n_j\}} Y^n_{ \z^n_j}   .  \label{eq:b133}
 \eea

 Similar to \eqref{eq:b123},  (H4), (H5), \eqref{eqn-d011} and \eqref{eqn-d011b} imply that
 $ \big| g(t,Y^n_t, Z^n_t   ) \big| \n \le \n  \k  \n + \n  (1 \n + \n \k) h_t  \n +    2 \k |Y^n_t |
     \n + \n  \k| Z^n_t |^\a $,   \dtp ~ \;
  It then follows from H\"older's inequality that  \pas
 \bea
   \int_{\ga_\nu \land \z^n_j}^{\ga^n_\nu \land \z^n_j} \n
   \big| g(s,Y^n_s \n , Z^n_s   ) \big|  ds
   & \tn \dn \le   & \tn  \dn   \dn   \int_{\ga_\nu }^{\ga^n_\nu }  \n
   \big| g(s,Y^n_s \n , Z^n_s   ) \big|  ds
       \n    \le  \n    C_0 \n   \int_{\ga_\nu }^{\ga^n_\nu } \n \( 1  \n  + \n     h_s
     \n  + \n     |Y^n_s |    \) ds
       \n  + \n   \k  (\ga^n_\nu  \n  - \n   \ga_\nu)^{1-\a/2}
       \( \int_{\ga_\nu }^{\ga^n_\nu } \n   | Z^n_s |^2   ds \)^{\a/2}  \qq \q  \label{eq:b126} \\
   & \tn    \le   & \tn     C_0 \n   \int_0^T \n \( 1  \n  + \n     h_s
     \n  + \n     |Y^n_s |    \) ds
       \n  + \n  C_\a
       \( \int_0^T \n   | Z^n_s |^2   ds \)^{\a/2} .     \label{eq:b127}
 \eea
 By Fubini's Theorem and the uniform integrability of $\{Y^n_\z \}_{\z \in \cT}$,
 $ \hE \n \int_0^T  \n  |Y^n_s | ds  \n = \n  \int_0^T  \n  \hE \big[ |Y^n_s |\big] ds
  \n \le \n  T \underset{s \in [0,T]}{\sup} \hE \big[ |Y^n_s |\big]  \n < \n  \infty $,
 which together with $Z^n \in \hH^{2,\a}$
 shows that the last term in \eqref{eq:b127} is integrable.
 As $Z^n \in \underset{p \in (0,1)}{\cap} \hH^{2,p} \subset \hH^{2,0} $ shows that
  $ \big\{ \z^n_j  \big\}_{j   \in   \hN}$ is stationary,
  it holds \pas ~ that $ \lmt{j \to \infty} Y_{\ga_\nu \land \z^n_j} \n = \n  Y_{\ga_\nu} $
  though we have not yet shown whether $Y$ is a continuous process.
  Letting $j \to \infty$ in \eqref{eq:b131} and \eqref{eq:b133},
  we  can deduce from  the uniform integrability of $\{Y^n_\z \}_{\z \in \cT}$
 and the conditional-expectation version of  dominated convergence theorem  that
 \bea
  Y_{\ga_\nu} & \tn \ge & \tn  \b1_{\{\ga_\nu = T\}} Y^n_T +
 \lmt{j \to \infty} I^{n,j}_2 = \b1_{\{\ga_\nu = T\}} \xi + \b1_{\{\ga_\nu < T \}} \hE \n  \[   Y^n_{\ga^n_\nu  }
  \n + \dn  \int_{\ga_\nu  }^{\ga^n_\nu  }  \n
  g(t,Y^n_t, Z^n_t   )  dt \bigg| \cF_{\ga_\nu  } \] \nonumber \\
   & \tn  =  & \tn  \b1_{\{\ga_\nu = T\}} \xi + \b1_{\{\ga_\nu < T \}}
   \hE \n  \[ \b1_{\{\ga^n_\nu=T\}} \xi \n + \n \b1_{\{\ga^n_\nu<T\}}  U_{\ga^n_\nu  }
   \n + \n  \int_{\ga_\nu  }^{\ga^n_\nu  }  \n
  g(t,Y^n_t, Z^n_t   )  dt \bigg| \cF_{\ga_\nu  } \]  , \q \pas , \label{eq:b141}
 \eea
 where we used in the last equality the fact that $Y^n_{\ga^n_\nu} = U_{\ga^n_\nu}$, \pas ~ on
 $ \{ \ga^n_\nu < T \} $ by the continuity of $Y^n$ and $U$.

 \ss  Since $\lmtd{n \to \infty}  \b1_{\{\ga^n_\nu=T\}} =  \b1_{\{\ga_\nu=T\}}  $
 and since $\xi  \in L^1(\cF_T)   $,
 applying the conditional-expectation version of  dominated convergence theorem  yields that
 \bea  \label{eq:b143}
 \lmt{n \to \infty} \b1_{\{\ga_\nu < T \}}
  \hE \[   \b1_{\{\ga^n_\nu=T\}} \xi \big| \cF_{\ga_\nu  } \] =
 \lmt{n \to \infty}
  \hE \[   \b1_{\{\ga_\nu < T \}} \b1_{\{\ga^n_\nu=T\}} \xi \big| \cF_{\ga_\nu  } \] =  0 , \q \pas
 \eea
 As $\b1_{\{\ga^n_\nu<T\}} \big| U_{\ga^n_\nu  } \big| \n = \n  \b1_{\{\ga^n_\nu<T\}}  \big|  Y^n_{\ga^n_\nu  } \big|
 \n \le \n  \big|  Y^1_{\ga^n_\nu  } \big|  \n + \n  \big|  Y_{\ga^n_\nu  } \big|$, \pas,
 the uniform integrability of $ \{ Y^1_\z \}_{\z \in \cT} $ and $ \{ Y_\z \}_{\z \in \cT} $ implies that
 of $ \big\{ \b1_{\{\ga^n_\nu<T\}}   U_{\ga^n_\nu  }   \big\}_{n \in \hN}$,
 and it then  follows from the continuity of $U$ that
 \bea   \label{eq:b145}
 \lmt{n \to \infty} \hE \[   \b1_{\{\ga^n_\nu<T\}}  U_{\ga^n_\nu  } \big| \cF_{\ga_\nu  } \]
 = \hE \[   \b1_{\{\ga_\nu<T\}}  U_{\ga_\nu  } \big| \cF_{\ga_\nu  } \]
 = \b1_{\{\ga_\nu<T\}}  U_{\ga_\nu  } , \q \pas
 \eea

 Set $\wt{\a}:=\frac12 (1+\a) \in (0,1)$. Given   $\e > 0$,  with
 $A^n_\e \n : = \n  \Big\{ \hE \Big[  \int_{\ga_\nu  }^{\ga^n_\nu  }
 \big|  g(s,Y^n_s, Z^n_s   ) \big|  ds \Big| \cF_{\ga_\nu  } \Big]  \n > \n  \e \Big\}  \n \in \n  \cF_{\ga_\nu  } $,
    \eqref{eq:b126}, H\"older's inequality and \eqref{eq:b137} imply that
 \beas
 \hspace{-3mm}  \hP ( A^n_\e ) & \tn \dn \le  & \tn \dn
 \frac{1}{\e} \hE \[ \b1_{A^n_\e} \hE \bigg[  \int_{\ga_\nu  }^{\ga^n_\nu  }
 \big|  g(t,Y^n_t, Z^n_t   ) \big|  dt \bigg| \cF_{\ga_\nu  } \bigg]  \]
 = \frac{1}{\e}  \hE \[ \b1_{A^n_\e} \int_{\ga_\nu  }^{\ga^n_\nu  }  \big| g(t,Y^n_t, Z^n_t   ) \big| dt   \]  \\
   & \tn \dn    \le  & \tn \dn  \frac{C_0}{\e}    \hE \n   \int_{\ga_\nu }^{\ga^n_\nu } \n \( 1  \n  + \n     h_t
     \n  + \n     |Y^n_t |   \) dt
       \n  + \n  \frac{\k}{\e}  \hE \[ (\ga^n_\nu  \n  - \n   \ga_\nu)^{1-\a/2}
       \( \int_{\ga_\nu }^{\ga^n_\nu } \n   | Z^n_t |^2   dt \)^{\a/2}  \]    \\
   & \tn \dn    \le  & \tn \dn  \frac{C_0}{\e} \hE \n   \int_{\ga_\nu }^{\ga^n_\nu } \n \( 1  \n  + \n     h_t
     \n  + \n     |Y^1_t | \n  + \n     |Y_t |   \) dt   \n  + \n
    \frac{\k}{\e}  \Big\{ \hE \[  (\ga^n_\nu -\ga_\nu )^{\frac{(2-\a)\wt{\a}}{2(\wt{\a}-\a)}} \] \Big\}^{1-\a/\wt{\a}}
     \Bigg\{ \hE \Bigg[ \bigg( \int_{\ga_\nu  }^{\ga^n_\nu  }
   |    Z^n_t    |^2 dt \bigg)^{\wt{\a}/2}   \Bigg] \Bigg\}^{\a/\wt{\a}} \\
   & \tn \dn    \le & \tn \dn  \n \frac{C_0}{\e} \hE \n   \int_{\ga_\nu }^{\ga^n_\nu } \n \( 1  \n  + \n     h_t
     \n  + \n     |Y^1_t | \n  + \n     |Y_t |   \) dt   \n  + \n \frac{C_\a}{\e}
     \Big\{ \hE \[  (\ga^n_\nu  \n - \n \ga_\nu )^{\frac{(2-\a)\wt{\a}}{2(\wt{\a}-\a)}} \] \Big\}^{1-\a/\wt{\a}}
     \Bigg\{ \hE \n  \[ 1 \n + \n \underset{s \in [0, T]}{\sup} |Y^1_t|^{\wt{\a}}     \n + \n
   \underset{s \in [0, T]}{\sup} |Y_t|^{\wt{\a}}
    \n + \dn     \int_0^T \n  h_t d t \] \Bigg\}^{\a/\wt{\a}} .
 \eeas
 Since the Fubini's Theorem and the uniform integrability of $\{Y^1_\z\}_{\z \in \cT}$, $\{Y_\z\}_{\z \in \cT}$
 show  that
    \beas
    \q  \hE \n \int_0^T \dn \( 1  \n  + \n     h_t
     \n  + \n     |Y^1_t | \n  + \n     |Y_t |   \) dt \n \le \dn
     \int_0^T \n  ( 1  \n  + \n     h_t ) dt \n + \dn  \int_0^T \n  \hE \[ |Y^1_t | \n  + \n     |Y_t |\] dt
     \n \le \dn
     \int_0^T \n  ( 1  \n  + \n     h_t ) dt \n + \n  T   \underset{t \in [0,T]}{\sup} \hE [ |Y^1_t | ]
     \n + \n  T \underset{t \in [0,T]}{\sup} \hE [ |Y_t | ] \n < \n  \infty  ,
     \eeas
   letting $n \to \infty$, we can deduce from
     the dominated convergence theorem and the bounded convergence theorem   that
     \beas
    \lmt{n \to \infty} \, \hP \Bigg\{ \hE \bigg[  \int_{\ga_\nu  }^{\ga^n_\nu  }
 \big|  g(s,Y^n_s, Z^n_s   ) \big|  ds \bigg| \cF_{\ga_\nu  } \bigg]  \n > \n  \e \Bigg\} = 0 , \q \pas
     \eeas
 Thus, $\hE \[  \int_{\ga_\nu  }^{\ga^n_\nu  }
 g(s,Y^n_s, Z^n_s   )  ds \big| \cF_{\ga_\nu  } \]$ converges to $0$ in probability $\hP$.
 Up to a subsequence of $\{(Y^n,Z^n)\}_{n \in \hN}$, one has
 \beas
 \lmt{n \to \infty } \hE \[  \int_{\ga_\nu  }^{\ga^n_\nu  }
  g(s,Y^n_s, Z^n_s   )  ds \big| \cF_{\ga_\nu  } \] = 0, \q \pas ,
 \eeas
 which together with \eqref{eq:b141}$-$\eqref{eq:b145} leads to that
 $ Y_{\ga_\nu} \n \ge \n  \b1_{\{\ga_\nu = T\}} \xi  \n + \n  \b1_{\{\ga_\nu<T\}}  U_{\ga_\nu  } $,   \pas ~ \qed

\bibliographystyle{siam}
\bibliography{DRBSDE_UI}

\end{document}